\documentclass[12pt]{amsart}

\usepackage{amsmath}
\usepackage{amsmath,amsfonts,amssymb,amsthm}
\usepackage{graphicx,color}

\usepackage{mathrsfs}
\usepackage{hyperref}
\DeclareMathRadical{\sqrtsign}{symbols}{"70}{largesymbols}{"70}
\providecommand{\abs}[1]{\lvert#1\rvert}

\newlength{\figboxwidth}             
\newcommand{\infinity}{\infty}

\def\@ifundefined#1#2#3%
  {\expandafter\ifx\csname#1\endcsname\relax#2\else#3\fi}

\@ifundefined{theoremstyle}{
}{
\theoremstyle{plain} 
}
\newtheorem{theorem}{Theorem}[section]

\newtheorem{proposition}[theorem]{Proposition}
\newtheorem{lemma}[theorem]{Lemma}

\newtheorem{corollary}[theorem]{Corollary}

\@ifundefined{theoremstyle}{
}{
\theoremstyle{definition} 
}
\newtheorem{definition}[theorem]{Definition}

\newtheorem{remark}[theorem]{Remark}

\mathchardef\GG="321D
\newcommand{\mcc}[1]{{}}

\numberwithin{equation}{section}

\begin{document}

\title[Hausdorff dimension of some subsets near $3$]
{Hausdorff dimension of some subsets of the Lagrange and Markov spectra near $3$}

\author{Carlos Gustavo \uppercase{Moreira}}
\address {SUSTech International Center for Mathematics, Shenzhen, Guangdong, People’s Republic of China and IMPA, Estrada Dona Castorina 110, 22460-320, Rio de Janeiro, Brazil}
\email{gugu@impa.br}
\thanks{The first author was partially supported by CNPq and FAPERJ}

\author{Christian Camilo Silva Villamil}
\address {SUSTech International Center for Mathematics, Shenzhen, Guangdong, People’s Republic of China and IMPA, Estrada Dona Castorina 110, 22460-320, Rio de Janeiro, Brazil}
\email{ccsilvav@sustech.edu.cn}
\thanks{The second author was partially supported by CNPq}


\begin{abstract}
We study the sets $\mathcal{L}$ and $\mathcal{M}\setminus\mathcal{L}$ near $3$, where $\mathcal{L}$ and $\mathcal{M}$ are the classical Lagrange and Markov spectra. More specifically, we construct a strictly decreasing sequence $\{a_r\}_{r\in \mathbb{N}}$ converging to $3$, such that for any $r$ one can find a subset $\mathcal{B}_r\subset (a_{r+1},a_r)\cap \mathcal{L}^{'}$ with the property that the Hausdorff dimension of $((a_{r+1},a_r)\cap \mathcal{L})\setminus \mathcal{B}_r$ is less than the Hausdorff dimension of $\mathcal{B}_r$ and for $t\in \mathcal{B}_r$ the sets of irrational numbers with Lagrange value bounded by $t$ and exactly $t$ respectively, have the same Hausdorff dimension. We also show that, as $t$ varies in $\mathcal{B}_r$, this Hausdorff dimension is a strictly increasing function. Finally, in relation to $\mathcal{M}\setminus \mathcal{L}$, we find $C>0$ such that we can bound from above the Hausdorff dimension of $(\mathcal{M}\setminus \mathcal{L})\cap (-\infty,3+\rho)$ by 
$\frac{\log (\abs{\log \rho})-\log (\log(\abs{\log \rho}))+C}{\abs{\log \rho}}$ if $\rho>0$ is small.
 
\end{abstract}

\maketitle

\tableofcontents

\section{Introduction}

	\subsection{The Lagrange spectrum}
	
The Lagrange spectrum is a subset of the real line which appears naturally in the study of Diophantine approximations of real numbers. Consider an irrational real number $x \in \mathbb{R} \setminus \mathbb{Q}$. We define $\ell(x)$ as the supremum of the set of all $k> 0$ such that
$$\left|x - \frac{p}{q}\right| < \frac{1}{k q^2}$$
holds for infinitely many pairs of integers $p, q$ with $q>0$ (possibly with $\ell(x) = \infty$). The number $\ell(x)$ is known as the \emph{Lagrange value} of $x$, and the \emph{Lagrange spectrum} is defined as the set of all finite Lagrange values:
$$\mathcal{L} = \{ \ell(x) < \infty \ \mid\ x \in \mathbb{R} \setminus \mathbb{Q} \}.$$
	
By means of the continued fraction expansion of $x$, it is possible to obtain the symbolic-dynamical characterization of the Lagrange spectrum:
$$	\mathcal{L} = \left\{ \limsup_{n \to \infty} \lambda(\sigma^n(\omega))< \infty\ \mathbin{\Big|}\ \omega \in (\mathbb{N}^*)^{\mathbb{Z}} \right\},$$
where, for $\omega=(\omega_n)_{n\in\mathbb{Z}} \in (\mathbb{N}^*)^{\mathbb{Z}}$, $\lambda(\omega) = [\omega^+] + [0; \omega^-]$, where $\omega^+ = (\omega_n)_{n \geq 0}$ and $\omega^- = (\omega_{-n})_{n \ge 1}$ and $\sigma(\omega)=(\omega_{n+1})_{n\in\mathbb{Z}}$.
	
\subsection{The Markov spectrum}
	
The Markov spectrum is another fractal subset of the real line which is very closely related to the Lagrange spectrum. Using the symbolic-dynamical definition of the Lagrange spectrum as starting point, it can be defined similarly as
	
	$$\mathcal{M} = \left\{ \sup_{n \in\mathbb{Z}} \lambda(\sigma^n(\omega))< \infty\ \mathbin{\Big|}\ \omega \in (\mathbb{N}^*)^{\mathbb{Z}} \right\}.$$	
	We denote by $m(\omega)=\sup_{n \in \mathbb{Z}} \lambda(\sigma^n(\omega))$ the \emph{Markov value} of $\omega \in (\mathbb{N}^*)^{\mathbb{Z}}$.
	
This set is also related to some Diophantine approximation problems. Indeed, it encodes the (inverses of) minimal possible values of real indefinite quadratic forms with normalized discriminants (equal to $1$). Nevertheless, throughout this article we will only use the symbolic-dynamical definitions of $\mathcal{L}$ and $\mathcal{M}$.

We refer the reader to the expository article by Bombieri \cite{Bombi} and to the books by Cusick--Flahive \cite{CF}, and by Lima--Matheus--Moreira--Romaña \cite{LMMR} for a more detailed account on these sets.
		
\subsection{Structure of the Lagrange and Markov spectra}
	
Both the Lagrange and Markov spectra have been intensively studied since the seminal work of Markov \cite{M79}. In particular, it is well-known that
$$ \mathcal{L} \cap [0, 3) = \mathcal{M} \cap [0, 3) = \left\{ \sqrt{5} < \sqrt{8} < \frac{\sqrt{221}}{5} < \dotsb \right\},$$
that is, $\mathcal{L}$ and $\mathcal{M}$ coincide below $3$ and consist of a sequence of \emph{explicit} quadratic surds accumulating only at $3$. Moreover, it is also possible to explicitly characterize the sequences $\omega \in (\mathbb{N}^*)^{\mathbb{Z}}$ associated with Markov values less than or equal to $3$, \cite{M79} and \cite{Bombi}.

On the other hand, the behavior of these sets after $3$ remains somewhat mysterious. Indeed, it is known that $\mathcal{L} \subseteq \mathcal{M}$ and some authors conjectured that these sets are \emph{equal}; Freiman disproved this conjecture only in 1968 \cite{F}. But now, much more is known in this regard: In \cite{MM} it was proved that the Hausdorff dimension $HD(\mathcal{M}\setminus \mathcal{L})$ of the set difference between $\mathcal{M}$ and $\mathcal{L}$ is larger than $1/2$ but strictly smaller than $1$, and thus $int(\mathcal{M})=int(\mathcal{L)}$, i.e., the interior of the Markov spectrum coincides with the interior of the Lagrange spectrum.
	
Even if the previous paragraph suggests that these sets are somewhat different, they are known to coincide before $3$ and after large enough values. Indeed, Hall showed in 1947 that $\mathcal{L}$ and $\mathcal{M}$ contain a half-line $[c, \infty)$; any such ray is hence known as a \emph{Hall ray}. After several years, Freiman found the largest Hall ray to be $[c_{\text{F}}, \infty)$, where $c_{\text{F}} \approx 4.5278\ldots$ is an explicit quadratic surd known as Freiman's constant \cite{F2}. These results in turn imply that $\mathcal{L}$ and $\mathcal{M}$ coincide starting at $c_{\text{F}}$, so they both contain the half-line $[c_{\text{F}}, \infty)$.
	
There are more striking similarities between these two sets. In particular, their Hausdorff dimensions coincide when truncated: the first author showed that	
$$HD(\mathcal{L} \cap (-\infty, t)) = HD(\mathcal{M} \cap (-\infty, t))$$	
for every $t \in\mathbb{R}$, where $HD(X)$ denotes the Hausdorff dimension of the set $X$ \cite{M3}. Clearly, this result shows that, when studying the Hausdorff dimension of such truncated versions, one can choose to use either $\mathcal{L}$ or $\mathcal{M}$. Define the function $d:\mathbb{R} \rightarrow [0,1]$ given by
 $$d(t)= H(\mathcal{L} \cap (-\infty, t)) = H(\mathcal{M} \cap (-\infty, t)).$$	
Moreira also proved in \cite{M3} that $d$ is continuous, surjective and such that $d(3)=0$. Moreover, that $d(t)=\min \{1,2D(t)\},$ where 
$$D(t)=HD(\ell^{-1}(-\infty,t))=HD(\ell^{-1}(-\infty,t])$$ 
is also a continuous surjective function from $\mathbb{R}$ to $[0,1).$ 

Recently in \cite{ERMR}, more precise estimates of $d(t)$ were given for $t$ close to $3$. Specifically, if $H\colon[-1, +\infty)\to [-e^{-1},+\infty)$ is given by $H(x)=x e^x$ (its inverse, $H^{-1}\colon [-e^{-1},+\infty)\to [-1,+\infty)$ is the Lambert function), then for all sufficiently small $\rho>0$, we have
$$d(3+\rho)=2\cdot\frac{H^{-1}(e^{c_0}|\log \rho|)}{|\log \rho|}+O\left(\frac{\log (|\log \rho|)}{|\log \rho|^2}\right),$$
where $c_0=-\log\log((3+\sqrt{5})/2)$. In particular, we get for some constants $C_1, C_2>0$ and small $\rho>0$ the following inequalities that we will use in the proof of our main theorems:
\begin{equation}\label{cunocdos}
C_1\cdot \frac{\log (|\log \rho|)}{|\log \rho|} \leq d(3 + \rho) \leq C_2 \cdot \frac{\log (|\log \rho|)}{|\log \rho|}
\end{equation}
and 
\begin{equation}\label{cunocdos2}
d(3+\rho)\leq 2\cdot \frac{\log(\abs{\log \rho})-\log (\log (\abs{\log \rho}))+C_2}{\abs{\log \rho}}.
\end{equation}

In this article, we are interested in values of $t\in \mathbb{R}$ such that $t<t_1$ where $t_1:=\sup \{s\in \mathbb{R}:d(s)<1\}=3.334384...$ (see \cite{MMPV}).
\subsection{Dynamical spectra}
Let $\varphi:S\rightarrow S$ be a diffeomorphism of a $C^{\infty}$ compact surface $S$ with a mixing horseshoe $\Lambda$ and let $f:S\rightarrow \mathbb{R}$ be a differentiable function. Following the above characterization of the classical spectra, we define the maps $\ell_{\varphi,f}: \Lambda \rightarrow \mathbb{R}$ and $m_{\varphi,f}: \Lambda \rightarrow \mathbb{R}$ given by $\ell_{\varphi,f}(x)=\limsup\limits_{n\to \infty}f(\varphi^n(x))$ and $m_{\varphi,f}(x)=\sup\limits_{n\in\mathbb{Z}}f(\varphi^n(x))$ for $x\in \Lambda$ and call $\ell_{\varphi,f}(x)$ the \textit{Lagrange value} of $x$ associated to $f$ and $\varphi$ and also $m_{\varphi,f}(x)$ the \textit{Markov value} of $x$ associated to $f$ and $\varphi$. The sets
$$\mathcal{L}_{\varphi,\Lambda,f}=\ell_{\varphi,f}(\Lambda)=\{\ell_{\varphi,f}(x):x\in \Lambda\}\ \  
\text{and}\ \
\mathcal{M}_{\varphi,\Lambda, f}=m_{\varphi,f}(\Lambda)=\{m_{\varphi,f}(x):x\in \Lambda\}$$
are called \textit{Lagrange Spectrum} of $(\varphi,\Lambda,f)$ and \textit{Markov Spectrum} of $(\varphi,\Lambda,f)$ respectively.

Let us first fix a Markov partition $\{R_a\}_{a\in \mathcal{A}}$ consisting of rectangles $R_a$ with small diameter delimited by compact pieces $I_a^s$, $I_a^u$, of stable and unstable manifolds of certain points of $\Lambda$.
As usual, if $\Sigma=\left\{\underline{a}=(a_{n})_{n\in \mathbb{Z}}\in \mathcal{A}^{\mathbb{Z}}:\forall n\in \mathbb{Z},\ \ \varphi(R_{a_n})\cap R_{a_{n+1}}\neq \emptyset\right\}$ is equipped with the shift $\sigma:\Sigma\to\Sigma$ defined by $\sigma(\underline{a})_{n}=a_{n+1}$. The dynamics of $\varphi$ on $\Lambda$ is topologically conjugate to the shift on $\Sigma$, namely, there is a homeomorphism $\Pi: \Lambda \to \Sigma$ such that $\varphi\circ \Pi=\Pi\circ \sigma$.

Using the locally invariant $C^{1+\alpha}$ stable and unstable foliations (where $\alpha>0$), it is possible to define projections $\pi^u_a: R_a\rightarrow I^s_a$ and $\pi^s_a: R_a\rightarrow I^u_a$. Given $x\in R_a$, set $\pi^s(x)=\pi^u_a(x)$ and $\pi^u(x)=\pi^s_a(x)$. In this way, we have the stable and unstable Cantor sets
$$K^s=\pi^s(\Lambda)=\bigcup_{a\in \mathcal{A}}\pi^u_a(\Lambda\cap R_a)\ \  
\text{and}\ \
K^u=\pi^u(\Lambda)=\bigcup_{a\in \mathcal{A}}\pi^s_a(\Lambda\cap R_a),$$
which are $C^{1+\alpha}$ dynamically defined, associated to the expanding maps $\psi_s$ and $\psi_u$ defined by
$$\psi_s(\pi^s(y))=\pi^s(\varphi^{-1}(y))\ \  
\text{and}\ \
\psi_u(\pi^u(z))=\pi^u(\varphi(z)).$$

It turns out that dynamical Markov and Lagrange spectra associated to hyperbolic dynamics are closely related to the classical Markov and Lagrange spectra. Several results on the Markov and Lagrange dynamical spectra associated to horseshoes in dimension 2 which are analogous to previously known results on the classical spectra were obtained recently. We refer the reader to the book \cite{LMMR} for more information.

In our present work, it is important to mention that in \cite{GCD}, in the context of {\it conservative} diffeomorphism it is proven (as a generalization of the results in \cite{CMM16}) that for typical choices of the dynamic and of the function, the intersections of the corresponding dynamical Markov and Lagrange spectra with half-lines $(-\infty,t)$ have the same Hausdorff dimensions, and this defines a continuous function of $t$ whose image is $[0,\min \{1,\tau\}]$, where $\tau$ is the Hausdorff dimension of the horseshoe. 

Finally, in \cite{LM2} is showed that, for any $N\ge 2$ with $N\neq 3$, the initial segments of the classical spectra until $\sqrt{N^2+4N}$ (i.e., the intersection of the spectra with $(-\infty,\sqrt{N^2+4N}]$) are dynamical Markov and Lagrange spectra associated to a horseshoe $\Lambda(N)$ (naturally associated to continued fractions with coefficients bounded by $N$) of some smooth conservative diffeomorphism $\varphi_N$ of $\mathcal{S}^2$ and to some smooth real function $f_N$. Using this, in \cite{GC} it is proven that for any $t$ that belongs to the closure of the interior of the classical Markov and Lagrange spectra $D(t)=HD(\ell^{-1}(t))$ and that $D$ is strictly increasing when is restricted to the interior of the spectra. One related result is in \cite{C2}, where it was shown that unless some countable set, the set $\mathcal{J}$ of elements $t\in\mathcal{L}$ that satisfies that $D(t)=HD(\ell^{-1}(t))$, is the same as the set of $t\in\mathcal{L}$ where $d(t)=\lim \limits_{\epsilon\to 0^{+}}HD(\mathcal{L}\cap (t-\epsilon,t+\epsilon))$ i.e., the set of $t\in\mathcal{L}$ where $d(t)$ is equal to the local Hausdorff dimension of the Lagrange spectrum at $t$. Even more, it was proved that one can set $\mathcal{J}=\{\eta^-:\eta\in (0,1)\}$ where for $\eta\in (0,1)$, $\eta^-=\min\{t\in \mathbb{R}:D(t)=\eta\}$.

Here, we will explore again the dynamical nature of the classical spectra (at least the portion until $\sqrt{12}$) to study the sets $\mathcal{L}$ and $\mathcal{M}\setminus\mathcal{L}$ near $3$. Now, we can state our main results.

Our first theorem decomposes one interval of the form $(3,s)$ where $s<t_1$ into an enumerable collection of disjoint intervals in such a way that for each of these intervals one can find some subset of $\mathcal{L}$ contained on it with ``big'' Hausdorff dimension and such that $D$ restricted to this set has similar properties as $D$ has in the interior of the spectra. To be more precise, one has
\begin{theorem}\label{teo1}
There exists a decreasing sequence $\{a_r\}_{r\in \mathbb{N}}$ with $a_1<t_1$, $d(a_{r+1})<d(a_r)$ and $\lim \limits_{n\rightarrow \infty}a_r=3$, such that, given $r\in \mathbb{N}$ we can find a subset $\mathcal{B}_r\subset (a_{r+1},a_r)\cap \mathcal{L}$ with the following properties:
    \begin{itemize}
        \item $HD(((a_{r+1},a_r)\cap \mathcal{L})\setminus \mathcal{B}_r)<d(a_{r+1})$,
        \item $D(t)=HD(\ell^{-1}(t))$ for $t \in \mathcal{B}_r$, 
        \item $HD(\mathcal{B}_r)=HD((a_{r+1},a_r)\cap \mathcal{L})$,
        \item $D|_{\mathcal{B}_r}$ is strictly increasing,
        \item $\mathcal{B}_r\subset \mathcal{L}^{'}$.  \end{itemize}   
\end{theorem}

Our second theorem establishes a non-trivial upper bound for the Hausdorff dimension of the different set 
$\mathcal{M}\setminus \mathcal{L}$ close to $3$.

\begin{theorem}\label{teo2}
There is a constant $C>0$ such that 
$$HD((\mathcal{M}\setminus \mathcal{L})\cap (-\infty,3+\rho))\leq \frac{\log (\abs{\log \rho})-\log (\log (\abs{\log \rho}))+C}{\abs{\log \rho}}$$
for any $\rho>0$ small.
\end{theorem}

\section{Preliminares}

\subsection{Sets of finite type and connection of subhorseshoes}\label{tipofinito}

The following definitions and results can be found in \cite{GC}. Fix a horseshoe $\Lambda$ of some conservative diffeomorphism $\varphi:S\rightarrow S$ and $\mathcal{P}=\{R_a\}_{a\in \mathcal{A}}$ some Markov partition for $\Lambda$. Take a finite collection $X$ of finite admissible words $\theta=(a_{-n(\theta)},\dots,a_{-1},a_0,a_1,\dots,a_{n(\theta)})$ and set $R(\theta;0)=\Pi^{-1}\{(x_{n})\in \Sigma:(x_{-n(\theta)},\dots,x_0,\dots,x_{n(\theta)})=\theta    \}$. We say that the maximal invariant set 
$$M(X)=\bigcap \limits_{m \in \mathbb{Z}} \varphi ^{-m}(\bigcup \limits_{\theta \in X}  R(\theta;0))$$ 
is a \textit{hyperbolic set of finite type}. Even more, it is said to be a \textit{subhorseshoe} of $\Lambda$ if it is nonempty and $\varphi|_{M(X)}$ is transitive. Observe that a subhorseshoe need not be a horseshoe; indeed, it could be a periodic orbit in which case it will be called of trivial.

By definition, hyperbolic sets of finite type have local product structure. In fact, any hyperbolic set of finite type is a locally maximal invariant set of a neighborhood of a finite number of elements of some Markov partition of $\Lambda$.

\begin{definition}
Any $\tau \subset M(X)$ for which there are two different subhorseshoes $\Lambda^1$ and $\Lambda^2$ of $\Lambda$ contained in $M(X)$ with 
$$\mathcal{T}_{\Lambda^1, \Lambda^2}=\{x\in M(X):\ \omega(x)\subset \Lambda^1\ \text{and}\ \alpha(x)\subset \Lambda^2  \}$$
will be called a transient set or transient component of $M(X)$.
\end{definition}
Note that for any subhorseshoe $\tilde{\Lambda}\subset \Lambda$, being  $\varphi$ conservative, one has
\begin{equation}\label{transient2}
HD(\tilde{\Lambda})=HD(K^s(\tilde{\Lambda}))+HD(K^u(\tilde{\Lambda}))=2HD(K^u(\tilde{\Lambda})). 
\end{equation}
And also, by the local product structure and the previous equation, given a transient set $\tau$ as before, it is true that 
\begin{equation}\label{transient}
 HD(\mathcal{T}_{\Lambda^1, \Lambda^2})=HD(K^s(\Lambda^2))+HD(K^u(\Lambda^1))=\frac{HD(\Lambda^2)+HD(\Lambda^1)}{2}.
\end{equation}

\begin{proposition}\label{appendix}
Any hyperbolic set of finite type $M(X)$, associated with a finite collection of finite admissible words $X$ as before, can be written as
$$M(X)=\bigcup \limits_{i\in \mathcal{I}} \tilde{\Lambda}_i $$ 
where $\mathcal{I}$ is a finite index set (that may be empty) and for $i\in \mathcal{I}$,\ $\tilde{\Lambda}_i$ is a subhorseshoe or a transient set.
\end{proposition}

Fix $f:S\rightarrow \mathbb{R}$  differentiable. Given $t\in \mathbb{R}$, we define 
$$\Lambda_t=m_{\varphi, f}^{-1}((-\infinity,t])= \{x\in\Lambda: \forall n\in \mathbb{Z}, \ f(\varphi^n(x))\leq t\}.$$
A notion that plays an important role in the proof of Theorems $1$ and $2$ is the notion of \textit{connection of subhorseshoes}. 

\begin{definition}\label{conection of horseshoes1}
Given $\Lambda^1$ and $\Lambda^2$ subhorseshoes of $\Lambda$ and $t\in \mathbb{R}$, we said that $\Lambda^1$ \emph{connects} with $\Lambda^2$ or that $\Lambda^1$ and $\Lambda^2$ \emph{connect} before $t$ if there exist a subhorseshoe $\tilde{\Lambda}\subset \Lambda$ such that $\Lambda^1 \cup \Lambda^2 \subset \tilde{\Lambda}$ and $\max f|_{\tilde{\Lambda}}<t$.
\end{definition}

For our present purposes, the next criterion of connection will be important.

\begin{proposition}\label{connection11}
Suppose $\Lambda^1$ and $\Lambda^1$ are subhorseshoes of $\Lambda$ and for some $x,y \in \Lambda$ we have $x\in W^u(\Lambda^1)\cap W^s(\Lambda^2)$ and $y\in W^u(\Lambda^2)\cap W^s(\Lambda^1)$. If for some $t\in \mathbb{R}$, it is true that 
$$\Lambda^1 \cup \Lambda^2 \cup \mathcal{O}(x) \cup \mathcal{O}(y) \subset \Lambda_t,$$ then for every $\epsilon >0$,\ $\Lambda^1$ and $\Lambda^2$ connect before $t+\epsilon$. In addition, the subhorseshoe $\tilde{\Lambda}(\epsilon)$ in Definition \ref{conection of horseshoes1} can be chosen such that $\mathcal{O}(x) \cup \mathcal{O}(y)\subset \tilde{\Lambda}(\epsilon)$.
\end{proposition}

 \begin{corollary}\label{connection3}
 Let $\Lambda^1$,\ $\Lambda^2$ and $\Lambda^3$ subhorseshoes of $\Lambda$ and $t\in \mathbb{R}$. If $\Lambda^1$\ connects with $\Lambda^2$ before $t$ and $\Lambda^2$\ connects with $\Lambda^3$ before $t$. Then also $\Lambda^1$\ connects with $\Lambda^3$ before $t$.
 \end{corollary}

\subsection{The horseshoe $\Lambda(2)$}\label{Lambda}

Let $N\geq 2$ and $N\neq 3$ be an integer. In \cite{LM2} is proved that the portion of the classical spectra up to $\sqrt{N^2+4N}$ i.e, $\mathcal{L}\cap (-\infty, \sqrt{N^2+4N}]$ and $\mathcal{M}\cap (-\infty, \sqrt{N^2+4N}]$ are the dynamically defined Lagrange and Markov spectra  $\mathcal{L}_{\varphi,\Lambda(N),f}$ and $\mathcal{M}_{\varphi,\Lambda(N),f}$ associated with some $\varphi$, $\Lambda(N)$ and $f$.
More specifically, if $C_N=\{x=[0;a_1,a_2,...]: a_i\le N, \forall i\ge 1\}$ and $\tilde{C}_N=\{1,2,...,N\}+C_N$, we set $\Lambda(N)=C_N\times \tilde{C}_N$ and
then consider $\varphi:\Lambda(N) \rightarrow \Lambda(N)$ given by
$$\varphi([0;a_1,a_2,...],[a_0;a_{-1},a_{-2},...])=([0;a_2,a_3,...],[a_1;a_0,a_{-1},...]),$$
that can be extended to a $C^{\infty}$ conservative diffeomorphism on a diffeomorphic copy of the 2-dimensional sphere $\mathbb{S}^2$. Also, the real map is given by $f(x,y)=x+y$.

For $\Lambda(N)$ we have the Markov partition $\{R_a\}_{a\in \mathcal{A}}$ where $\mathcal{A}=\{1,2, \dots,N\}$ and $R_a$ is such that $R_a\cap \Lambda(N)=C_N\times(C_N+a)=C_N\times C_N+(0,a)$. By definition, $\varphi$ expands in the $x$-direction and contracts in the $y$-direction. Therefore, for $([0;a_1,a_2,...],[a_0;a_{-1},a_{-2},...])\in R_{a_0}$ we can set $\pi^s_{a_0}([0;a_1,a_2,...],[a_0;a_{-1},a_{-2},...])=([0;a_1,a_2,...],[a_0;\overline{N,1}])$ and then  
\begin{eqnarray*}
\psi_u([0;a_1,a_2,...],[a_0;\overline{N,1}])&=&\pi ^s_{a_1}(\phi([0;a_1,a_2,...],[a_0;\overline{N,1}]))\\
&=& \pi ^s_{a_1}([0;a_2,a_3,...],[a_1;a_0,\overline{N,1}])\\
&=&([0;a_2,a_3,...],[a_1;\overline{N,1}]).
	\end{eqnarray*}
Thus we can identify $(K^u(\Lambda(N)),\psi_u)$ with $(\tilde{C}_N,\tilde{G})$ where for $[a_0;a_1,a_2,...]\in \tilde{C}_N$ one has $\tilde{G}([a_0;a_1,a_2,...])=[a_1;a_2,a_3,...]$. A similar identification can be made for $(K^s(\Lambda(N)),\psi_s)$.

As in the introduction, if $\Sigma_N=\{1,2,\dots,N \}^{\mathbb{Z}}$, one has that $\varphi|_{\Lambda(N)}$ is topologically conjugated to $\sigma:\Sigma_N\rightarrow \Sigma_N$ (via the function $\Pi: \Lambda(N) \to \Sigma_N$ given by $\Pi([0;a_1,a_2,...],[a_0;a_{-1},...])=(...,a_{-2},a_{-1},a_0,a_1,a_2,...)$ ), and that in sequences, $f$ becomes the restriction of $\lambda$ to $\Sigma_N$. 

Our goal is to study the structure of the sets $\mathcal{L}$ and $\mathcal{M}\setminus \mathcal{L}$ near $t = 3$. Note that, if a sequence $\omega\in(\mathbb{N}^*)^{\mathbb{Z}}$ contains any letter greater or equal than $3$, then $\lambda(\omega)>3.52$, which is ``much larger'' than $t_1$, so we can ignore such sequences. Thus, throughout the entire article, a \emph{word} is made up of letters in the alphabet $\{1, 2\}$. That is, we can restrict our attention to the triple $(\varphi,\Lambda(2),f)$ or, in sequences, to the triple $(\sigma,\Sigma_2,\lambda)$. As we will use both points of view, it is convenient to define for $t\in \mathbb{R}$ the set $\Sigma(t)=\Pi((\Lambda(2))_t)=\{\omega\in\Sigma_2: \forall n\in \mathbb{Z}, \ \lambda(\sigma^n(\omega))\leq t\}.$ 

\subsection{Cuts, alphabets and renormalization of words}
In this article, we will use freely some notations, definitions and theorems of the first sections of \cite{ERMR}. Through this subsection we will introduce them. Hence, most of the content will not be neither proved nor referenced. 
\subsubsection{Sequences in $\Sigma(3)$}  
Bombieri in \cite{Bombi} showed that bi-infinite words $\omega\in \Sigma(3)$ have to follow very special patterns (which is essentially a restatement of much older results by Markov \cite{M79}). Indeed, he showed that $\omega$ must be a word in the letters $a = 22$ and $b = 11$, that is, the number of consecutive ones or twos is always even or infinite. And if $U$ and $V$ are the Nielsen substitutions given by 
    $$  U \colon \begin{matrix}
		a &\mapsto &ab \\
		b &\mapsto &b
	\end{matrix}, \qquad
	V \colon \begin{matrix}
		a &\mapsto &a \\
		b &\mapsto &ab.
	\end{matrix}
    $$
and extended to finite or infinite words in the alphabet $\{a,b \}$ in the obvious way, then the words $\omega$ with Markov value less than $3$ are exactly the periodic sequences with period either $a$ or $b$ or $W(ab)$ for some unique word $W$ in the alphabet $\{U,V\}$. 
    
Given a pair of words $(u, v)$ in the alphabet $\{a,b\}$, we also define the operations $\overline{U}(u, v) = (uv, v)$ and $\overline{V}(u, v) = (u, uv)$. Let $T$ be the tree obtained by successive applications of the functions $\overline{U}$ and $\overline{V}$, starting at the root $(a, b)$, let $\mathscr{A}$ be the set of vertices of $T$ and let $\mathscr{A}_n$, for $n \geq 0$, be the set of elements of $\mathscr{A}$ whose distance to the root $(a, b)$ is exactly $n$. Also, define $\mathsf{c}$ as the concatenation operator, that is, $\mathsf{c}(u, v) = uv$. The following lemma gives us an alternative description of the periods of words with Markov value less than $3$. 
	
\begin{lemma}\label{lem:alphabets}
		Let $(\alpha, \beta) \in \mathscr{A}$. Then, there exists a word $W$ in the alphabet $\{U,V\}$ such that $\alpha = W(a)$ and $\beta = W(b)$. In particular, $\mathsf{c}(\mathscr{A})$ is the set of periods (different from $a$ and $b$) of words with Markov value less than $3$.
	\end{lemma}   
We say that $(\alpha,\beta)\in\mathscr{A}$ is an \emph{ordered alphabet} and that $\{\alpha,\beta \}$ is the alphabet associated to $(\alpha,\beta)$.  
We finish this subsubsection with the following technical lemma. 
\begin{lemma}\label{lem:determine_alphabet}
Suppose that a word $w$ can be written as a concatenation $\tau\alpha\beta\tau'$ for some words $\tau$, $\tau'$, $\alpha$ and $\beta$, with $(\alpha,\beta)\in \mathscr{A}_n$ and $n \in \mathbb{N}$. If there exist $(A,B)\in \mathscr{A}_n$, $k\ge 1$ and $w_1,\dotsc,w_k\in \{A,B\}$ such that $w=w_1\ldots w_k$, then $(A,B)=(\alpha,\beta)$ and there exists $1\le j<k$ such that $w_1\ldots w_{j-1}=\tau$, $w_j=\alpha$, $w_{j+1}=\beta$ and $w_{j+2}\ldots w_k=\tau'$.  
\end{lemma}    

\subsubsection{Cuts and intervals determined by words}  Through the text, words can be either finite or bi-infinite. If $w$ is a finite word, we denote its \emph{size} by $\abs{w}$, that is, the number of letters $1$ or $2$ that are needed to write $w$. If $w=a_1\dots a_n\in \{1,2 \}^n$ we denote by $w^*=a_n\dots a_1$ its \emph{transposed}. We will also consider \emph{cuts} of finite words, which consist of a word $w$ together with a choice of a pair of letters marked with a vertical segment. We usually write cuts as $\omega = w_1|w_2$, where $w_1$ and $w_2$ are finite non-empty words. We said that a cut in $w$ is a \emph{good cut} if for any bi-infinite word containing $w$ as a factor, $\lambda$ has value less than $3$ in any of the two positions determined by the cut. Similarly, we said that a cut in $w$ is a \emph{bad cut} if for any bi-infinite word containing $w$ as a factor, $\lambda$ has value greater than $3$ in some of the positions determined by the cut.

Now, recall that $\Sigma(t) = \{\omega \in (\mathbb{N}^*)^{\mathbb{Z}} \ \mid\ m(\omega) \leq t\}$ and given $n\in\mathbb{N}$ define $\Sigma(t, n)$ as the set of length-$n$ subwords of sequences in $\Sigma(t)$. The following theorem is important for us.	
\begin{theorem} \label{palabras}
For $n \geq 68$, one has		
$$\Sigma(3 + 6^{-3n}, n) = \Sigma(3, n) = \Sigma(3 - 6^{-3n}, n).$$
\end{theorem}	
This theorem can be interpreted as follows: given a bi-infinite word, whose Markov value is exponentially close to $3$ (smaller than $3+6^{-3n}$), then its size-$n$ subwords are indistinguishable from those in $\Sigma(3, n)$. That is to say, a length-$n$ word cannot detect the patterns of symbols that make their Markov values different from $3$; they are only present when considering words of larger lengths.

We will prove several lemmas that allow us to understand the structure of the bi-infinite words in $\Sigma(3+6^{-3n})$ and their finite subwords. But for this, it is necessary to show first some basic facts about the function $\lambda$, in terms of the intervals determined by finite words.

For a finite word $w$ consider the closed interval $I(w)$ consisting of the numbers in $[0,1]$ whose continued fractions start with $w$. We define $r(\omega)=\lfloor\log (\abs{I(w)}^{-1}) \rfloor$, which controls the order of magnitude of the size of $I(w)$. Let us recall the following properties that will be useful for us. 

\begin{itemize}
    \item If $w$ is a non-empty finite word in the alphabet $\{1,2\}$ then    
    $$  (\abs{w}-3) \log\left( \frac{3 + \sqrt{5}}{2} \right) \leq r(w) \leq (\abs{w}+1) \log(3 + 2 \sqrt{2}).$$
    \item For $n\in \mathbb{N}^{*}$ one has
    $$ \abs{I(1^n)}^{-1}= -\frac{1}{5}(-1)^{n+1} + \frac{\sqrt{5}+1}{10}\left( \frac{3 + \sqrt{5}}{2} \right)^{n+1} - \frac{\sqrt{5}-1}{10}\left(\frac{3 - \sqrt{5}}{2}\right)^{n+1} $$ and       
$$ \abs{I(2^n)}^{-1} = \frac{(3 + 2\sqrt{2})^n - (3 - 2\sqrt{2})^n}{4\sqrt{2}}.$$
\item For any finite words $w_1$ and $w_2$ one has 
    $$\frac{1}{2} \abs{I(w_1)}\cdot\abs{I(w_2)}<\abs{I(w_1w_2)}<2\abs{I(w_1)}\cdot\abs{I(w_2)}.$$
\end{itemize}
We have the following lemmas.
\begin{lemma}\label{lem:121_212_forbidden}
Let $\omega \in \Sigma(3.06)$. Then, $\omega$ does not contain $121$ or $212$ as subwords.
\end{lemma}

\begin{lemma} \label{intervals}
Let $\omega\in \Sigma_2$ not containing $121$ and $212$ as subwords and such that $\omega = \dots r_2r_1w^*b|aws_1s_2\dots$, where $w$ is a finite word, $r_1 \neq s_1$, with $r_i, s_i \in \{1,2\}$ for each $i$. Then
$$\abs{I(bwb)}< sign([w,s_1,s_2,\dots]-[w,r_1,r_2\dots])(\lambda(\omega)-3)<\abs{I(bw1)}.$$
In particular if $w$ has even length, $r_1=1$ and $s_1=2$, then 
$$ \abs{I(bwb)}< \lambda(\omega)-3<\abs{I(bw1)}.$$	
\end{lemma}

\begin{corollary}\label{simetria}
Let $w$ a finite word in the alphabet $\{a,b \}$. For any bi-infinite word $\omega=\dots bw^*b|aw a \dots $ one has
$$3+\frac{1}{144}\abs{I(w)}<\lambda(\omega)<3+\frac{1}{3}\abs{I(w)}.$$
\end{corollary}

\begin{proof}
   Using the previous equations, one gets  
    $\abs{I(b)}^{-1}=6$ and $\abs{I(1)}^{-1}=2.$ Then,
   $$\frac{1}{144}\abs{I(w)}=\frac{1}{4}\cdot\frac{1}{36}\abs{I(w)}=\frac{1}{4}\abs{I(b)}^2\cdot\abs{I(w)}\leq \abs{I(bwb)}$$
and 
   $$\abs{I(bw1)}\leq 4\abs{I(b)}\cdot\abs{I(1)}\cdot \abs{I(w)}=\frac{1}{3}\abs{I(w)}.$$
The result follows from Lemma \ref{intervals}.
\end{proof}

\begin{corollary}\label{mas5}
In the conditions of Corollary \ref{simetria} 
$$\lambda(\omega)>3+6^{-(\abs{w}+5)}.$$
In particular, if $\omega\in \Sigma(3+6^{-3n})$ then $\abs{w}+5\geq 3n.$

\end{corollary}\label{8}
\begin{proof}
  As $r(w)<\log (3+2\sqrt{2})^{\abs{w}+1}$, by definition, we have
  $$\log \abs{I(w)}^{-1}<\log (3+2\sqrt{2})^{\abs{w}+1}+1<\log e\cdot 6^{\abs{w}+1}$$
  and then $I(w)>\frac{1}{18}6^{-\abs{w}}$. By Lemma \ref{simetria}
  $$\lambda(\omega)>3+\frac{1}{144\times 18}6^{-\abs{w}}>3+6^{-(\abs{w}+5)}.$$ 
  Additionally, if $\abs{w}+5<3n$ then
  $$\lambda(\omega)>3+6^{-(\abs{w}+5)}>3+6^{-3n}.$$
  
 \end{proof}

The following lemma allows us to compare bad cuts: if one bad cut gives value less than $t$, then other bad cut of the same kind, that is ``worse'' than the first one, also gives value less than $t$.
\begin{lemma}\label{compare}
Let $\omega$ and $\tilde{\omega}$ be two different words in the alphabet $\{1,2\}$ such that $\tilde{\omega}$ begins with $\omega$. If for some $x,y\in \{1,2\}$ with $x\neq y$ one has the bad cut $ x\omega^*b|a\omega y$ in some word of $\Sigma(t)$, where $t>3$. Then for any sequence $X=\dots x\tilde{\omega}^*b|a\tilde{\omega} y \dots$ one has $\lambda(X)<t$, where the zero position of the sequence $X$ is in the $2$ of the cut.

\end{lemma}

\begin{proof}
If $\tilde{\omega}=\omega \tilde{X}$, we have:
\begin{eqnarray*}
 \lambda(X) &=& \lambda(\dots x\tilde{\omega}^*b|a\tilde{\omega} y \dots)=[22\tilde{\omega} y\dots ]+[0;11\tilde{\omega}x\dots ]\\ &=&3-[22\tilde{\omega}y\dots]-[0;11\tilde{\omega}y\dots ]+[22\tilde{\omega}y\dots ]+[0;11\tilde{\omega}x\dots ]\\ &=& 3+ [0;11\tilde{\omega}x\dots]-[0;11\tilde{\omega}y\dots ]= 3+ [0;11\omega\tilde{X}x\dots]-[0;11\omega\tilde{X}y\dots ] \\ &\leq& 
 3+ \abs{[0;11\omega\tilde{X}x\dots]-[0;11\omega\tilde{X}y\dots ]}< 3+ \abs{[0;11\omega x\dots]-[0;11\omega y\dots ]}\\ &=& 3+ [0;11\omega x\dots]-[0;11\omega y\dots ]=\lambda(\dots x\omega^*b|a\omega y \dots)\leq t
 \end{eqnarray*} 
where we use that $\lambda(\dots x\omega^*b|a\omega y \dots)>3$.
\end{proof}




\subsubsection{Weakly renormalizable words}

The following lemma shows that bi-infinite words with Markov value exponentially close to $3$ (relative to the size of the interval they induce) cannot contain both $\alpha\alpha$ and $\beta\beta$ if $(\alpha, \beta) \in \mathscr{A}$. 
	
\begin{lemma}\label{lem:aa_bbr}
Let $(\alpha,\beta)\in \mathscr{A}$. If $w$ is a finite word in the associated alphabet $\{\alpha,\beta\}$ starting with $\alpha\alpha$ and ending by $\beta\beta$ such that $r(w)\leq r$, then the Markov value of any bi-infinite word containing $w$ as a factor is larger than $3+e^{-r}$.
\end{lemma}

We now define the notion of a weakly renormalizable word, which is central to our methods as it is used to find suitable alphabets in which words can be written.
	
\begin{definition}\label{def:weaklyrenormalizable}
Let $(\alpha,\beta)\in \mathscr{A}$ and $w$ be a finite word in the alphabet $\{a,b\}$. We say that $w$ is \emph{$(\alpha,\beta)$-weakly renormalizable} if we can write $w=w_1 \gamma w_2$ where $\gamma$ is a word (called the \emph{renormalization kernel}) in the alphabet $\{\alpha,\beta\}$ and $w_1, w_2$ are (possibly empty) finite words with $\abs{w_1}, \abs{w_2}<\max\{\abs{\alpha},\abs{\beta}\}$ such that $w_2$ is a prefix of $\alpha\beta$ and $w_1$ is a suffix of $\alpha\beta$, with the following restrictions: If $(\alpha,\beta)=(u,uv)$ for some $(u,v)\in\mathscr{A}$ and $\gamma$ ends with $\alpha$, then $\abs{v}\leq \abs{w_2}$. If $(\alpha,\beta)=(uv,v)$ for some $(u,v)\in\mathscr{A}$ and $\gamma$ starts with $\beta$, then $\abs{u}\leq \abs{w_1}$.
\end{definition}
    
Definition \ref{def:weaklyrenormalizable} is motivated by the following ideas. Given an alphabet $\{\alpha, \beta\}$ with $(\alpha, \beta) \in \mathscr{A}$, it may not be possible to write a word $w$ in terms of $\alpha$ and $\beta$. Nevertheless, it may very well be possible to write ``most'' of $w$ in terms of $\alpha$ and $\beta$, preceded by and followed by some short trailing words. These words are $w_1$ and $w_2$ in the previous definition, and the condition ensuring that they are short is that $\abs{w_1}, \abs{w_2} < \max\{\abs{\alpha}, \abs{\beta}\}$. To further ensure that $w_1$ and $w_2$ are well-adjusted to the chosen alphabet, we also require them to be a prefix or suffix of $\alpha\beta$; then $w$ is contained in $\alpha\beta \gamma \alpha\beta$, where the renormalization kernel $\gamma$ can be written in the alphabet $\{\alpha, \beta\}$. 

Exhibiting a word as being $(\alpha, \beta)$-weakly renormalizable is nontrivial in general and, to complicate matters even further, the choice of alphabet $(\alpha, \beta) \in \mathscr{A}$ is not clear to begin with. Nevertheless, any word in the alphabet $\{a, b\}$ is trivially $(a, b)$-weakly renormalizable (by setting the renormalization kernel equal to the entire word). With these considerations, we will now present the renormalization algorithm.

    \begin{lemma}[Renormalization algorithm]\label{lem:renormalizingr}
        Let $w \in \Sigma(3 + e^{-r}, \abs{w})$ satisfying $r(w)\leq r$. If for some $(u,v)\in \mathscr{A}$, $w$ is $(u, v)$-weakly renormalizable as $w = w_1 \gamma w_2$ with $\gamma \neq \emptyset$, then $w$ is $(\alpha, \beta)$-weakly renormalizable for some $(\alpha, \beta) \in \{(uv, v), (u, uv)\}$. Moreover, if $\gamma$ starts with $u$ or ends with $v$, then $w_1$ or $w_2$, respectively, does not change for the renormalization with alphabet $(\alpha,\beta)$.
    \end{lemma}
    
   Let us explain the renormalization algorithm: Suppose that we have a word $w$ that is $(u, v)$-weakly renormalizable, then it is of the form $w = w_1 \gamma w_2$, where $\gamma$ is written in terms of $u$ and $v$. The word $\gamma$ cannot contain factors of the form $uu \ldots vv$ or $vv \ldots uu$ as discussed in the proof of that lemma, so it is written as powers of $u$ (respectively, $v$) followed by single instances of $v$ (respectively, $u$). Hence, we can choose a new alphabet $(\alpha, \beta) = (u, uv)$ (respectively, $(\alpha, \beta) = (uv, v)$) so that all exponents are now reduced by $1$ when $\gamma$ is written in the new alphabet $(\alpha, \beta)$. This simplifies the structure of the renormalization kernel at the cost of making the alphabet more complex. We will apply the renormalization algorithm inductively a certain number of times to ensure that the complexity of both the renormalization kernel and the alphabet remain reasonable.

        \begin{definition}\label{def:semirenormalizable}
        Let $(\alpha,\beta)\in \mathscr{A}$ and $w$ be a finite word in the alphabet $\{1,2\}$. We say that $w$ is \emph{$(\alpha,\beta)$-semi renormalizable} if there is an extension $\tilde{w}$ of $w$ of at most two digits, one to the left and one to the right such that $\tilde{w}$ is \emph{$(\alpha,\beta)$-weakly renormalizable}.
        \end{definition}

    Let us try to justify this definition: There are subwords of words 
    written in the alphabet $\{a,b\}$ that can fail to be weakly renormalizable for any alphabet with nontrivial kernel, because they are missing one digit at one (or both) of their ends. For example, the word of even length $w=21\ldots 1$ is a subword of $b^\infty ab^\infty$, and hence it belongs to $\Sigma(3,n)$. However, as is easy to check, it can only be exhibited as an $(\alpha,\beta)$-weakly renormalizable word by $w=w_1w_2$. But $w$ is $(a, b)$-semi renormalizable, with nontrivial kernel, since $2w1$ is a word in $\{a,b\}$.

\section{Farey sequences and subhorseshoes}
There is a well-known way to list all the rational numbers in the interval $[0,1]$. Given $n\in \mathbb{N}$ let us define first the Farey sequence of order $n$, $F_n$ as the increasing sequence of positive fractions in which numerator and denominator have no common positive divisor other than one and in which the denominator is less than or equal to $n$. The Farey sequence of order $n$ contains all of the members of the Farey sequences of lower orders. In particular $F_n$ contains all of the members of $F_{n-1}$ and also contains an additional fraction for each number that is less than $n$ and relatively prime with $n$. Thus, one has $\abs{F_n}=\abs{F_{n-1}}+\Phi(n)$ where $\Phi(n)$ denotes the usual Euler function that counts the number of non-negative integers less than $n$ that are relatively prime with $n$. As $\abs{F_1}=\abs{\{0,1\}}=2$, then the formula for the length of the Farey sequence of order $n$ is:
$$\abs{F_n}=1+\sum\limits_{k=1}^{n}\Phi(k).$$

We said that two elements $p/q, r/s\in F_n$ with $p/q<r/s$ are \emph{consecutive} if they are successive elements of the sequence $F_n$. One has the following properties:
\begin{itemize}
    \item If $p/q<r/s$ are consecutive in $F_n$ for some $n$ then $ps-qr=1$.
    \item If $p/q<r/s$ and $ps-qr=1$, then $p/q$ and $r/s$ are consecutive in $F_n$ for $\max \{q,s\}\leq n<q+s$. And in $F_{q+s}$, the fractions $p/q$ and $r/s$ are separated just for the element $\frac{p}{q}\oplus\frac{r}{s}=\frac{p+r}{q+s}$, with $p+r$ and $q+s$ relatively primes, which is called the \emph{mediant} of $p/q$ and $r/s$.
    \end{itemize}
Then, starting with the rational numbers $0$ and $1$ and recursively, taking mediants, one shows that this process generates all the rational numbers in the interval $(0,1)$ as a mediant exactly once. 

Given $\kappa\in \mathsf{c}(\mathscr{A})\cup \{a,b\}$ we associate the rational number $\theta(\kappa)\in [0,1]$, given by the proportion of letters $b$ in the word $\kappa$, that is $\theta(\kappa)=\frac{\abs{\kappa}_b}{\abs{\kappa}}$ where $\abs{\kappa}_b$ is the number of letters $b$ in $\kappa$, and similarly, $\abs{\kappa}_a$ is the number of letters $a$ in $\kappa$. Note that for any $(\alpha,\beta)\in \mathscr{A}$ we have $\theta(\alpha\beta)=\theta(\alpha)\oplus\theta(\beta)$. 

Given any $p/q\in (0,1)$ with $p$ and $q$ relatively primes, we have that $q-p$ and $p$ are also relatively primes, then using Theorem $18$ of \cite{Bombi} one concludes that for some unique ordered alphabet $(\alpha,\beta)$, one has that $\abs{\alpha\beta}_a=q-p$ and $\abs{\alpha\beta}_b=p$, so we conclude that $\theta(\alpha\beta)=\frac{\abs{\alpha\beta}_b}{\abs{\alpha\beta}}= \frac{\abs{\alpha\beta}_b}{\abs{\alpha\beta}_a+\abs{\alpha\beta}_b} =\frac{p}{p+(q-p)}=\frac{p}{q}$ and then $\theta$ is a bijection. In this way, we have the
\begin{lemma}\label{unique}
For each $n\in \mathbb{N}$ we can find unique words $\kappa_1,\dots,\kappa_{\abs{F_n}}\in \mathsf{c}(\mathscr{A})\cup \{a,b\}$ such that $F_n=\{\theta(\kappa_1),\dots,\theta(\kappa_{\abs{F_n}})\}$ and $(\kappa_i,\kappa_{i+1})\in \mathscr{A}$
for each $1\leq i\leq \abs{F_n}-1$. 
\end{lemma}

\begin{proof}
The proof is by induction: $F_1=\{\theta(a),\theta(b)\}$ and $(a,b)\in \mathscr{A}$. Suppose the lemma is true for $k=n$, by hypothesis, for each element of $x\in F_{n+1}\setminus F_n$ one can find $1\leq i\leq \abs{F_n}-1$ such that $(\kappa_i,\kappa_{i+1})\in \mathscr{A}$
and $x=\theta(\kappa_i)\oplus\theta(\kappa_{i+1})=\theta(\kappa_i\kappa_{i+1})$. This is enough because by definition $(\kappa_i, 
\kappa_i\kappa_{i+1})$ and $(\kappa_i\kappa_{i+1},\kappa_{i+1})$ are both ordered alphabets. 
\end{proof}

\begin{remark}
From the proof of the previous lemma, one can perceive the correspondence between taking mediants and the construction of the vertex of the tree $T$ i.e., the elements of $\mathscr{A}.$  
\end{remark}

\begin{lemma}
For any $(\alpha,\beta)\in\mathscr{A}$ we have $\theta(\alpha)<\theta(\beta)$.    
\end{lemma}
\begin{proof}
  The proof is by induction: $0=\theta(a)<\theta(b)=1$, suppose that for some ordered alphabet $\theta(\alpha)=\frac{\abs{\alpha}_b}{\abs{\alpha}}<\frac{\abs{\beta}_b}{\abs{\beta}}=\theta(\beta)$ then 
  $$\theta(\alpha\beta)= \frac{\abs{\alpha}_b+\abs{\beta}_b}{\abs{\alpha}+\abs{\beta}}  <\frac{\abs{\beta}_b}{\abs{\beta}}=\theta(\beta)$$
and 
$$ \theta(\alpha)=\frac{\abs{\alpha}_b}{\abs{\alpha}}<\frac{\abs{\alpha}_b+\abs{\beta}_b}{\abs{\alpha}+\abs{\beta}}  =\theta(\alpha\beta).$$
This finishes the proof.  
\end{proof}

\begin{corollary}\label{farey1}
 Given $\alpha,\beta\in \mathsf{c}(\mathscr{A})\cup \{a,b\}$ one has: $(\alpha, \beta)\in\mathscr{A}$ if, and only if, $\theta(\alpha)$ and $\theta(\beta)$ are consecutive in $F_{\max\{\abs{\alpha},\abs{\beta}\}}$.    
\end{corollary}

\begin{proof}
Let us show first by induction in $\max\{\abs{\alpha},\abs{\beta}\}$ that $(\alpha,\beta)$ is an ordered alphabet implies that $\theta(\alpha)$ and $\theta(\beta)$ are consecutive in $F_{\max\{\abs{\alpha},\abs{\beta}\}}$: If $\abs{\alpha}=\abs{\beta}=1$ then $\alpha=a$, $\beta=b$, $\theta(\alpha)=0$ and $\theta(\beta)=1$ and $a$ and $b$ are consecutive in $F_{\max\{\abs{\alpha},\abs{\beta}\}}=F_1$. Suppose the affirmation holds for $k\leq n$ and let $\alpha, \beta\in \mathsf{c}(\mathscr{A})\cup \{a,b\}$ such that $\max\{\abs{\alpha},\abs{\beta}\}=n+1$. Suppose first $(\alpha,\beta)=(\tilde{\alpha}\beta,\beta)$ where $(\tilde{\alpha},\beta)\in \mathscr{A}$, then $\max\{\abs{\tilde{\alpha}},\abs{\beta}\}<\max\{\abs{\alpha},\abs{\beta}\}=n+1$ and, by hypothesis, we have that $\theta(\tilde{\alpha})<\theta(\beta)$ are consecutive in $F_{\max\{\abs{\tilde{\alpha}},\abs{\beta}\}}$, and then $\theta(\tilde{\alpha})$,   $\theta(\alpha)=\theta(\tilde{\alpha}\beta)$ and $\theta(\beta)$ are consecutive in $F_{\abs{\tilde{\alpha}\beta}}=F_{\abs{\alpha}}=F_{\max \{\abs{\alpha},\abs{\beta}\}}=F_{n+1}$. Now, if $(\alpha,\beta)=(\alpha,\alpha\tilde{\beta})$, where $(\alpha,\tilde{\beta})\in \mathscr{A}$, then $\max\{\abs{\alpha},\abs{\tilde{\beta}}\}<\max\{\abs{\alpha},\abs{\beta}\}=n+1$ and,
by hypothesis, one has that $\theta(\alpha)<\theta(\tilde{\beta})$ are consecutive in $F_{\max\{\abs{\alpha},\abs{\tilde{\beta}}\}}$, and then $\theta(\alpha)$,  $\theta(\beta)=\theta(\alpha\tilde{\beta})$ and $\theta(\tilde{\beta})$ are consecutive in $F_{\abs{\alpha\tilde{\beta}}}=F_{\abs{\beta}}=F_{\max \{\abs{\alpha},\abs{\beta}\}}=F_{n+1}$ which concludes the induction step.

The other implication, is a direct consequence of Lemma \ref{unique}.
\end{proof}

Given a finite word $w=a_1\dots a_n$ with $n\geq2$ and $a_i\in \{a,b \}$ for each $i$, we define $w^+=a_2\dots a_n$ and $w^-=a_1\dots a_{n-1}$; naturally, if $n=1$ then we set $a_1^+=a_1^-=\emptyset$, the empty word.

\begin{lemma}\label{comienzo}
For every $(\alpha,\beta)\in \mathscr{A}$, $\alpha$ starts with $a$, $\beta$ ends with $b$. Moreover, every word $\alpha^k\beta$, with $k\ge 1$, starts with $\beta^-a$, and every word $\alpha\beta^k$, with $k\ge 1$, ends with $b\alpha^+$. In particular, we have the equality $\alpha\beta=\beta^-ab\alpha^+$.
\end{lemma}

\begin{proof}
Lemma $3.8$ of \cite{ERMR}.
\end{proof}

\begin{lemma}\label{farey}
Let $\alpha, \beta, \gamma\in \mathsf{c}(\mathscr{A})\cup \{a,b\}$ such that the rational numbers $\theta(\alpha),\theta(\beta), \theta(\gamma)$ are consecutive in $F_{\max\{\abs{\alpha},\abs{\beta},\abs{\gamma}\}}$. Then, $\alpha\beta\gamma\notin \Sigma(3,\abs{ \alpha\beta\gamma})$, but  $(\alpha\beta\gamma)^+, (\alpha\beta\gamma)^- \in \Sigma(3,\abs{ \alpha\beta\gamma}-2)$.

\end{lemma} 

\begin{proof} 
We will consider some cases:
\textbf{i)} $\abs{\alpha}, \abs{\gamma}<\abs{\beta}:$ In this case, $\theta(\alpha)$ and $\theta(\gamma)$ are consecutive in $F_{\max\{\abs{\alpha},\abs{\gamma} \}}$ and then by Lemma \ref{farey1} $(\alpha,\gamma)$ is an ordered alphabet. As the rational that is between $\theta(\alpha)$ and $\theta(\gamma)$ in $F_{\max\{\abs{\alpha},\abs{\beta},\abs{\gamma}\}}$ is $\theta(\alpha\gamma)$ then by injectivity of $\theta$ one concludes that
$\beta=\alpha\gamma$ and $\alpha\beta\gamma=\alpha\alpha\gamma\gamma\notin \Sigma(3,\abs{\alpha\beta\gamma})$ by Lemma \ref{lem:aa_bbr}. On the other hand, $\alpha\alpha\gamma\gamma^-a$ is a prefix of $\alpha\alpha\gamma\alpha\gamma\in \Sigma(3,\abs{\alpha\alpha\gamma\alpha\gamma})$, and $b\alpha^+\alpha\gamma\gamma$ is a suffix of $\alpha\gamma\alpha\gamma\gamma\in \Sigma(3,\abs{ \alpha\gamma\alpha\gamma\gamma})$ because $(\alpha\alpha\gamma,\alpha\gamma)$ and $(\alpha\gamma,\alpha\gamma\gamma)$ are ordered alphabets.

\textbf{ii)} $\abs{\alpha}<\abs{\beta}<\abs{\gamma}:$ In this case, as $\theta(\alpha)$ and $\theta(\beta)$ are consecutive in $F_{\max\{\abs{\alpha},\abs{\beta} \}}$ and $\theta(\beta)$ and $\theta(\gamma)$ are consecutive in $F_{\max\{\abs{\beta},\abs{\gamma} \}}$, by Lemma \ref{farey1} one gets that $(\alpha,\beta)$ and $(\beta,\gamma)$ are ordered alphabets. One concludes then that $\beta=\alpha\tilde{\beta}$ and $\gamma=\beta\tilde{\gamma}$ where $(\alpha,\tilde{\beta})$ and $(\beta,\tilde{\gamma})$ are ordered alphabets and as $\abs{\alpha}+\abs{\beta}>\abs{\gamma}$ (in other case $\theta(\alpha\beta)$ is between $\theta(\alpha)$ and $\theta(\beta)$ in $F_{\abs{\gamma}}$), one has $\abs{\tilde{\gamma}}<\abs{\alpha}$. Now, again by Lemma \ref{farey1}, $\theta(\alpha)$,  $\theta(\beta)$ and $\theta(\tilde{\gamma})$ are consecutive in $F_{\abs{\beta}}$. As the rational that is between $\theta(\alpha)$ and $\theta(\tilde{\gamma})$ in $F_{\abs{\beta}}$ is $\theta(\alpha\tilde{\gamma})$ then by injectivity of $\theta$ one concludes that $\beta=\alpha\tilde{\gamma}$ and then $\tilde{\beta}=\tilde{\gamma}$. Then $\alpha\beta\gamma=\alpha\alpha\tilde{\gamma}\alpha\tilde{\gamma}\tilde{\gamma}\notin \Sigma(3,\abs{\alpha\beta\gamma})$ but, $\alpha\alpha\tilde{\gamma}\alpha\tilde{\gamma}\tilde{\gamma}^-a$ is a prefix of $\alpha\alpha\tilde{\gamma}\alpha\tilde{\gamma}\alpha \tilde{\gamma}\in \Sigma(3,\abs{ \alpha\alpha\tilde{\gamma}\alpha\tilde{\gamma}\alpha \tilde{\gamma}})$, and $b\alpha^+\alpha\tilde{\gamma}\alpha\tilde{\gamma}\tilde{\gamma}$ is a suffix of $\alpha\tilde{\gamma}\alpha\tilde{\gamma}\alpha\tilde{\gamma}\tilde{\gamma}\in \Sigma(3,\abs{ \alpha\tilde{\gamma}\alpha\tilde{\gamma}\alpha\tilde{\gamma}\tilde{\gamma}})$ because $(\alpha\alpha\tilde{\gamma}\alpha\tilde{\gamma},\alpha \tilde{\gamma})$ and $(\alpha\tilde{\gamma},\alpha\tilde{\gamma}\alpha\tilde{\gamma}\tilde{\gamma})$ are alphabets.


\textbf{iii)} $\abs{\gamma}<\abs{\beta}<\abs{\alpha}:$ In this case, as $\theta(\alpha)$ and $\theta(\beta)$ are consecutive in $F_{\abs{\alpha}}$ and $\theta(\beta)$ and $\theta(\gamma)$ are consecutive in $F_{\abs{\beta}}$, by Lemma \ref{farey1} one gets that $(\alpha,\beta)$ and $(\beta,\gamma)$ are ordered alphabets. Then $\alpha=\tilde{\alpha}\beta$, $\gamma=\beta\tilde{\gamma}$ where $(\tilde{\alpha},\beta)$ and $(\tilde{\beta},\gamma)$ are ordered alphabets and as $\abs{\beta}+\abs{\gamma}>\abs{\alpha}$ (in other case $\theta(\beta\gamma)$ is between $\theta(\beta)$ and $\theta(\gamma)$ in $F_{\abs{\alpha}}$), one has $\abs{\tilde{\alpha}}<\abs{\gamma}$. Now, again by Lemma \ref{farey1}, $\theta(\tilde{\alpha})$, $\theta(\beta)$ and $\theta(\gamma)$ are consecutive in $F_{\abs{\beta}}$. As the rational that is between $\theta(\tilde{\alpha})$ and $\theta(\gamma)$ in $F_{\abs{\beta}}$ is $\theta(\tilde{\alpha}\gamma)$ then one concludes that $\beta=\tilde{\alpha}\gamma$. Then $\alpha\beta\gamma=\tilde{\alpha}\tilde{\alpha}\gamma\tilde{\alpha}\gamma\gamma\notin \Sigma(3,\abs{\alpha\beta\gamma})$ but, $\tilde{\alpha}\tilde{\alpha}\gamma\tilde{\alpha}\gamma\gamma^-a$ is a prefix of $\tilde{\alpha}\tilde{\alpha}\gamma\tilde{\alpha}\gamma\tilde{\alpha}\gamma\in \Sigma(3,\abs{ \tilde{\alpha}\tilde{\alpha}\gamma\tilde{\alpha}\gamma\tilde{\alpha}\gamma})$, and $b\tilde{\alpha}^+\tilde{\alpha}\gamma\tilde{\alpha}\gamma\gamma$ is a suffix of $\tilde{\alpha}\gamma\tilde{\alpha}\gamma\tilde{\alpha}\gamma\gamma\in \Sigma(3,\abs{ \tilde{\alpha}\gamma\tilde{\alpha}\gamma\tilde{\alpha}\gamma\gamma})$ because $(\tilde{\alpha}\tilde{\alpha}\gamma\tilde{\alpha}\gamma,\tilde{\alpha}\gamma)$ and $(\tilde{\alpha}\gamma,\tilde{\alpha}\gamma\tilde{\alpha}\gamma\gamma)$ are ordered alphabets.


\textbf{iv)} $\abs{\beta}<\abs{\alpha},\abs{\gamma}:$ In this case, $\alpha=\tilde{\alpha}\beta^k$ for some $k\geq 1$ with $\abs{\tilde{\alpha}}<\abs{\beta}$ and then by Lemma \ref{farey1}, $\theta(\tilde{\alpha})$ is consecutive with $\theta(\beta)$ in $F_{\abs{\beta}}$. Analogously, $\gamma=\beta^r\tilde{\gamma}$ for some $r\geq 1$ with $\abs{\tilde{\gamma}}<\abs{\beta}$ and then $\theta(\beta)$ is consecutive with $\theta(\tilde{\gamma})$ in $F_{\abs{\beta}}$. Then we must have $\beta=\tilde{\alpha}\tilde{\gamma}$ and $\alpha\beta\gamma=\tilde{\alpha}\beta^{k+r+1}\tilde{\gamma}=\tilde{\alpha}(\tilde{\alpha}\tilde{\gamma})^{k+r+1}\tilde{\gamma}\notin \Sigma(3,\abs{\alpha\beta\gamma})$ but, $\tilde{\alpha}(\tilde{\alpha}\tilde{\gamma})^{k+r+1}\tilde{\gamma}^-a$ is a prefix of 
$\tilde{\alpha}(\tilde{\alpha}\tilde{\gamma})^{k+r+2}\in \Sigma(3,\abs{\tilde{\alpha}(\tilde{\alpha}\tilde{\gamma})^{k+r+2}})$ and $b\tilde{\alpha}^+(\tilde{\alpha}\tilde{\gamma})^{k+r+1}\tilde{\gamma}$ is a suffix of 
 $(\tilde{\alpha}\tilde{\gamma})^{k+r+2}\tilde{\gamma}\in \Sigma(3,\abs{(\tilde{\alpha}\tilde{\gamma})^{k+r+2}\tilde{\gamma}},)$ and $(\tilde{\alpha}(\tilde{\alpha}\tilde{\gamma})^{k+r+1},\tilde{\alpha}\tilde{\gamma})$ and $(\tilde{\alpha}\tilde{\gamma},(\tilde{\alpha}\tilde{\gamma})^{k+r+1}\tilde{\gamma})$ are alphabets. This finishes the proof of the lemma.
\end{proof}

Given an ordered alphabet $(\alpha,\beta)$ consider the periodic orbit, $\psi_{\alpha\beta}$, determined by the periodic point $\Pi^{-1}(\overline{\alpha\beta})$. Similarly, define $\psi_a$ and $\psi_b$ as the fixed orbits given by the points $\Pi^{-1}(\overline{a})$ and  $\Pi^{-1}(\overline{b})$ respectively.

\begin{proposition}\label{ferraduradelcomienzo}
 Given $\epsilon>0$ there exist a subhorseshoe $\tilde{\Lambda}^{\epsilon}\subset (\Lambda(2))_{3+\epsilon}$ with the property that for any $p\in \mathsf{c}(\mathscr{A})\cup \{a,b\}$ one has $\psi_p\subset\tilde{\Lambda}^{\epsilon}$.
    
\end{proposition}

\begin{proof}

Given $n\in \mathbb{N}$, suppose that the associated words by $\theta^{-1}$ to the Farey sequence $F_n$ are $\alpha_0,\alpha_1,\dots,\alpha_N$ in that order. If $x$ is a letter in the word $\alpha_j=c_1c_2\dots c_r$, lets say $x=c_s$ where $1\leq s\leq r$. We have the following cases, where $j$ is such that the adjacent words considered are defined in each case:

\textbf{i)} $s-1\geq n:$ In this case, as $\abs{\alpha_{j+1}\alpha_{j+2}}\geq 2(n+1)$, $x$ belongs to the word $\alpha_j\alpha_{j+1}\alpha_{j+2}^-\in \Sigma(3,\abs{\alpha_j\alpha_{j+1}\alpha_{j+2}}-2)$ by Lemma \ref{farey} and the distance of $x$ to both extremes of that word is greater than $n$.

\textbf{ii)} $s-1\geq n$ and $\abs{\alpha_{j+1}}> n:$ In this case $x$ belongs to the word $\alpha_j\alpha_{j+1}\in \Sigma(3,\abs{\alpha_j\alpha_{j+1}})$ by Lemma \ref{farey} and the distance of $x$ to both extremes of that word is greater than $n$.

\textbf{iii)} $r-s\geq n:$ In this case, as $\abs{\alpha_{j-2}\alpha_{j-1}}\geq 2(n+1)$, $x$ belongs to the word $\alpha_{j-2}^+\alpha_{j-1}\alpha_j\in \Sigma(3,\abs{\alpha_{j-2}\alpha_{j-1}\alpha_j}-2)$ by Lemma \ref{farey} and the distance of $x$ to both extremes of that word is greater than $n$.

\textbf{iv)} $r-s\geq n$ and $\abs{\alpha_{j-1}}>n:$ In this case $x$ belongs to the word $\alpha_{j-1}\alpha_j\in \Sigma(3,\abs{\alpha_{j-1}\alpha_j})$ by Lemma \ref{farey} and the distance of $x$ to both extremes of that word is greater than $n$.

\textbf{v)} $s-1, r-s< n:$ In this case, as $\abs{\alpha_j\alpha_{j+1}}\geq 2(n+1)$ and $\abs{\alpha_{j-1}\alpha_j} \geq 2(n+1)$, deleting the first or the last letter of $\alpha_{j-1}\alpha_j\alpha_{j+1}$ by Lemma \ref{farey} we obtain a word of $\kappa\in \Sigma(3,\abs{\alpha_{j-1}\alpha_j\alpha_{j+1}}-2)$ and the distance of $x$ to both extremes of $\beta$ that word is again greater than $n$.

Given an alphabet $(\alpha,\beta)$ and $\tilde{n}\in \mathbb{N}$, consider the Farey sequence $F_n$, where $n=\abs{(\alpha\beta)^{\tilde{n}}\beta}/2$ (remember that $(\alpha\beta,(\alpha\beta)^{\tilde{n}-1}\beta)$ is an alphabet). If the associated words to $F_n$ are $\alpha_0,\alpha_1,\dots,\alpha_N$ and $\alpha_i= 
(\alpha\beta)^{\tilde{n}}\beta$, then we affirm that for $\tilde{n}$ large enough, the sequence $\theta_{\alpha\beta,b}=\overline{\alpha\beta}^-\alpha_i\alpha_{i+1}\dots \alpha_N \overline{b}^+\in \Pi((\Lambda(2))_{3+\epsilon/4})=\Sigma(3+\epsilon/4)$ where $\overline{\alpha\beta}^-$ is the infinite periodic sequence to the left with period $\alpha\beta$ and $\overline{b}^+$ is the infinite periodic sequence to the right with period $b$.  Given any letter $x$ in any of the words $\alpha_j$ with $j=i+2,\dots, N-2$ by the items \textbf{i},\textbf{iii} and \textbf{v}, one has that the word of size $2n+1$ contained in $\alpha_i\alpha_{i+1}\dots \alpha_N$ and centered in that letter is a word in $\Sigma(3,2n+1)$. Additionally, if $x=c_s$ is a letter in the last $\alpha\beta\beta$ of $\alpha_i$ where $\alpha_i=c_1c_2\dots c_r$ then by choosing $\tilde{n}$ large, one has that $s-1>n$ and then by item \textbf{i}, the previous affirmation also holds for that $x$. In the other positions of $\alpha_i$ and any position corresponding to $\overline{\alpha\beta}^-$ the Markov value is less than $3$ because by Theorem $15$ of \cite{Bombi} one has  
that for any ordered alphabet $(\alpha,\beta)$, the number $\abs{\alpha^+\beta^-}$ is equal to the supremum of sizes of words $w$ such that the cut $bw^*a|bwa$ appears in the sequence $\overline{\alpha\beta}$. In particular, $\lambda$ has value less than $3$ in any position of the $\alpha\beta$ in the middle of $\alpha\beta\alpha\beta\alpha\beta$ (we will show later a much more general result). If $x$ is a letter of $\alpha_{i+1}$, as $\abs{\alpha_i}=2n$ one concludes by items \textbf{i}, \textbf{iv} and \textbf{v} again the same. Finally, as $\alpha_{N-1}=\theta^{-1}(\frac{n-1}{n})=ab^{n-1}=c_1c_2\dots c_r$ (which corresponds to the alphabet $(ab^{n-2},b)$) then in any position $x=c_s$ of the $a$ of $\alpha_{N-1}$ one has $r-s>n$, then by item \textbf{iii} one has that the word of size $2n+1$ contained in $\alpha_i\alpha_{i+1}\dots \alpha_N$ and centered in that letter is a word in $\Sigma(3,2n+1)$. As in any position of any $b$, $\lambda$ has value less than $3$, the affirmation follows for $\tilde{n}$ large. 

Now, again consider the alphabet $(\alpha,\beta)$ and the Farey sequence $F_n$, where now $n=\abs{(\alpha(\alpha\beta)^{\tilde{n}}}/2$ and $\tilde{n}\in\mathbb{N}$ (remember that $(\alpha(\alpha\beta)^{\tilde{n}-1},\alpha\beta)$ is an alphabet). If the associated words to $F_n$ are $\alpha_0,\alpha_1,\dots,\alpha_N$ and $\alpha_i= 
\alpha(\alpha\beta)^{\tilde{n}}$, then we affirm that for $\tilde{n}$ large enough, the sequence $\theta_{a,\alpha\beta}=\overline{a}^-\alpha_0\alpha_1\dots \alpha_{i-1}\alpha_i \overline{\alpha\beta}^+\in \Sigma(3+\epsilon/4)$. Given any letter $x$ in any of the words $\alpha_j$ with $j=2,\dots, i-2$ by the items \textbf{i},\textbf{iii} and \textbf{v}, one has that the word of size $2n+1$ contained in $\alpha_0\alpha_1\dots \alpha_i$ and centered in that letter is a word in $\Sigma(3,2n+1)$. Additionally, If $x=c_s$ is a letter in the first $\alpha\alpha\beta$ of $\alpha_i$ where $\alpha_i=c_1c_2\dots c_r$ then by choosing $\tilde{n}$ large, one has that $r-s>n$ and then by item \textbf{iii}, the previous affirmation also holds for that $x$. In the other positions of $\alpha_i$ and any position corresponding to $\overline{\alpha\beta}^+$ the Markov value is less than $3$ as before. If $x$ is a letter of $\alpha_{i-1}$, as $\abs{\alpha_i}=2n$ one concludes by items \textbf{ii}, \textbf{iii} and \textbf{v} again the same. Finally, as $\alpha_1=\theta^{-1}(\frac{1}{n})=a^{n-1}b=c_1c_2\dots c_r$ (which corresponds to the alphabet $(a,a^{n-2}b)$), then in the position $2=c_{r-2}=c_s$ of $\alpha_1$ one has $s-1>n$, then by item \textbf{i} one has that the word of size $2n+1$ contained in $\alpha_0\alpha_1\dots \alpha_{i-1}\alpha_i$ and centered in that letter is a word in $\Sigma(3,2n+1)$. As in any $2$ in the middle of $222$, $\lambda$ has value less than $3$, the affirmation follows for $\tilde{n}$ large.

Analogously given $n\in \mathbb{N}$, if we consider the Farey sequence $F_n$ with associated words $\alpha_0,\alpha_1,\dots,\alpha_N$, in that order, then one can show for $n$ large, using the same arguments as before, that the sequence $\theta_{a,b}=\overline{a}^-\alpha_0\alpha_1\dots \alpha_N \overline{b}^+\in \Sigma(3+\epsilon/4)$. Of course, this also implies that $\theta_{b,a}=\theta_{a,b}^*=\overline{b}^-\alpha_N^*\alpha_{N-1}^*\dots \alpha_0^* \overline{a}^+\in \Sigma(3+\epsilon/4)$.

The sequences $\theta_{a,b}$ and $\theta_{b,a}$ determine two points $x,y \in \Lambda$ such that $x\in W^u(\psi_a)\cap W^s(\psi_b)$, $y\in W^u(\psi_b)\cap W^s(\psi_a)$ and 
$\mathcal{O}(x) \cup \mathcal{O}(y) \subset (\Lambda(2))_{3+\epsilon/4}$. Analogously, given one alphabet $(\alpha_1,\beta_1)$, the sequence $\theta_{a,\alpha_1\beta_1}$ determine a point $x\in W^u(\psi_a)\cap W^s(\psi_{\alpha_1\beta_1})$ with $m_{\varphi,f}(x)\leq 3+\epsilon/2$ and we can use the sequences $\theta_{\alpha_1\beta_1,b}$ and $\theta_{b,a}$ to find some point $y\in W^u(\psi_{\alpha_1\beta_1})\cap W^s(\psi_a)$ with $m_{\varphi,f}(y)\leq 3+\epsilon/2$. In any case, Proposition \ref{connection11} let us find subhorseshoes $\Lambda^1,\Lambda^2\subset (\Lambda(2))_{3+\frac{3\epsilon}{4}}$ such that $\psi_a\cup\psi_b\subset \Lambda^1$ and  $\psi_a\cup\psi_{\alpha_1\beta_1}\subset \Lambda^2$. 


Now, let $r(\epsilon)\in \mathbb{N}$ sufficiently large such that if $\tilde{\alpha}=(a_{-r(\epsilon)},\dots,a_0\dots, a_{r(\epsilon)})\in \{1,2\}^{2r(\epsilon)+1}$
and $\tilde{x}, \tilde{y}\in R(\tilde{\alpha};0)\ \mbox{then}\ \abs{f(\tilde{x})-f(\tilde{y})}<\epsilon/4.$ Consider the set
$$P(\epsilon)=\bigcap \limits_{n \in \mathbb{Z}} \varphi ^{-n}(\bigcup \limits_{\tilde{\alpha} \in C(\epsilon)}  R(\tilde{\alpha};0)),$$
where 
$$C(\epsilon)=\{\tilde{\alpha}\in \{1,2\}^{2r(\epsilon)+1}:R(\tilde{\alpha};0)\cap (\Lambda(2))_{3+\frac{3\epsilon}{4}}\neq \emptyset \}.$$
Note that by construction, $(\Lambda(2))_{3+\frac{3\epsilon}{4}}\subset P(\epsilon)\subset (\Lambda(2))_{3+\epsilon}$ and $P(\epsilon)$ is a hyperbolic set of finite type. As the periodic orbits $\psi_p\subset\tilde{\Lambda}^{\epsilon}$ with $p\in \mathsf{c}(\mathscr{A})\cup \{a,b\}$ all belong to the same transitive component of $P(\epsilon)$ it follows the existence of the subhorseshoe $\tilde{\Lambda}^{\epsilon}$ as in the statement of the proposition. 
\end{proof}

\section{Some results about cuts and finite words}

We will usually start with cuts written in the alphabet $\{a,b\}$ and then get more sophisticated cuts by applying the Nielsen substitutions. Then, for our purposes, it is necessary to know exactly the effect of applying any $W$ in the alphabet $\{U,V \}$ to good and bad cuts. We start with some preliminary results.


\begin{corollary}\label{9}
 Let $(u,v)\in \mathscr{A}$, an ordered alphabet. If for some $j\geq 2$, $\abs{u^j}\geq \abs{uv}$ then $u^j$ begins with $v^-a$. Analogously, if $\abs{v^j}\geq \abs{uv}$ then $v^j$ ends with $bu^+$.
 \end{corollary}
\begin{proof}
    For the first part, if $v=b$ then $v^-a=a$ and the result holds. If $u=\tilde{u}v$ where $(\tilde{u},v)$ is an alphabet, then $u=v^-ab\tilde{u}^+$. If $v=u^i\tilde{v}$ where $(u,\tilde{v})$ is an alphabet with $\abs{\tilde{v}}<\abs{u}$, one has $\abs{u^{i+1}}<\abs{uv}\leq \abs{u^j}$ and then $u^j$ has the same beginning of size $\abs{u^{i+1}}$ than $uu^i\tilde{v}=uv=v^-abu^+$ and $\abs{v^-a}=\abs{v}=\abs{u^i\tilde{v}}<\abs{u^{i+1}}$.

     For the second part, if $u=a$ then $bu^+=b$ and the result holds. If $v=u\tilde{v}$ where $(u,\tilde{v})$ is an alphabet, then $v=\tilde{v}^-abu^+$. Finally, if $u=\tilde{u}v^i$ where $(\tilde{u},v)$ is an alphabet with $\abs{\tilde{u}}<\abs{v}$, one has $\abs{v^{i+1}}<\abs{uv}\leq \abs{v^j}$ and then $v^j$ has the same end of size $\abs{v^{i+1}}$ than $\tilde{u}v^iv=uv=v^-abu^+$ and $\abs{bu^+}=\abs{u}=\abs{\tilde{u}v^i}<\abs{v^{i+1}}$.
\end{proof}

\begin{lemma} \label{lem:reverse}
    For any finite word $w$ in the alphabet $\{a, b\}$, we have the identities $bU(w^*) = U(w)^*b$ and $V(w^*)a = aV(w)^*$. 
\end{lemma}

\begin{proof}
Lemma $3.14$ of \cite{ERMR}.
\end{proof}

\begin{lemma}\label{transpuesto}
Given any words $w$ in the alphabet $\{a,b \}$ and $W$ in the alphabet $\{U,V \}$ one has 
$$(u^+ W(w) v^-)^*=u^+ W(w^*) v^-$$
where $u=W(a)$ and $v=W(b)$. In particular, $(W(awb))^*=bu^+ W(w^*) v^-a.$

\end{lemma}

\begin{proof}
The proof is by induction on $W$. If $W=\emptyset$, there is nothing to show. Assume that we have $(u^+ W(w) v^-)^*=u^+ W(w^*) v^-$. If $\tilde{W}=UW$, $\tilde{u}=\tilde{W}(a)=U(W(a))=U(u)=U(au^+)=abU(u^+)$ and $\tilde{v}=\tilde{W}(b)=U(W(b))=U(v)=U(v^-b)=U(v^-)b$ then $\tilde{u}^+=bU(u^+)$ and $\tilde{v}^-=U(v^-)$. Therefore
\begin{eqnarray*}
\tilde{u}^+\tilde{W}(w^*)\tilde{v}^-&=&bU(u^+)UW(w^*)U(v^-)=bU(u^+W(w^*)v^-)=bU((u^+W(w)v^-)^*)\\ &=& (U(u^+W(w)v^-))^*b=(bU(u^+W(w)v^-))^*=(bU(u^+)UW(w)U(v^-))^*\\ &=& (\tilde{u}^+\tilde{W}(w)\tilde{v}^-)^*.
\end{eqnarray*}
On the other hand, if now $\tilde{W}=VW$, $\tilde{u}=\tilde{W}(a)=V(W(a))=V(u)=V(au^+)=aV(u^+)$ and $\tilde{v}=\tilde{W}(b)=V(W(b))=V(v)=V(v^-b)=V(v^-)ab$ then $\tilde{u}^+=V(u^+)$ and $\tilde{v}^-=V(v^-)a$. Therefore
\begin{eqnarray*}
\tilde{u}^+\tilde{W}(w^*)\tilde{v}^-&=&V(u^+)VW(w^*)V(v^-)a=V(u^+W(w^*)v^-)a=V((u^+W(w)v^-)^*)a\\ &=& a(V(u^+W(w)v^-))^*=(V(u^+W(w)v^-)a)^*=(V(u^+)VW(w)V(v^-)a)^*\\ &=& (\tilde{u}^+\tilde{W}(w)\tilde{v}^-)^*.
\end{eqnarray*}
This finishes the proof by induction. Finally,
$$(W(aw b))^*=(u W(w) v)^*=(au^+ W(w) v^-b)^*=b(u^+ W(w) v^-)^*a=bu^+ W(w^*) v^-a.$$

\end{proof}
\begin{corollary}\label{1}
If there is a bad cut of the form $Xbw^*b|aw a Y$ where $w$, $X$ and $Y$ are words in the alphabet $\{a,b \}$ and $W$ is a word in $\{U,V \}$ such that $\abs{W(X)}\geq \abs{u}$ and $\abs{W(Y)}\geq \abs{v}$, then $W(Xbw^* b|a w a Y )$ contains the bad cut 
$$b u^+ W (w^*) v | u W(w) v^- a=b u^+ W (w^*) v^- b | a u^+ W(w) v^- a$$
where $u=W(a)$ and $v=W(b)$. 
\end{corollary}

\begin{proof}
As $\abs{W(X)}\geq \abs{u}$ and $\abs{W(Y)}\geq \abs{v}$ we have
$\abs{W(Xb)}\geq \abs{uv}$ and $\abs{W(aY)}\geq \abs{uv}$. If $X$ contains one $a$ then for some $i\geq 1$ we have that $W(Xb)$ ends with $uv^i$ which ends with $bu^+$ by Lemma \ref{comienzo}. If for some $j\geq 1$ we have $X=b^j$ then $\abs{W(Xb)}=\abs{v^{j+1}}\geq \abs{uv}$ and Corollary \ref{9}, implies that $W(Xb)$ ends with $bu^+$ again. On the other hand, if $Y$ contains one $b$ then for some $i\geq 1$ we have that $W(aY)$ begins with $u^iv$ which begins with $v^-a$ by Lemma \ref{comienzo}. If for some $j\geq 1$ we have $Y=a^j$ then $\abs{W(aY)}=\abs{u^{j+1}}\geq \abs{uv}$ and Corollary \ref{9}, implies that $W(aY)$ begins with $v^-a$ again. Then, we conclude that 
\begin{eqnarray*}
    W(Xbw^* b|a w a Y)&\supseteq& W(Xb) W (w^*) v | u W(w) W(aY)\\ &\supseteq& b u^+ W (w^*) v^- b | a u^+ W(w) v^- a.
    \end{eqnarray*}
Finally, by Lemma \ref{transpuesto} $(u^+ W(w) v^-)^*=u^+ W(w^*) v^-$, from this follows the result. 
\end{proof}

\begin{remark}
 Note that the hypothesis of the corollary holds if $X$ contains at least one $a$ and $Y$ contains at least one $b$.
 \end{remark}

\begin{corollary}\label{2}
If there is a bad cut of the form $aw^*a|bw b$ where $w$ is a word in the alphabet $\{a,b \}$ and $W$ is a word in $\{U,V \}$, then the cut determined by $W$ is a bad cut:
$$W(aw^* a|b w b)=au^+W(w^*)v^-a|bu^+W(w)v^-b$$  
where $u=W(a)$ and $v=W(b)$.         
\end{corollary}
\begin{proof}
This is a consequence of Lemma \ref{transpuesto} because
$$W(aw^* ab w b)=uW(w^*)uvW(w)v=au^+W(w^*)v^-abu^+W(w)v^-b$$
and $(u^+ W(w) v^-)^*=u^+ W(w^*) v^-$.
\end{proof}

\begin{corollary}\label{3}
If there is a good cut of the form $Xbw^*a|bw a Y$ where $w$, $X$ and $Y$ are words in the alphabet $\{a,b \}$ and $W$ is a word in $\{U,V \}$ such that $\abs{W(X)}\geq \abs{u}$ and $\abs{W(Y)}\geq \abs{v}$, then $W(Xbw^* a|b w a Y )$ contains the good cut determined by $W$:
$$b u^+ W (w^*) v^- a | b u^+ W(w) v^- a$$
where $u=W(a)$ and $v=W(b)$. 
\end{corollary}

\begin{proof}
As in the proof of Corollary \ref{1}, as $\abs{W(X)}\geq \abs{u}$ and $\abs{W(Y)}\geq \abs{v}$ we have that $W(Xb)$ ends with $bu^+$  $W(aY)$ begins with $v^-a$. Then, we conclude that 
\begin{eqnarray*}
    W(Xbw^* a|b w a Y )&\supseteq& W(Xb) W (w^*) v^- a|b u^+ W(w) W(aY)\\ &\supseteq& b u^+ W (w^*) v^- a|b u^+ W(w) v^- a.
    \end{eqnarray*}
Finally, by Lemma \ref{transpuesto} $(u^+ W(w) v^-)^*=u^+ W(w^*) v^-$, from this follows the result. 
\end{proof}

\begin{remark}
 As in Lemma \ref{1}, the hypothesis of the corollary holds if $X$ contains at least one $a$ and $Y$ contains at least one $b$.
 \end{remark}

\begin{corollary}\label{4}
If there is a good cut of the form $aw^*b|aw b$ where $w$ is a word in the alphabet $\{a,b \}$ and $W$ is a word in $\{U,V \}$, then the cut determined by $W$ is a good cut:
$$W(aw^* b|a w b)=au^+W(w^*)v^-b|au^+W(w)v^-b$$  
where $u=W(a)$ and $v=W(b)$.         
\end{corollary}
\begin{proof}
This is a consequence of Lemma \ref{transpuesto} because
$$W(aw^* b|a w b)=uW(w^*)v|uW(w)v=au^+W(w^*)v^-b|au^+W(w)v^-b$$
and $(u^+ W(w) v^-)^*=u^+ W(w^*) v^-$.
\end{proof}


\begin{proposition}\label{control}
Let $X,R$ and $Y$ be words in the alphabet $\{a,b \}$ and $W$ be a word in $\{U,V \}$ such that $\abs{W(X)}\geq \abs{u}$ and $\abs{W(Y)}\geq \abs{v}$. Additionally, suppose that all cuts of $XRY$ that includes letters of $R$ are good cuts. Then, all cuts that includes letters of $W(R)$ in $W(XRY)=W(X)W(R)W(Y)$ are good cuts.

\end{proposition}

\begin{proof}
Note first that if we apply $U$ or $V$ to one of the cuts $a|a$ or $b|b$ and after this is introduced some $a$, then in that $a$ the cuts are good. This is because $U(a|a)=ab|ab$ and $V(b|b)=ab|ab$. Additionally, note that $\abs{W^{'}(J(X))}\geq \abs{W(a)}\geq \abs {W^{'}(a)}$ and $\abs{W^{'}(J(Y))}\geq \abs{W(b)}\geq \abs {W^{'}(b)}$, where $W^{'}$ is the word obtained from $W$ by deleting its last letter $J$. The result follows by induction, using corollaries \ref{3} and \ref{4}.
\end{proof}




The following lemma shows the existence of the ordered alphabet with reasonable size that we were looking for. 
\begin{lemma}\label{12}
   Let $n\geq 68$ and let $w\in \Sigma(3+6^{-3n},3n)$. Then, there exits an ordered alphabet $(\alpha,\beta)\in \mathscr{A}$ satisfying $\abs{\alpha},\abs{\beta}<n$ and $\abs{\alpha \beta}\geq n$ such that $w$ is $(\alpha,\beta)$-semi renormalizable. 
\end{lemma}
\begin{proof}
Corollary $3.23$ of \cite{ERMR}.
\end{proof}

\begin{lemma}\label{algoritmo}
Let $w\in \Sigma(3+6^{-3n},3n)$, where $n\geq 61$, be a finite word. If $w$ is $(u, v)$-weakly renormalizable as $w = w_1 \gamma w_2$ with $\gamma \neq \emptyset$, then $w$ is $(\alpha, \beta)$-weakly renormalizable for some $(\alpha, \beta) \in \{(uv, v), (u, uv)\}$. 
\end{lemma}
\begin{proof}
As $$r(\omega)\leq (\abs{w}+1)\log(3+\sqrt{2})\leq (3n+1)\log(3+\sqrt{2})<(3\log 6)n.$$
Taking $r=(3\log 6)n$ and applying Lemma \ref{lem:renormalizingr} we obtain the result.
\end{proof}


\begin{lemma}\label{14}
If $(\alpha,\beta)$ is an ordered alphabet with $\abs{\alpha \beta}<n$ and $w\in \Sigma(3+6^{-3n},3n)$ contains $\alpha \beta$, then $w$ is $(\alpha,\beta)$-semi renormalizable as $\tilde{w}=w_1\gamma w_2$. If $w$ begins (or ends) with $\alpha \beta$ then $w_1=\emptyset$ ($w_2=\emptyset$).
\end{lemma}
\begin{proof}

First note that $w$ is trivially $(a,b)$-semi renormalizable, say $\tilde{w}=\gamma_0$ where $\gamma_0$ is a word in the alphabet $\{a,b\}$.  Now we apply inductively Lemma \ref{algoritmo} to obtain a sequence of alphabets $(A_j,B_j)\in\mathscr{A}_j$ such that for all $0\leq j\leq m$, the word $\tilde{w}$ is $(A_j,B_j)$-weakly renormalizable for each $j$ and $|A_mB_m|\geq n$. 
    
On the other hand, since $(\alpha,\beta)\in\mathscr{A}$ there exists a sequence of alphabets $(\alpha_i,\beta_i)\in\mathscr{A}_i$ such that $\alpha\beta$ can be written in the alphabet $\{\alpha_i,\beta_i \}$ for all $0\leq i\leq n$ and $(\alpha_n,\beta_n)=(\alpha,\beta)$. Since $\alpha\beta$ starts with $a=\alpha_0$ and ends with $b=\beta_0$, inductively we obtain that $\alpha\beta$ starts with $\alpha_i$ and ends with $\beta_i$. In particular $\alpha\beta$ contains $\alpha_i\beta_i$.  
    
Write $\tilde{w}=w_1\gamma_jw_2$ as in the definition of $(A_j,B_j)$-weakly renormalizable. Using the fact that $\alpha\beta$ contains $\alpha_j\beta_j$, gluing some words $\tau$ and $\tau'$ we get that $\tau\alpha_j\beta_j\tau'=A_jB_j\gamma_jA_jB_j$ is a word in the alphabet $\{A_j,B_j\}$. Hence, by Lemma \ref{lem:determine_alphabet}, we obtain that $(A_j,B_j)=(\alpha_j,\beta_j)$ for all $0\leq j\leq n$. In particular $m>n$, because otherwise  $n\leq|A_mB_m|=|\alpha_m\beta_m|<n$. This shows that $\tilde{w}$ is $(\alpha,\beta)$-weakly renormalizable.

Now assume that $w$ starts with $\alpha\beta$ (the other case is analogous). Observe that there is no need to complete the word to the left. We will show that $w_1=\emptyset$ for all $0\leq j\leq n$. Note that we already showed that $w_1$ is empty for $(\alpha_0,\beta_0)=(a,b)$. If $w_1$ becomes nonempty for $k+1$ for some $0\leq k\leq n$, it must happen that $\tilde{w}=\gamma_k w_2$ starts with $\beta_k$ (because of the renormalization algorithm, Lemma \ref{lem:renormalizingr}). But $w$ starts with $\alpha\beta$, which in turn starts with $\alpha_k^s\beta_k$, which leads to a contradiction because it starts with $\beta_k^-a$. Since $(\alpha_n,\beta_n)=(\alpha,\beta)$ this finishes the proof.
\end{proof}


The following lemma is the version of Lemma \ref{lem:aa_bbr} that we will use through the text.
\begin{lemma}\label{alfacuadrado}
In any word $w\in\Sigma(3+6^{-3n},3n+2)$, where $n\geq 61$, it cannot appear a subword $\tilde{w}$ written in the alphabet $\{\alpha,\beta\}$, associated to $(\alpha,\beta)\in \mathscr{A}$, beginning with $\alpha^2$ and ending with $\beta^2$. 
\end{lemma}
\begin{proof}
Let $\tilde{w}$ as in the statement of the lemma. As $r(\tilde{w})\leq (\abs{\tilde{w}}+1)\log(3+\sqrt{2})\leq (3n+3)\log(3+\sqrt{2})$, by Lemma \ref{lem:aa_bbr}, $\tilde{w}$ does not appear as a subword of a word of $$\Sigma(3+\exp({-(3n+1) \log(3+\sqrt{2})}))=\Sigma(3+(2+2\sqrt{2})^{-(3n+1)})\supset \Sigma(3+6^{-3n}).$$
\end{proof}

\begin{corollary}\label{force}
Let $(\alpha,\beta)$ some alphabet and $w\in \Sigma(3+6^{-3n},3n+2)$. Then, in any of the cases
\begin{itemize}
    \item $\alpha^{r_1}\beta\alpha^{r_2}\beta$ is a subword of $w$
    \item $\alpha\beta^{r_1}\alpha\beta^{r_2}\alpha\beta$ is a subword of $w$
    \end{itemize}
where $r_1,r_2\geq 1$. We conclude that $r_2\geq r_1-1$.
\end{corollary}

\begin{proof}
 In the first case, if $(\tilde{\alpha},\tilde{\beta})=(\alpha,\alpha^{r_2}\beta)$, one has
 $$\alpha^{r_1}\beta\alpha^{r_2}\beta= \alpha^{r_1-r_2}\alpha^{r_2}\beta\alpha^{r_2}\beta=\tilde{\alpha}^{r_1-r_2}\tilde{\beta}^2.$$
 In the second case, if $(\tilde{\alpha},\tilde{\beta})=(\alpha\beta^{r_2},\beta)$, by Lemma \ref{transpuesto} one has
\begin{eqnarray*}
     (\tilde{\alpha}^2\tilde{\beta}^{r_1-r_2})^*=(\tilde{\alpha}\tilde{\alpha}\tilde{\beta}^{r_1-r_2-1}\tilde{\beta})^*&=&b\tilde{\alpha}^+\tilde{\beta}^{r_1-r_2-1}\tilde{\alpha}\tilde{\beta}^-a =b\alpha^+\beta^{r_2}\beta^{r_1-r_2-1}\alpha\beta^{r_2}\beta^-a\\ &=&b\alpha^+\beta^{r_1-1}\alpha\beta^{r_2}\beta^-a\subseteq  \alpha\beta^{r_1}\alpha\beta^{r_2}\alpha\beta.
     \end{eqnarray*}
As $\Sigma(3+6^{-3n},3n+2)$ is closed by transpositions, Lemma \ref{alfacuadrado} implies in any case that $r_1-r_2\leq 1$.
\end{proof}
 
\section{Dimension of subhorseshoes that do not connect with $\psi_b$}

Remember that we are considering the horseshoe
$\Lambda(2)=C(2)\times\tilde{C}(2)$ equipped with the diffeomorphism $\varphi$ and the map $f$ as in Section \ref{Lambda}. Now, in \cite{M3}, it was proved for $s\leq \max f|_{\Lambda(2)}$ that $D(s)=HD(\ell^{-1}(-\infty,s])=HD(\pi^u((\Lambda(2))_s))$ where, as before,  $\pi^u((\Lambda(2))_s)$ is the projection of $(\Lambda(2))_s$ over the unstable Cantor set of $\Lambda(2)$ and in \cite{GCD} it was proved that 
$2HD(\pi^u((\Lambda(2))_s))=HD((\Lambda(2))_s).$ From this we conclude for $t<t_1=\sup \{s\in \mathbb{R}:d(s)<1\}$ that 
$$d(t)=\min \{1,2D(t)\}=2D(t)=HD((\Lambda(2))_t).$$
This section is devoted to the proof of the following proposition, which is the cornerstone of the proof of theorems \ref{teo1} and \ref{teo2}. 

\begin{proposition}\label{fundamental}
Let $n\in\mathbb{N}$ and $\tilde{\Lambda}$ some subhorseshoe of $\Lambda(2)$ such that $\max f|_{\tilde{\Lambda}}<3+6^{-3n}$. If for some $\epsilon>0$ with $\max f|_{\tilde{\Lambda}}+\epsilon< 3+6^{-3n}$, $\tilde{\Lambda}$ does not connect with $\psi_{b}=\mathcal{O}(\Pi^{-1}(\overline{b}))$ before $\max f|_{\tilde{\Lambda}}+\epsilon$, then for some constant $C_0>0$, that does not depend on $n$, one has  $HD(\tilde{\Lambda})\leq \frac{C_0}{n}$ provided that $n$ is large.
\end{proposition}

Let $\tilde{\Lambda}$ as in the statement of the proposition. If $\tilde{\Lambda}$ is trivial, there is nothing to prove. In other case, $\max f|_{\tilde{\Lambda}}>3$, otherwise one would have $0<HD(\tilde{\Lambda})\leq HD((\Lambda(2))_3)=d(3)=0$. Then, by propostion \ref{ferraduradelcomienzo} one can find some subhorseshoe $\bar{\Lambda}\subset (\Lambda(2))_{\max f|_{\tilde{\Lambda}}}$ such that for any $p\in \mathsf{c}(\mathscr{A})\cup \{a,b\}$ one has $\psi_p\subset\bar{\Lambda}$ and, in particular, as $\psi_b\subset\bar{\Lambda}$, by Corollary \ref{connection3}, $\tilde{\Lambda}$ does not connect with $\bar{\Lambda}$ before $\max f|_{\tilde{\Lambda}}+\epsilon$. Now, by Proposition \ref{connection11}, as $\tilde{\Lambda}\cup \bar{\Lambda} \subset (\Lambda(2))_{\max f|_{\tilde{\Lambda}}}$, we cannot have at the same time the existence of two points $x\in W^u(\tilde{\Lambda})\cap W^s(\bar{\Lambda})$ and $y\in W^u(\bar{\Lambda})\cap W^s(\tilde{\Lambda})$ such that $\mathcal{O}(x) \cup \mathcal{O}(y) \subset (\Lambda(2))_{\max f|_{\tilde{\Lambda}}+\epsilon/2}$. Without loss of generality suppose that there is no $x\in W^u(\tilde{\Lambda})\cap W^s(\bar{\Lambda})$ with $m_{\varphi,f}(x)\leq \max f|_{\tilde{\Lambda}}+\epsilon/2$ (the argument for the other case is similar). We will show that this condition forces the possible letters that may appear in the sequences that determine the unstable Cantor set of $\tilde{\Lambda}$.

Let us begin fixing $R\in \mathbb{N}$ large enough such that $R\geq 2n$ and $1/2^{R-2}<\epsilon/2$ and consider the set $\mathcal{C}_{2R+1}=\{I^u(a_0;a_1, \dots, a_{2R+1}): I^u(a_0;a_1, \dots, a_{2R+1})\cap K^u(\tilde{\Lambda})\neq \emptyset\}$, clearly $\mathcal{C}_{2R+1}$ is a covering of $K^u(\tilde{\Lambda})$. We will give a mechanism to construct coverings $\mathcal{C}_k$
with $k\geq 2R+1$ that can be used to \emph{efficiently} cover $K^u(\tilde{\Lambda})$ as $k$ goes to infinity.

Indeed, if for some $k\geq 2R+1$, and $I^u(a_0; a_1, \dots,a_{l(k)})\in \mathcal{C}_k$,  $(a_0, a_1, \dots,a_{l(k)})$ has continuations with forced first letter. That is, for every $\alpha=(\alpha_n)_{n\in \mathbb{Z}}\in \Pi(\tilde{\Lambda})$ with $\alpha_0,\alpha_1, \dots, \alpha_{l(k)}=a_0,a_1, \dots, a_{l(k)}$ one has $\alpha_{l(k)+1}=a_{l(k)+1}$ for some fixed $a_{l(k)+1}$, then we can refine the original cover $\mathcal{C}_k$, by replacing the interval $I^u(a_0;a_1, \dots, a_{l(k)})$ with the interval $I^u(a_0;a_1, \dots, a_{l(k)}, a_{l(k)+1})$. 

On the other hand, if $(a_0, a_1, \dots,a_{l(k)})$ has two continuations, said $ \gamma_{l(k)+1}=\\ (2, a_{l(k)+2},\dots)$ and $\beta_{l(k)+1}=(1, a^*_{l(k)+2},\dots)$. Take $\alpha=(\alpha_n)_{n\in \mathbb{Z}}\in \Pi(\tilde{\Lambda})$ and $\tilde{\alpha}=(\tilde{\alpha}_n)_{n\in \mathbb{Z}}\in \Pi(\tilde{\Lambda})$, such that $\alpha=(\dots,\alpha_{-2},\alpha_{-1};a_0, a_1, \dots,a_{l(k)},\gamma_{l(k)+1})$ and $\tilde{\alpha}=\\ (\dots,\tilde{\alpha}_{-2},\tilde{\alpha}_{-1};a_0, a_1, \dots,a_{l(k)},\beta_{l(k)+1})$. Fix $n\geq 68$ and consider the word $\tilde{\kappa}$ of size $3n$ just before the bifurcation. Lemma \ref{12} let us find an ordered alphabet $(\alpha,\beta)\in \mathscr{A}$ satisfying $\abs{\alpha},\abs{\beta}<n$ and $\abs{\alpha \beta}\geq n$ such that $\tilde{\kappa}$ is $(\alpha,\beta)$-semi renormalizable. Let $\kappa$ the word of biggest length in $\{\alpha,\beta\}$ and $\gamma$ the renormalization kernel, we have the following cases:

\begin{itemize}
    \item $\kappa=\alpha$ and $\kappa$ appears in $\gamma$
    
    \item $\kappa=\beta$ and $\kappa$ appears in $\gamma$
    
    \item $\kappa=\alpha$ and $\kappa$  does not appear in $\gamma$    
    \item $\kappa=\beta$ and $\kappa$ does not appear in $\gamma$
    
\end{itemize}
To deal with these cases, we will first prove some lemmas, the \emph{connection schemes}, which will let us to force at least $n/5$ letters in both continuations $\gamma_{l(k)+1}$ and $\beta_{l(k)+1}$ in any case. That is, once we know $(a_0,a_1,\dots,a_{l(k)})$ then $(2,a_{l(k)+2},\dots,a_{l(k)+n/5})$ and $(2,a^*_{l(k)+2},\dots,a^*_{l(k)+n/5})$ are determined.

\subsection{Connection schemes}
In the following subsection, in most of the cases, we will not use directly the ordered alphabet $(\alpha,\beta)$ of the previous subsection, instead, we will consider some auxiliary alphabet $(u,v)$ and express in terms of it the concatenation of the end of $(a_0,a_1,\dots, a_{l(k)})$ with the beginning of $\gamma_{l(k)+1}$ and $\beta_{l(k)+1}$ (at least up to size $n$). In this context, if $w$ is the biggest word written in the alphabet $\{u,v\}$ associated to $(u,v)$, that comes before $a_{l(k)}$, we said that one has an $(u,v)$-\emph{bifurcation} in $w$, that the continuation determined by $\gamma_{l(k)+1}$ is the $u$-\emph{continuation} and the continuation determined by $\beta_{l(k)+1}$ is the $v$-\emph{continuation}. To explain this terminology, note that by Lemma \ref{comienzo} and Corollary \ref{9}, the $u$-continuation begins with $v^-a$ and as trivially, the $v$-continuation begins with $v^-b$, then the end of   
$(a_0,a_1,\dots, a_{l(k)})$ is actually $wv^-$ in this case. 
\begin{remark}
In the same way, we can define one $(a,1)$-\emph{bifurcation} and one $(2,b)$-\emph{bifurcation}.
\end{remark}

\subsubsection{Previous lemmas}

We present here some lemmas that establish constraints on words of $\Sigma(3+6^{-3n})$ and will be used in the proof of the schemes. 

\begin{lemma}\label{bwords}
Let $(u,v)\in \mathscr{A}$ be an ordered alphabet. If $\abs{u^2v}\leq n$, then, in any word $w\in \Sigma(3+6^{-3n},3n)$ it cannot appear neither $w_0=bu^+u^2vu^3v^-a$ nor $w_1=bu^+uvu^2v^-a$ nor $w_2=bu^+uvu^2vu^2v^-a$ as subwords. The same holds for $w_3=bu^+uvuvu^2vuv^-a$ if $\abs{uv}\leq n/2$.
\end{lemma}

\begin{proof}
Observe first that $\abs{w_1}<\abs{w_0}<\abs{w_2}=\abs{u^6v^3}=3\abs{u^2v}\leq 3n$ and $\abs{w_3}=\abs{u^6v^4}<\abs{u^6v^6}=6\abs{uv}\leq 3n$. Using Proposition \ref{transpuesto} we conclude that (here the internal transposed is taken respect to the alphabet $\{u,v\}$)
$$w_0=(u(u^2vu^3)^*v)^*=(uu^3vu^2v)^*=(u^2uuvuuv)^*=(\alpha^2\beta^2)^*,$$
$$w_1=(u(uvu^2)^*v)^*=(uu^2vuv)^*=(u^2uvuv)^*=(\bar{\alpha}^2\bar{\beta}^2)^*,$$
$$w_2=(u(uvu^2vu^2)^*v)^*=(uu^2vu^2vuv)^*=(u^2uvuuvuv)^*=(\bar{\alpha}^2\bar{\beta}\bar{\alpha}\bar{\beta}^2)^*$$ 
and 
$$w_3=(u(uvuvu^2vu)^*v)^*=(uuvu^2vuvuv)^*=(uuvuuvuvuv)^*=(\tilde{\alpha}^2\tilde{\beta}^2)^*,$$
where $(\alpha,\beta)=(u,uuv)$, $(\bar{\alpha},\bar{\beta})=(u,uv)$ and $(\tilde{\alpha},\tilde{\beta})=(uuv,uv)$ are ordered alphabets. As $\Sigma(3+6^{-3n},3n)$ is closed by transpositions, the result follows from Lemma \ref{alfacuadrado}.
\end{proof}

\begin{lemma}\label{s}
Let $(u,v)\in \mathscr{A}$ be an ordered alphabet and let $\omega \in \Sigma(3+6^{-3n})$. Assume $\omega$ contains the subword $(uv)^su(uv)^{\tilde{s}}uuv$, where $s\geq 3$ and $1\leq \tilde{s}\leq s-2$. Then, $\abs{u(uv)^{\tilde{s}+1}}\geq 3(n-1).$ In particular, $\omega$ cannot contain such a word if $\abs{u(uv)^s}<3n.$
\end{lemma}

\begin{proof}
    Lemma \ref{1} applied to the cut 
    $$ab(ab)^{\tilde{s}}ab|a(ab)^{\tilde{s}}aab=ab(a(ba)^{\tilde{s}})^*b|a(a(ba)^{\tilde{s}})ab=Xb\omega^*b|a\omega a Y$$
let us see that if $W$ is the word in $\{U,V \}$ such that $u=W(a)$ and $v=W(b)$ then one has the bad cut 
\begin{eqnarray*}
b u^+ ((uv)^{\tilde{s}}u) v^- b | a u^+ (u(vu)^{\tilde{s}}) v^- a&=&b u^+ ((uv)^{\tilde{s}}u) v | u (u(vu)^{\tilde{s}}) v^- a\\ 
&\subseteq& W(ab(ab)^{\tilde{s}}ab|a(ab)^{\tilde{s}}aab) \\ &=&
uv(uv)^{\tilde{s}}uv|u(uv)^{\tilde{s}}uuv.
\end{eqnarray*}
As $\abs{u^+ (u(vu)^{\tilde{s}}) v^-}=\abs{u(uv)^{\tilde{s}+1}}-2$, it follows from Corollary \ref{mas5} that  $\abs{u(uv)^{\tilde{s}+1}}+3\geq 3n$. If $\abs{u(uv)^s}<3n$ then, 
$$3n-3\leq \abs{u(uv)^{\tilde{s}+1}}\leq \abs{u(uv)^{s-1}}=\abs{u(uv)^s}-\abs{uv}<3n-4$$
that is a contradiction because $\abs{uv}\geq 4$ for any $(u,v)\in \mathscr{A}$. 
\end{proof}

\begin{lemma}\label{rere}
Let $(u,v)\in \mathscr{A}$ be an ordered alphabet such that $\abs{uv}<n$ and let $\omega\in \Sigma(3+6^{-3n})$ which has as subword the word $(uv)^2$. The word of size $3n$ beginning in the last $uv$ of $(uv)^2$ is $(u,v)$-semi renormalizable. Then
\begin{itemize}
    \item If $\abs{u^2v}\leq n$ then the renormalization kernel cannot be $\gamma=uvu^s$, where $s>0$.
    \item If $\abs{uv}\leq \frac{3}{7}n$ then the renormalization kernel cannot be $\gamma=uvu^2vu^s$, where $s>0$.
    \end{itemize}
\end{lemma}
\begin{proof}
  For the first item
  $$\abs{u^{s-2}}=\abs{u^s}-\abs{u^2}= \abs{\gamma}-\abs{uv}-\abs{u^2}>3n-\abs{uv}-\abs{uv}-\abs{u^2}=3n-2\abs{u^2v}\geq \abs{u^2v}>\abs{uv}.$$
Similarly, for the second item
$$\abs{u^{s-2}}=\abs{u^s}-\abs{u^2}= \abs{\gamma}-\abs{uvuuv}-\abs{u^2}>3n-\abs{uv}-\abs{uvuuv}-\abs{u^2}>3n-6\abs{uv}\geq \abs{uv}.$$
Lemma \ref{9} let us conclude that $u^{s-2}$ begins with $v^-a$, but by Lemma \ref{bwords} this is an absurd because in the first case we would have the subword $bu^+uvu^2v^-a \subseteq uvuvu^2u^{s-2}$ and in the second case $bu^+uvu^2vu^2v^-a \subseteq uvuvu^2vu^2u^{s-2}$.

\end{proof}

\begin{lemma}\label{rere2}
Let $(u,v)\in \mathscr{A}$ be an ordered alphabet such that $\abs{uv}<n$ and let $\omega\in \Sigma(3+6^{-3n})$ which has as subword the word $uvu^2v$. The word of size $3n$ beginning in the last $uv$ of $uvu^2v$ is $(u,v)$-semi renormalizable. Then, if $\abs{u^2v}\leq n$, the renormalization kernel cannot be $\gamma=uvu^s$, where $s>0$.
\end{lemma}
\begin{proof}
First observe that 
\begin{eqnarray*}
      \abs{u^{s-3}}=\abs{u^s}-\abs{u^3}&\geq& \abs{\gamma}-\abs{uv}-\abs{u^3}>3n-\abs{uv}-\abs{uv}-\abs{u^3}\\ &=& 3n-2\abs{u^2v}-\abs{u}\geq \abs{u^2v}-\abs{u}=\abs{uv}.
\end{eqnarray*}
Lemma \ref{9} let us conclude that $u^{s-3}$ begins with $v^-a$, but by Lemma \ref{bwords} this is an absurd because we would have the subword $bu^+u^2vu^3v^-a \subseteq uvu^2vu^3u^{s-3}$.

\end{proof}

\subsubsection{The schemes}

The following lemma is from \cite{C1} and will be used in the proof of all schemes.
\begin{lemma}\label{lemao}
Given $T \in \mathbb{N}$, let $\beta^1,\beta^2,\beta^3\in \Sigma_2^+:=\{1,2\}^{\mathbb{N}}$ such that $[0;\beta^1]<[0;\beta^2]<[0;\beta^3]$. If for two sequences $\alpha=(\alpha_n)_{n\in \mathbb{Z}}\ \mbox{and}\ \tilde{\alpha}=(\tilde{\alpha}_n)_{n\in \mathbb{Z}}\ \mbox{in}\ \Sigma_2$ it is true that  $\alpha_{0},\dots ,\alpha_{2T+1}=\tilde{\alpha}_{0},\dots ,\tilde{\alpha}_{2T+1}$, then, for all $j\leq 2T+1$ we have
\begin{eqnarray*}
\lambda(\sigma^j(\dots,\alpha_{-2},\alpha_{-1};\alpha_{0},\dots ,\alpha_{2T+1},\beta^2))< \max \{m(\dots,\alpha_{-2},\alpha_{-1};\alpha_{0},\dots ,\alpha_{2T+1},\beta^1), \\ m(\dots ,\tilde{\alpha}_{-2},\tilde{\alpha}_{-1};\tilde{\alpha}_{0},\dots ,\tilde{\alpha}_{2T+1},\beta^3)\}+1/2^{T-1}.\end{eqnarray*}   
\end{lemma}

\begin{lemma}[Scheme 1]\label{esquema1}
Let $(u,v)\in \mathscr{A}$ be an ordered alphabet. Suppose one has an $(u,v)$-bifurcation in $u(uv)^iu$ where $i\geq 2$ and $\abs{u^2v}<n$. Then, the $u$-continuation begins with $uvu$ and the $v$-continuation with $vuuvu$. Additionally, if $\abs{uv}\leq \frac{3}{7}n$, then the $u$-continuation begins with $uvuv$ and the $v$-continuation with $vuuvuv$ and if $\abs{(uv)^i}\leq 3n$ the $u$-continuation begins with $(uv)^i$.
\end{lemma}

\begin{proof}


Suppose we have an $(u,v)$-bifurcation in $u(uv)^iu$ with $\abs{u^2v}<n$. By Lemma \ref{rere} the continuation that begins with $u$, begins with $u^{\tilde{j}}v$ for some $\tilde{j}\geq 1$. If $\tilde{j}>1$ then, as $u^{\tilde{j}} v$ and $uv$ begin with $v^-a$ we would have 
$$[v\dots]<[(uv)^i\overline{u(uv)^i}]<[uu^{\tilde{j}-1}v\dots ]$$
and therefore by Lemma \ref{lemao} for all $j\leq l(k)$
\begin{eqnarray*}
\lambda(\sigma^j(\dots,\tilde{\alpha}_{-1};\tilde{\alpha}_{0},\dots ,\tilde{\alpha}_{l(k)},(uv)^i\overline{u(uv)^i}))&\leq&\max \{m(\dots, \alpha_{-1};\alpha_{0},\dots ,\alpha_{l(k)},\gamma_{l(k)+1}),\\&&  m(\dots ,\tilde{\alpha}_{-1};\tilde{\alpha}_{0},\dots ,\tilde{\alpha}_{l(k)},\beta_{l(k)+1})\}+\frac{1}{2^{R-1}} \\ &<& \max f|_{\tilde{\Lambda}}+\epsilon/2.
\end{eqnarray*}
Proposition \ref{control} applied to $XRY=ababab$ where $R=ab$ implies that for any alphabet $(\alpha,\beta)$, $\lambda$ has value less than $3$ in any position of the $\alpha\beta$ in the middle of $\alpha\beta\alpha\beta\alpha\beta$. 
Then, for $j\geq l(k)+1$, as $(u(uv)^{i-1},uv)$ is an alphabet we conclude that 
\begin{eqnarray*}
\lambda(\sigma^j(\dots,\tilde{\alpha}_{-2},\tilde{\alpha}_{-1};\tilde{\alpha}_{0},\dots ,\tilde{\alpha}_{l(k)-\abs{v^-}},(uv)^i\overline{u(uv)^i}))&=&\\
\lambda(\sigma^j(\dots,\tilde{\alpha}_{-1};\tilde{\alpha}_{0},\dots,u(uv)^iu,(uv)^i\overline{u(uv)^i}))&<&3< \max f|_{\tilde{\Lambda}}+\epsilon/2.
\end{eqnarray*}
Then taking $x=\Pi^{-1}((\dots,\tilde{\alpha}_{-2},\tilde{\alpha}_{-1};\tilde{\alpha}_{0},\dots ,\tilde{\alpha}_{l(k)-\abs{v^-}},(uv)^i\overline{u(uv)^i}))$ one would have
$$x\in W^u(\tilde{\Lambda})\cap W^s(\psi_{u(uv)^i})\subset W^u(\tilde{\Lambda})\cap W^s(\bar{\Lambda})\ \mbox{and}\ m_{\varphi,f}(x)\leq \max f|_{\tilde{\Lambda}}+\epsilon/2$$
that is a contradiction. We conclude that $\tilde{j}=1$ in this case. Now, the letter that follows $uv$ in the $u$-continuation cannot be $v$ because in other case we would have $u^2v^2$ and $\abs{u^2v^2}=2\abs{uv}<2n$ contradicts Lemma \ref{alfacuadrado}. 

If, additionally $\abs{uv}\leq \frac{3}{7}n$, again by Lemma \ref{rere}, for some $\tilde{j}\geq 1$ the $u$-continuation begins with $uvu^{\tilde{j}}v$ and if $\tilde{j}>1$, using that 
$$[v\dots]<[(uv)^i\overline{u(uv)^i}]<[uvuu^{\tilde{j}-1}v\dots ],$$
we conclude again that 
$$\lambda(\sigma^j(\dots,\tilde{\alpha}_{-2},\tilde{\alpha}_{-1};\tilde{\alpha}_{0},\dots ,\tilde{\alpha}_{l(k)-\abs{v^-}},(uv)^i\overline{u(uv)^i}))< \max f|_{\tilde{\Lambda}}+\epsilon/2$$
and that $\tilde{j}=1$. Then, the $u$-continuation begins with $(uv)^2$ in this case. If $\abs{(uv)^i}\leq 3n$, continuing in this way (here we need to use Lemma \ref{rere} again), we can force $(uv)^i$ since after $(uv)^r$ with $2\leq r\leq i-1$ always follows $u$ because in other case we would have $u(uv)^rv=u^2(vu)^{r-1}v^2$ and $\abs{u(uv)^rv}\leq \abs{(uv)^i}\leq 3n$
which contradicts Lemma \ref{alfacuadrado}.

For the $v$-continuation, if after $v$ we have other $v$, then we would have 
$$[vv\dots]<[\overline{vu}]<[(uv)^i\dots ]$$
and therefore, as before, by Lemma \ref{lemao} for all $j\leq l(k)$
\begin{eqnarray*}
\lambda(\sigma^j(\dots,\tilde{\alpha}_{-2},\tilde{\alpha}_{-1};\tilde{\alpha}_{0},\dots ,\tilde{\alpha}_{l(k)-\abs{v^-}},\overline{vu}))< \max f|_{\tilde{\Lambda}}+\epsilon/2.
\end{eqnarray*}
As $(u,v)$ is an alphabet, as before, $\lambda$ has value less than $3$ in any position of the $uv$ in the middle of $uvuvuv$.
Then, for $j\geq l(k)+1$ we conclude that 
\begin{eqnarray*}
\lambda(\sigma^j(\dots,\tilde{\alpha}_{-2},\tilde{\alpha}_{-1};\tilde{\alpha}_{0},\dots ,\tilde{\alpha}_{l(k)-\abs{v^-}},\overline{vu}))&=&
\lambda(\sigma^j(\dots,\tilde{\alpha}_{-1};\tilde{\alpha}_{0},\dots,u(uv)^iu,\overline{vu}))\\ &<&3< \max f|_{\tilde{\Lambda}}+\epsilon/2.
\end{eqnarray*}
Then taking $x=\Pi^{-1}((\dots,\tilde{\alpha}_{-2},\tilde{\alpha}_{-1};\tilde{\alpha}_{0},\dots ,\tilde{\alpha}_{l(k)-\abs{v^-}},\overline{vu}))$ one would have
$$x\in W^u(\tilde{\Lambda})\cap W^s(\psi_{uv})\subset W^u(\tilde{\Lambda})\cap W^s(\bar{\Lambda})\ \mbox{and}\ m_{\varphi,f}(x)\leq \max f|_{\tilde{\Lambda}}+\epsilon/2$$
that is a contradiction. We conclude that the $v$-continuation begins with $vu$.

Now, if after $vu$ it follows one $v$, then we would have 
$$[vuv\dots]<[v\overline{u(uv)^i}]<[(uv)^i\dots ]$$
and therefore, as before, by Lemma \ref{lemao} for all $j\leq l(k)$
\begin{eqnarray*}
\lambda(\sigma^j(\dots,\tilde{\alpha}_{-2},\tilde{\alpha}_{-1};\tilde{\alpha}_{0},\dots ,\tilde{\alpha}_{l(k)-\abs{v^-}},v\overline{u(uv)^{i+1}}))< \max f|_{\tilde{\Lambda}}+\epsilon/2.
\end{eqnarray*}

As $(u(uv)^i,uv)$ is an alphabet, $\lambda$ has value less than $3$ in any position of the $u(uv)^i$ in the middle of $u(uv)^{i+1}u(uv)^{i+1}u(uv)^{i+1}$. Additionally, as both cuts $a(ab)^iab|a(ab)^i=a(ab)^iab|aa(ba)^ib$ and $Xba|baY=abba|baab$ are good, corollaries \ref{3} and \ref{4} imply that we have the good cuts $u(uv)^iuv|u(uv)^i$ and $b u^+ v^- a | b u^+ v^- a\subseteq uvvuvuuv$. Then, for $j\geq l(k)+1$ we conclude that 
\begin{eqnarray*}
\lambda(\sigma^j(\dots,\tilde{\alpha}_{-2},\tilde{\alpha}_{-1};\tilde{\alpha}_{0},\dots ,\tilde{\alpha}_{l(k)-\abs{v^-}},v\overline{u(uv)^{i+1}}))&=&\\
\lambda(\sigma^j(\dots,\tilde{\alpha}_{-1};\tilde{\alpha}_{0},\dots,u(uv)^iu,v\overline{u(uv)^{i+1}}))&<&3< \max f|_{\tilde{\Lambda}}+\epsilon/2.
\end{eqnarray*}
Then taking $x=\Pi^{-1}((\dots,\tilde{\alpha}_{-2},\tilde{\alpha}_{-1};\tilde{\alpha}_{0},\dots ,\tilde{\alpha}_{l(k)-\abs{v^-}},v\overline{u(uv)^{i+1}}))$ one would have the contradiction 
$$x\in W^u(\tilde{\Lambda})\cap W^s(\psi_{u(uv)^{i+1}})\subset W^u(\tilde{\Lambda})\cap W^s(\bar{\Lambda})\ \mbox{and}\ m_{\varphi,f}(x)\leq \max f|_{\tilde{\Lambda}}+\epsilon/2$$
We conclude that the $v$-continuation begins with $vuu$.

Lemma \ref{rere} let us conclude that for some $\tilde{j}\geq 0$ the $v$-continuation begins with $vuuu^{\tilde{j}}v$. But if $\tilde{j}\neq 0$ then one would have $bu^+uvu^2v^-a\subseteq uvuvu^2u^{\tilde{j}}v$ that contradicts Lemma \ref{bwords}. Then the $v$-continuation begins with $vuuvu$ (the last $u $ to avoid the sequence $u^2v^2$). If, additionally $\abs{uv}\leq \frac{3}{7}n$, Lemma \ref{rere} let us conclude again that for some $\tilde{j}\geq 1$ the $v$-continuation begins with $vuuvu^{\tilde{j}}v$. But if $\tilde{j}>1$ then one would have $bu^+uvuvu^2vuv^-a\subseteq uvuvuvu^2vu^{\tilde{j}}v$ that contradicts Lemma \ref{bwords}. This concludes the proof of the lemma.   
\end{proof}

\begin{lemma}[Scheme 2]\label{muchos2}
Suppose one has a $(2,b)$-bifurcation in $b(2)^i$ where $i\geq 2$. Then, the $2$-continuation begins with $ab$. 
\end{lemma}
\begin{proof}
If the $2$-continuation begins with $b$ then we would have 
$$[b\dots]<[\overline{a}]<[2b\dots ]$$
and therefore, as before, by Lemma \ref{lemao} for all $j\leq l(k)$
\begin{eqnarray*}
\lambda(\sigma^j(\dots,\tilde{\alpha}_{-2},\tilde{\alpha}_{-1};\tilde{\alpha}_{0},\dots ,\tilde{\alpha}_{l(k)},\overline{a}))< \max f|_{\tilde{\Lambda}}+\epsilon/2.
\end{eqnarray*}
As $\lambda$ has value less than $3$ in any $2$ in the middle of $222$ then, for $j\geq l(k)+1$ we conclude that 
\begin{eqnarray*}
\lambda(\sigma^j(\dots,\tilde{\alpha}_{-2},\tilde{\alpha}_{-1};\tilde{\alpha}_{0},\dots ,\tilde{\alpha}_{l(k)},\overline{a}))< \max f|_{\tilde{\Lambda}}+\epsilon/2.
\end{eqnarray*}
Then taking $x=\Pi^{-1}((\dots,\tilde{\alpha}_{-2},\tilde{\alpha}_{-1};\tilde{\alpha}_{0},\dots ,\tilde{\alpha}_{l(k)},\overline{a}))$ one would have
$$x\in W^u(\tilde{\Lambda})\cap W^s(\psi_{a})\subset W^u(\tilde{\Lambda})\cap W^s(\bar{\Lambda})\ \mbox{and}\ m_{\varphi,f}(x)\leq \max f|_{\tilde{\Lambda}}+\epsilon/2$$
that is a contradiction. We conclude that the $2$-continuation begins with $a$.

Now suppose that $(2)^i=a^s$ for some $s\geq 1$. If after the $a$ in the $2$-continuation there is other $2$, then we would have  
$$[b\dots]<[ab\overline{a^{s+1}b}]<[a2\dots ]$$
and therefore, as before, by Lemma \ref{lemao} for all $j\leq l(k)$
\begin{eqnarray*}
\lambda(\sigma^j(\dots,\tilde{\alpha}_{-2},\tilde{\alpha}_{-1};\tilde{\alpha}_{0},\dots ,\tilde{\alpha}_{l(k)},ab\overline{a^{s+1}b}))< \max f|_{\tilde{\Lambda}}+\epsilon/2.
\end{eqnarray*}
As $(a,a^sb)$ is an alphabet, $\lambda$ has value less than $3$ in any position of the $ab^{s+1}$ in the middle of $a^{s+1}ba^{s+1}ba^{s+1}b$ and in any $2$ in the middle of $222$. Note also that the cut $ba^sa|ba^sa$ is good. Then, for $j\geq l(k)+1$ we conclude that 
\begin{eqnarray*}
\lambda(\sigma^j(\dots,\tilde{\alpha}_{-2},\tilde{\alpha}_{-1};\tilde{\alpha}_{0},\dots ,\tilde{\alpha}_{l(k)},ab\overline{a^{s+1}b}))&=&
\lambda(\sigma^j(\dots,\tilde{\alpha}_{-1};\tilde{\alpha}_{0},\dots,ba^s,ab\overline{a^{s+1}b}))\\&<& \max f|_{\tilde{\Lambda}}+\epsilon/2.
\end{eqnarray*}
Then taking $x=\Pi^{-1}((\dots,\tilde{\alpha}_{-2},\tilde{\alpha}_{-1};\tilde{\alpha}_{0},\dots ,\tilde{\alpha}_{l(k)},ab\overline{a^{s+1}b}))$ one would have
$$x\in W^u(\tilde{\Lambda})\cap W^s(\psi_{a^{s+1}b})\subset W^u(\tilde{\Lambda})\cap W^s(\bar{\Lambda})\ \mbox{and}\ m_{\varphi,f}(x)\leq \max f|_{\tilde{\Lambda}}+\epsilon/2$$
that is a contradiction. We conclude that the $2$-continuation begins with $ab$ in this case.  

Suppose that $(2)^i=a^s2$ for some $s\geq 1$. If the $2$-continuation begins with $a2$ and the $b$-continuation begins with $ba^rb$ where $r\leq s-1$ then we would have
$$[ba^rb\dots]<[ab\overline{a^sb}]<[a2\dots ]$$
and therefore, as before, by Lemma \ref{lemao} for all $j\leq l(k)$
\begin{eqnarray*}
\lambda(\sigma^j(\dots,\tilde{\alpha}_{-2},\tilde{\alpha}_{-1};\tilde{\alpha}_{0},\dots ,\tilde{\alpha}_{l(k)},ab\overline{a^sb}))< \max f|_{\tilde{\Lambda}}+\epsilon/2.
\end{eqnarray*}
As $(a,a^{s-1}b)$ is an alphabet, $\lambda$ has value less than $3$ in any position of the $a^sb$ in the middle of $a^sba^sba^sb$ and the same holds for any $2$ in the middle of $222$. Further, as the cut $aa^{s-1}b|aa^{s-1}b$ is good, then, for $j\geq l(k)+1$ and $j\neq l(k)+2 $ we conclude that
\begin{eqnarray*}
\lambda(\sigma^j(\dots,\tilde{\alpha}_{-2},\tilde{\alpha}_{-1};\tilde{\alpha}_{0},\dots ,\tilde{\alpha}_{l(k)},ab\overline{a^sb}))&=&\\
\lambda(\sigma^j(\dots,\tilde{\alpha}_{-1};\tilde{\alpha}_{0},\dots,ba^s2,ab\overline{a^sb}))&<&3< \max f|_{\tilde{\Lambda}}+\epsilon/2.
\end{eqnarray*}

For $j=l(k)+2$, we have the bad cut $2a^sa|ba^sa$ which can be compared with the bad cut $2a^ra|ba^rb$ in the $b$-continuation, then by Lemma \ref{compare} one has 
$$\lambda(\sigma^{l(k)+2}(\dots,\tilde{\alpha}_{-2},\tilde{\alpha}_{-1};\tilde{\alpha}_{0},\dots ,\tilde{\alpha}_{l(k)},ab\overline{a^sb}))<\max f|_{\tilde{\Lambda}}+\epsilon/2.$$
Then taking $x=\Pi^{-1}((\dots,\tilde{\alpha}_{-2},\tilde{\alpha}_{-1};\tilde{\alpha}_{0},\dots ,\tilde{\alpha}_{l(k)},ab\overline{a^sb}))$ one would have
$$x\in W^u(\tilde{\Lambda})\cap W^s(\psi_{a^sb})\subset W^u(\tilde{\Lambda})\cap W^s(\bar{\Lambda})\ \mbox{and}\ m_{\varphi,f}(x)\leq \max f|_{\tilde{\Lambda}}+\epsilon/2$$
that is a contradiction. We conclude that the $2$-continuation begins with $ab$ in this case. 

If the $2$-continuation begins with $a2$ and the $b$-continuation begins with $ba^rb$ or $ba^r2b$ where $r\geq s$ then we would have
$$[ba^r\dots]<[ab\overline{a^{s+1}b}]<[a2\dots ]$$
and therefore, by Lemma \ref{lemao} for all $j\leq l(k)$
\begin{eqnarray*}
\lambda(\sigma^j(\dots,\tilde{\alpha}_{-2},\tilde{\alpha}_{-1};\tilde{\alpha}_{0},\dots ,\tilde{\alpha}_{l(k)},ab\overline{a^{s+1}b}))< \max f|_{\tilde{\Lambda}}+\epsilon/2.
\end{eqnarray*}
As $(a,a^sb)$ is an alphabet, $\lambda$ has value less than $3$ in any position of the $a^{s+1}b$ in the middle of $a^{s+1}ba^{s+1}ba^{s+1}b$ and the same holds for any $2$ in the middle of $222$. Further, as the cut $aa^sb|aa^sb$ is good, then, for $j\geq l(k)+1$ and $j\neq l(k)+2 $ we conclude that
\begin{eqnarray*}
\lambda(\sigma^j(\dots,\tilde{\alpha}_{-2},\tilde{\alpha}_{-1};\tilde{\alpha}_{0},\dots ,\tilde{\alpha}_{l(k)},ab\overline{a^{s+1}b}))&=&\\
\lambda(\sigma^j(\dots,\tilde{\alpha}_{-1};\tilde{\alpha}_{0},\dots,ba^s2,ab\overline{a^{s+1}b}))&<&3< \max f|_{\tilde{\Lambda}}+\epsilon/2.
\end{eqnarray*}

For $j=l(k)+2$, we have the bad cut $b2a^sa|ba^s22$ which can be compared with the bad cut $b2a^{s-1}a|ba^{s-1}22$ in the $b$-continuation, then by Lemma \ref{compare} one has 
$$\lambda(\sigma^{l(k)+2}(\dots,\tilde{\alpha}_{-2},\tilde{\alpha}_{-1};\tilde{\alpha}_{0},\dots ,\tilde{\alpha}_{l(k)},ab\overline{a^{s+1}b}))<\max f|_{\tilde{\Lambda}}+\epsilon/2.$$
Then taking $x=\Pi^{-1}((\dots,\tilde{\alpha}_{-2},\tilde{\alpha}_{-1};\tilde{\alpha}_{0},\dots ,\tilde{\alpha}_{l(k)},ab\overline{a^{s+1}b}))$ one would have
$$x\in W^u(\tilde{\Lambda})\cap W^s(\psi_{a^{s+1}b})\subset W^u(\tilde{\Lambda})\cap W^s(\bar{\Lambda})\ \mbox{and}\ m_{\varphi,f}(x)\leq \max f|_{\tilde{\Lambda}}+\epsilon/2$$
that is a contradiction. We conclude that the $2$-continuation begins with $ab$ in these cases. 

Finally, if the $2$-continuation begins with $a2$ and the $b$-continuation begins with $ba^r2b$ where $r\geq s-1$ then we would have
$$[ba^r2b\dots]<[aba^{s+1}2ba^{s+2}b\overline{a^{s+2}b}]<[a2\dots ]$$
and therefore, by Lemma \ref{lemao} for all $j\leq l(k)$
\begin{eqnarray*}
\lambda(\sigma^j(\dots,\tilde{\alpha}_{-2},\tilde{\alpha}_{-1};\tilde{\alpha}_{0},\dots ,\tilde{\alpha}_{l(k)},aba^{s+1}2ba^{s+2}b\overline{a^{s+2}b}))< \max f|_{\tilde{\Lambda}}+\epsilon/2.
\end{eqnarray*}
As $(a,a^{s+1}b)$ is an alphabet, $\lambda$ has value less than $3$ in any position of the $a^{s+2}b$ in the middle of $a^{s+2}ba^{s+2}ba^{s+2}b$ and the same holds for any 2 in the middle of $222$. Further, as the cuts $2a^{s+1}b|aa^{s+1}b$ and $ba^{s+1}a|ba^{s+1}a$ are good, then, for $j\geq l(k)+1$ and $j\neq l(k)+2, l(k)+5, l(k)+2(s+4)+2 $ we conclude that
\begin{eqnarray*}
\lambda(\sigma^j(\dots,\tilde{\alpha}_{-2},\tilde{\alpha}_{-1};\tilde{\alpha}_{0},\dots ,\tilde{\alpha}_{l(k)},aba^{s+1}2ba^{s+2}b\overline{a^{s+2}b}))&=&\\
\lambda(\sigma^j(\dots,\tilde{\alpha}_{-1};\tilde{\alpha}_{0},\dots,ba^s2,aba^{s+1}2ba^{s+2}b\overline{a^{s+2}b}))&<&3< \max f|_{\tilde{\Lambda}}+\epsilon/2.
\end{eqnarray*}

For $j= l(k)+3, l(k)+4,l(k)+2s+10 $, we have the bad cuts $b2a^sa|ba^s22$ (two times) and $22a^sb|aa^s2b$ (the transposed of the previous cut) which can be compared with the bad cut $22a^{r-1}b|aa^{r-1}2b$ in the $b$-continuation, then as $\Sigma(t+\epsilon)$ is closed by transpositions, by Lemma \ref{compare} one has 
$$\lambda(\sigma^j(\dots,\tilde{\alpha}_{-2},\tilde{\alpha}_{-1};\tilde{\alpha}_{0},\dots ,\tilde{\alpha}_{l(k)},aba^{s+1}2ba^{s+2}b\overline{a^{s+2}b}))<\max f|_{\tilde{\Lambda}}+\epsilon/2.$$
Then taking $x=\Pi^{-1}((\dots,\tilde{\alpha}_{-2},\tilde{\alpha}_{-1};\tilde{\alpha}_{0},\dots ,\tilde{\alpha}_{l(k)},aba^{s+1}2ba^{s+2}b\overline{a^{s+2}b}))$ one would have
$$x\in W^u(\tilde{\Lambda})\cap W^s(\psi_{a^{s+2}b})\subset W^u(\tilde{\Lambda})\cap W^s(\bar{\Lambda})\ \mbox{and}\ m_{\varphi,f}(x)\leq \max f|_{\tilde{\Lambda}}+\epsilon/2$$
that is a contradiction. We conclude that the $2$-continuation begins with $ab$ in these cases too. This finishes the proof of the lemma. 
\end{proof}

\begin{lemma}[Scheme 3]\label{esquema2}

Let $(u,v)\in \mathscr{A}$ be an ordered alphabet. Suppose one has an $(u,v)$-bifurcation in $vv(uv)^i$ where $i\geq 2$ and $\abs{u^2v}\leq n$. Then, the $u$-continuation begins with $uvv$. 
\end{lemma}

\begin{proof}
By Lemma \ref{rere} the continuation that begins with $u$, begins with $u^{\tilde{j}}v$ for some $\tilde{j}\geq 1$. If $\tilde{j}>1$ then, as $u^{\tilde{j}} v$ and $uv$ begin with $v^-a$ we would have 
$$[v\dots]<[\overline{uv}]<[uu^{\tilde{j}-1}v\dots ]$$
and therefore, as before, by Lemma \ref{lemao} for all $j\leq l(k)$
\begin{eqnarray*}
\lambda(\sigma^j(\dots,\tilde{\alpha}_{-2},\tilde{\alpha}_{-1};\tilde{\alpha}_{0},\dots ,\tilde{\alpha}_{l(k)-\abs{v^-}},\overline{uv}))< \max f|_{\tilde{\Lambda}}+\epsilon/2.
\end{eqnarray*}
As $(u,v)$ is an alphabet, as before, $\lambda$ has value less than $3$ in any position of the $uv$ in the middle of $uvuvuv$.
Then, for $j\geq l(k)+1$ we conclude that 
\begin{eqnarray*}
\lambda(\sigma^j(\dots,\tilde{\alpha}_{-2},\tilde{\alpha}_{-1};\tilde{\alpha}_{0},\dots ,\tilde{\alpha}_{l(k)-\abs{v^-}},\overline{uv}))&=&
\lambda(\sigma^j(\dots,\tilde{\alpha}_{-1};\tilde{\alpha}_{0},\dots,vv(uv)^i,\overline{uv}))\\ &<&3<\max f|_{\tilde{\Lambda}}+\epsilon/2.
\end{eqnarray*}
Then taking $x=\Pi^{-1}((\dots,\tilde{\alpha}_{-2},\tilde{\alpha}_{-1};\tilde{\alpha}_{0},\dots ,\tilde{\alpha}_{l(k)-\abs{v^-}},\overline{uv}))$ one would have
$$x\in W^u(\tilde{\Lambda})\cap W^s(\psi_{uv})\subset W^u(\tilde{\Lambda})\cap W^s(\bar{\Lambda})\ \mbox{and}\ m_{\varphi,f}(x)\leq \max f|_{\tilde{\Lambda}}+\epsilon/2$$
that is a contradiction. We conclude that the $u$-continuation begins with $uv$.

Again, by Lemma \ref{rere} the $u$-continuation begins with $uvu^{\tilde{j}}v$ for some $\tilde{j}\geq 0$. If $\tilde{j}>0$ then, we would have 
$$[v\dots]<[uv\overline{v(uv)^{i+1}}]<[u^{\tilde{j}}v\dots ]$$
and therefore by Lemma \ref{lemao} for all $j\leq l(k)$
\begin{eqnarray*}
\lambda(\sigma^j(\dots,\tilde{\alpha}_{-1};\tilde{\alpha}_{0},\dots ,\tilde{\alpha}_{l(k)-\abs{v^-}},uv\overline{v(uv)^{i+1}}))&\leq&\max \{m(\dots, \alpha_{-1};\alpha_{0},\dots ,\alpha_{k},\gamma_{k+1}),\\&&  m(\dots ,\tilde{\alpha}_{-1};\tilde{\alpha}_{0},\dots ,\tilde{\alpha}_{k},\beta_{k+1})\}+\frac{1}{2^{R-1}}\\ &<& \max f|_{\tilde{\Lambda}}+\epsilon/2.
\end{eqnarray*}
As $(uv,(uv)^iv)$ is an alphabet, $\lambda$ has value less than $3$ in any position of the $(uv)^{i+1}v$ in the middle of $(uv)^{i+1}v(uv)^{i+1}v(uv)^{i+1}v$. Additionally, as both cuts $(ab)^iabb|ab(ab)^ib=a(ba)^ibb|ab(ab)^ib$ and   $abb(ab)^ia|bb(ab)^{i+1}=ab(ba)^iba|bb(ab)^iab$ are good, corollaries \ref{3} and \ref{4} imply that we have the good cuts $u(vu)^ivv|uv(uv)^iv$ and $b u^+(vu)^ivv^- a |b \\  u^+v(uv)^iv^- a\subseteq   uv(vu)^ivuvv(uv)^iuv$. Then, for $j\geq l(k)+1$ we conclude that 
\begin{eqnarray*}
\lambda(\sigma^j(\dots,\tilde{\alpha}_{-2},\tilde{\alpha}_{-1};\tilde{\alpha}_{0},\dots ,\tilde{\alpha}_{l(k)-\abs{v^-}},uv\overline{v(uv)^{i+1}}))&=&\\
\lambda(\sigma^j(\dots,\tilde{\alpha}_{-1};\tilde{\alpha}_{0},\dots,vv(uv)^i,uv\overline{v(uv)^{i+1}}))&<&3< \max f|_{\tilde{\Lambda}}+\epsilon/2.
\end{eqnarray*}
Then taking $x=\Pi^{-1}((\dots,\tilde{\alpha}_{-2},\tilde{\alpha}_{-1};\tilde{\alpha}_{0},\dots ,\tilde{\alpha}_{l(k)-\abs{v^-}},uv\overline{v(uv)^{i+1}}))$ one would have
$$x\in W^u(\tilde{\Lambda})\cap W^s(\psi_{(uv)^{i+1}v})\subset W^u(\tilde{\Lambda})\cap W^s(\bar{\Lambda})\ \mbox{and}\ m_{\varphi,f}(x)\leq \max f|_{\tilde{\Lambda}}+\epsilon/2$$
that is a contradiction. We conclude that $\tilde{j}=0$ and then that the $u$-continuation begins with $uvv$.
    
\end{proof}

\begin{lemma}[Scheme 4]\label{muchos1}

Suppose one has an $(a,1)$-bifurcation in $a(1)^i$ where $i\geq 1$. Then, the $1$-continuation begins with $ba$. 
\end{lemma}
\begin{proof}
If the $1$-continuation begins with $a$ then we would have 
$$[1a\dots]<[\overline{b}]<[a\dots ]$$
and therefore, as before, by Lemma \ref{lemao} for all $j\leq l(k)$
\begin{eqnarray*}
\lambda(\sigma^j(\dots,\tilde{\alpha}_{-2},\tilde{\alpha}_{-1};\tilde{\alpha}_{0},\dots ,\tilde{\alpha}_{l(k)},\overline{b}))< \max f|_{\tilde{\Lambda}}+\epsilon/2.
\end{eqnarray*}
As $\lambda$ has value less than $3$ in any position of any $b$ then, for $j\geq l(k)+1$ we conclude that 
\begin{eqnarray*}
\lambda(\sigma^j(\dots,\tilde{\alpha}_{-2},\tilde{\alpha}_{-1};\tilde{\alpha}_{0},\dots ,\tilde{\alpha}_{l(k)},\overline{b}))< \max f|_{\tilde{\Lambda}}+\epsilon/2.
\end{eqnarray*}
Then taking $x=\Pi^{-1}((\dots,\tilde{\alpha}_{-2},\tilde{\alpha}_{-1};\tilde{\alpha}_{0},\dots ,\tilde{\alpha}_{l(k)},\overline{b}))$ one would have
$$x\in W^u(\tilde{\Lambda})\cap W^s(\psi_{b})\subset W^u(\tilde{\Lambda})\cap W^s(\bar{\Lambda})\ \mbox{and}\ m_{\varphi,f}(x)\leq \max f|_{\tilde{\Lambda}}+\epsilon/2$$
that is a contradiction. We conclude that the $1$-continuation begins with $b$.

Now suppose that $(1)^i=b^s$ for some $s\geq 1$. If after the $b$ in the $1$-continuation there is other $1$, then we would have  
$$[b1\dots]<[b\overline{ab^{s+1}}]<[a\dots ]$$
and therefore, as before, by Lemma \ref{lemao} for all $j\leq l(k)$
\begin{eqnarray*}
\lambda(\sigma^j(\dots,\tilde{\alpha}_{-2},\tilde{\alpha}_{-1};\tilde{\alpha}_{0},\dots ,\tilde{\alpha}_{l(k)},b\overline{ab^{s+1}}))< \max f|_{\tilde{\Lambda}}+\epsilon/2.
\end{eqnarray*}
As $(ab^s,b)$ is an alphabet, $\lambda$ has value less than $3$ in any position of the $ab^{s+1}$ in the middle of $ab^{s+1}ab^{s+1}ab^{s+1}$ and the same holds for any position of any $b$.
Then, for $j\geq l(k)+1$ we conclude that 
\begin{eqnarray*}
\lambda(\sigma^j(\dots,\tilde{\alpha}_{-1};\tilde{\alpha}_{0},\dots ,\tilde{\alpha}_{l(k)},b\overline{ab^{s+1}}))=
\lambda(\sigma^j(\dots,\tilde{\alpha}_{-1};\tilde{\alpha}_{0},\dots,ab^s,b\overline{ab^{s+1}}))< \max f|_{\tilde{\Lambda}}+\frac{\epsilon}{2}.
\end{eqnarray*}
Then taking $x=\Pi^{-1}((\dots,\tilde{\alpha}_{-2},\tilde{\alpha}_{-1};\tilde{\alpha}_{0},\dots ,\tilde{\alpha}_{l(k)},b\overline{ab^{s+1}}))$ one would have
$$x\in W^u(\tilde{\Lambda})\cap W^s(\psi_{ab^{s+1}})\subset W^u(\tilde{\Lambda})\cap W^s(\bar{\Lambda})\ \mbox{and}\ m_{\varphi,f}(x)\leq \max f|_{\tilde{\Lambda}}+\epsilon/2$$
that is a contradiction. We conclude that the $1$-continuation begins with $ba$ in this case.  

Suppose that $(1)^i=b^s1$ for some $s\geq 1$. If the $1$-continuation begins with $b1$ and the $a$-continuation begins with $ab^ra$ where $r\leq s-1$ then we would have
$$[b1\dots]<[b\overline{ab^s}]<[ab^ra\dots ]$$
and therefore, as before, by Lemma \ref{lemao} for all $j\leq l(k)$
\begin{eqnarray*}
\lambda(\sigma^j(\dots,\tilde{\alpha}_{-2},\tilde{\alpha}_{-1};\tilde{\alpha}_{0},\dots ,\tilde{\alpha}_{l(k)},b\overline{ab^s}))< \max f|_{\tilde{\Lambda}}+\epsilon/2.
\end{eqnarray*}
As $(ab^{s-1},b)$ is an alphabet, $\lambda$ has value less than $3$ in any position of the $ab^s$ in the middle of $ab^sab^sab^s$ and the same holds for any position of any $b$. Further, as the cut $bb^{s-1}a|bb^{s-1}a$ is good, then, for $j\geq l(k)+1$ and $j\neq l(k)+3 $ we conclude that
\begin{eqnarray*}
\lambda(\sigma^j(\dots,\tilde{\alpha}_{-2},\tilde{\alpha}_{-1};\tilde{\alpha}_{0},\dots ,\tilde{\alpha}_{l(k)},b\overline{ab^s}))&=&\\
\lambda(\sigma^j(\dots,\tilde{\alpha}_{-1};\tilde{\alpha}_{0},\dots,ab^s1,b\overline{ab^s}))&<&3< \max f|_{\tilde{\Lambda}}+\epsilon/2.
\end{eqnarray*}

For $j=l(k)+3$, we have the bad cut $1b^sb|ab^sa$ which can be compared with the bad cut  $1b^rb|ab^ra$ in the $a$-continuation, then by Lemma \ref{compare} one has 
$$\lambda(\sigma^{l(k)+3}(\dots,\tilde{\alpha}_{-2},\tilde{\alpha}_{-1};\tilde{\alpha}_{0},\dots ,\tilde{\alpha}_{k},b\overline{ab^s}))<\max f|_{\tilde{\Lambda}}+\epsilon/2.$$
Then taking $x=\Pi^{-1}((\dots,\tilde{\alpha}_{-2},\tilde{\alpha}_{-1};\tilde{\alpha}_{0},\dots ,\tilde{\alpha}_{l(k)},b\overline{ab^s}))$ one would have
$$x\in W^u(\tilde{\Lambda})\cap W^s(\psi_{ab^s})\subset W^u(\tilde{\Lambda})\cap W^s(\bar{\Lambda})\ \mbox{and}\ m_{\varphi,f}(x)\leq \max f|_{\tilde{\Lambda}}+\epsilon/2$$
that is a contradiction. We conclude that the $1$-continuation begins with $ba$ in this case. 

If the $1$-continuation begins with $b1$ and the $a$-continuation begins with $ab^ra$ or $ab^r1a$ where $r\geq s$ then we would have
$$[b1\dots]<[b\overline{ab^{s+1}}]<[ab^r\dots ]$$
and therefore, by Lemma \ref{lemao} for all $j\leq l(k)$
\begin{eqnarray*}
\lambda(\sigma^j(\dots,\tilde{\alpha}_{-2},\tilde{\alpha}_{-1};\tilde{\alpha}_{0},\dots ,\tilde{\alpha}_{l(k)},b\overline{ab^{s+1}}))< \max f|_{\tilde{\Lambda}}+\epsilon/2.
\end{eqnarray*}
As $(ab^s,b)$ is an alphabet, $\lambda$ has value less than $3$ in any position of the $ab^{s+1}$ in the middle of $ab^{s+1}ab^{s+1}ab^{s+1}$ and the same holds for any position of any $b$. Further, as the cut $bb^sa|bb^sa$ is good, then, for $j\geq l(k)+1$ and $j\neq l(k)+3 $ we conclude that
\begin{eqnarray*}
\lambda(\sigma^j(\dots,\tilde{\alpha}_{-2},\tilde{\alpha}_{-1};\tilde{\alpha}_{0},\dots ,\tilde{\alpha}_{l(k)},b\overline{ab^{s+1}}))&=&\\
\lambda(\sigma^j(\dots,\tilde{\alpha}_{-1};\tilde{\alpha}_{0},\dots,ab^s1,b\overline{ab^{s+1}}))&<&3< \max f|_{\tilde{\Lambda}}+\epsilon/2.
\end{eqnarray*}

For $j=l(k)+3$, we have the bad cut $221b^sb|ab^s11$ which can be compared with the bad cut $221b^{s-1}b|ab^{s-1}11$ in the $a$-continuation, then by Lemma \ref{compare} one has 
$$\lambda(\sigma^{l(k)+3}(\dots,\tilde{\alpha}_{-2},\tilde{\alpha}_{-1};\tilde{\alpha}_{0},\dots ,\tilde{\alpha}_{l(k)},b\overline{ab^{s+1}}))<\max f|_{\tilde{\Lambda}}+\epsilon/2.$$
Then taking $x=\Pi^{-1}((\dots,\tilde{\alpha}_{-2},\tilde{\alpha}_{-1};\tilde{\alpha}_{0},\dots ,\tilde{\alpha}_{l(k)},b\overline{ab^{s+1}}))$ one would have
$$x\in W^u(\tilde{\Lambda})\cap W^s(\psi_{ab^{s+1}})\subset W^u(\tilde{\Lambda})\cap W^s(\bar{\Lambda})\ \mbox{and}\ m_{\varphi,f}(x)\leq \max f|_{\tilde{\Lambda}}+\epsilon/2$$
that is a contradiction. We conclude that the $1$-continuation begins with $ba$ in these cases. 

Finally, if the $1$-continuation begins with $b1$ and the $a$-continuation begins with $ab^r1a$ where $r\geq s-1$ then we would have
$$[b1\dots]<[bab^{s+1}1ab^{s+2}\overline{ab^{s+2}}]<[ab^r1a\dots ]$$
and therefore, by Lemma \ref{lemao} for all $j\leq l(k)$
\begin{eqnarray*}
\lambda(\sigma^j(\dots,\tilde{\alpha}_{-2},\tilde{\alpha}_{-1};\tilde{\alpha}_{0},\dots ,\tilde{\alpha}_{l(k)},bab^{s+1}1ab^{s+2}\overline{ab^{s+2}}))< \max f|_{\tilde{\Lambda}}+\epsilon/2.
\end{eqnarray*}
As $(ab^{s+1},b)$ is an alphabet, $\lambda$ has value less than $3$ in any position of the $ab^{s+2}$ in the middle of $ab^{s+2}ab^{s+2}ab^{s+2}$ and the same holds for any position of any $b$. Further, as the cut 
$1b^{s+1}a|bb^{s+1}a$ is good, then, for $j\geq l(k)+1$ and $j\neq l(k)+3, l(k)+4,l(k)+2(s+3)+2 $ we conclude that
\begin{eqnarray*}
\lambda(\sigma^j(\dots,\tilde{\alpha}_{-2},\tilde{\alpha}_{-1};\tilde{\alpha}_{0},\dots ,\tilde{\alpha}_{l(k)},bab^{s+1}1ab^{s+2}\overline{ab^{s+2}}))&=&\\
\lambda(\sigma^j(\dots,\tilde{\alpha}_{-1};\tilde{\alpha}_{0},\dots,ab^s1,bab^{s+1}1ab^{s+2}\overline{ab^{s+2}}))&<&3< \max f|_{\tilde{\Lambda}}+\epsilon/2.
\end{eqnarray*}

For $j= l(k)+3, l(k)+4,l(k)+2s+8 $, we have the bad cuts $a1b^sb|ab^s11$ (two times) and $11b^sa|bb^s1a$ (the transposed of the previous cut) which can be compared with the bad cut $11b^{r-1}a|bb^{r-1}1a$ in the $a$-continuation, then as $\Sigma(t+\epsilon)$ is closed by transpositions, by Lemma \ref{compare} one has 
$$\lambda(\sigma^j(\dots,\tilde{\alpha}_{-2},\tilde{\alpha}_{-1};\tilde{\alpha}_{0},\dots ,\tilde{\alpha}_{l(k)},bab^{s+1}1ab^{s+2}\overline{ab^{s+2}}))<\max f|_{\tilde{\Lambda}}+\epsilon/2.$$
Then taking $x=\Pi^{-1}((\dots,\tilde{\alpha}_{-2},\tilde{\alpha}_{-1};\tilde{\alpha}_{0},\dots ,\tilde{\alpha}_{l(k)},bab^{s+1}1ab^{s+2}\overline{ab^{s+2}}))$ one would have
$$x\in W^u(\tilde{\Lambda})\cap W^s(\psi_{ab^{s+2}})\subset W^u(\tilde{\Lambda})\cap W^s(\bar{\Lambda})\ \mbox{and}\ m_{\varphi,f}(x)\leq \max f|_{\tilde{\Lambda}}+\epsilon/2$$
that is a contradiction. We conclude that the $1$-continuation begins with $ba$ in these cases too. This finishes the proof of the lemma. 

\end{proof}

Observe that if for some ordered alphabet, $(u,v)\in \mathscr{A}$ and $x\in \tilde{\Lambda}$, in the kneading sequence of $x$ appears the power $(uv)^i$ and $\abs{(uv)^i}\geq R$, where $R$ was given in the introduction of this section, lets say $\Pi(x)=(\dots,x_{n-1},x_n,(uv)^i, x_{n+\abs{(uv)^i}+1},\dots)$ then the point $\tilde{x}\in \Lambda$ with $ \Pi(\tilde{x})=(\dots,x_{n-1},x_n,\overline{uv})$ satisfies that   
$$\tilde{x}\in W^u(\tilde{\Lambda})\cap W^s(\psi_{uv})\subset W^u(\tilde{\Lambda})\cap W^s(\bar{\Lambda})\ \mbox{and}\ m_{\varphi,f}(\tilde{x})\leq \max f|_{\tilde{\Lambda}}+\epsilon/2$$
that is a contradiction.

\begin{lemma} [Scheme 5]\label{esquema3}
Let $(u,v)\in \mathscr{A}$ be an ordered alphabet. Suppose one has an $(u,v)$-bifurcation in $vv(uv)^iu$ where $i\geq 2$ and $\abs{u^2v}<n$. Then, the $u$-continuation begins with $uvu$ and the $v$-continuation with $vuuvu$. Additionally, if $\abs{uv}\leq \frac{3}{7}n$ and $i\geq 3$, then the $u$-continuation begins with $uvuv$ and the $v$-continuation with $vuuvuv$. 
\end{lemma}

\begin{proof}
Suppose we have an $(u,v)$-bifurcation in $vv(uv)^iu$ with $\abs{u^2v}<n$. By Lemma \ref{rere} the $u$-continuation begins with $u^{\tilde{j}}v$ for some $\tilde{j}\geq 1$. If $\tilde{j}>1$ then we have $bu^+uvuuv^-a\subseteq uvuvuuu^{\tilde{j}}v$ that contradicts Lemma \ref{bwords}. Then, the $u$-continuation begins with $uvu$. where the last $u$ is to avoid the word $u^2v^2$. If, additionally $\abs{uv}\leq \frac{3}{7}n$, $i\geq 3$ and the $u$-continuation begins with $uvu^{\tilde{j}}v$ where $\tilde{j}>1$ then we have $bu^+uvuvuuvuv^-a\subseteq uvuvuvuuvu^{\tilde{j}}v$ that contradicts, again, Lemma \ref{bwords}. Then, the $u$-continuation begins with $uvuv$ as we claimed.

For the $v$-continuation, if after $v$ there is other $v$, then we would have 
$$[vv\dots]<[v\overline{uv}]<[uvuv\dots ]$$
and therefore, as before, by Lemma \ref{lemao} for all $j\leq l(k)$
\begin{eqnarray*}
\lambda(\sigma^j(\dots,\tilde{\alpha}_{-2},\tilde{\alpha}_{-1};\tilde{\alpha}_{0},\dots ,\tilde{\alpha}_{l(k)-\abs{v^-}},v\overline{uv}))< \max f|_{\tilde{\Lambda}}+\epsilon/2.
\end{eqnarray*}
As $(u,v)$ is an alphabet, $\lambda$ has value less than $3$ in any position of the $uv$ in the middle of $uvuvuv$.
Then, for $j\geq l(k)+1$ we conclude that 
\begin{eqnarray*}
\lambda(\sigma^j(\dots,\tilde{\alpha}_{-2},\tilde{\alpha}_{-1};\tilde{\alpha}_{0},\dots ,\tilde{\alpha}_{l(k)-\abs{v^-}},v\overline{uv}))&=&
\lambda(\sigma^j(\dots,\tilde{\alpha}_{-1};\tilde{\alpha}_{0},\dots, vv(uv)^iu ,v\overline{uv}))\\&<& \max f|_{\tilde{\Lambda}}+\epsilon/2.
\end{eqnarray*}
Then taking $x=\Pi^{-1}((\dots,\tilde{\alpha}_{-2},\tilde{\alpha}_{-1};\tilde{\alpha}_{0},\dots ,\tilde{\alpha}_{l(k)-\abs{v^-}},v\overline{uv}))$ one would have
$$x\in W^u(\tilde{\Lambda})\cap W^s(\psi_{uv})\subset W^u(\tilde{\Lambda})\cap W^s(\bar{\Lambda})\ \mbox{and}\ m_{\varphi,f}(x)\leq \max f|_{\tilde{\Lambda}}+\epsilon/2$$
that is a contradiction. We conclude that the $v$-continuation begins with $vu$.

If after the $vu$ there is other $v$ then we would have 
$$[vuv\dots]<[v\overline{u(uv)^{i+1}}]<[uvuv\dots ]$$

and therefore, as before, by Lemma \ref{lemao} for all $j\leq l(k)$
\begin{eqnarray*}
\lambda(\sigma^j(\dots,\tilde{\alpha}_{-2},\tilde{\alpha}_{-1};\tilde{\alpha}_{0},\dots ,\tilde{\alpha}_{l(k)-\abs{v^-}},v\overline{u(uv)^{i+1}}))< \max f|_{\tilde{\Lambda}}+\epsilon/2.
\end{eqnarray*}

As $(u,v)$ and $(u(uv)^i,uv)$ are alphabets, $\lambda$ has value less than $3$  
in any position of the $uv$ in the middle of $uvuvuv$ and 
in any position of the $u(uv)^{i+1}$ in the middle of $u(uv)^{i+1}u(uv)^{i+1}u(uv)^{i+1}$. Also, as the cuts $Xba|baY=aba|baab$, $ (ab)^iabaa|b(ab)^i aab=ab(ab)^{i-1}abaa|b(ab)^iaab=aba(ba)^ia|b(ab)^iaab=Xba(ba)^ia|b(ab)^i\\aaY$, $ab|ab$ and $a(ab)^{i+1}|a(ab)^{i+1}=aa(ba)^ib|a(ab)^iab$ are good, corollaries \ref{3} and \ref{4} imply that we have the good cuts $b u^+ v^- a | b u^+ v^- a$, $b u^+ u(vu)^i v^- a | b u^+ (uv)^i u v^- a$, $uv|uv$ and $u(uv)^{i+1}|u(uv)^{i+1}$. 

Additionally, suppose $\abs{u}\geq \abs{v}$ and note that the word of size $3n$ ending with the beginning, $vvuv$, of $vv(uv)^iu$ is $(u,v)$-semi renormalizable with renormalization kernel $\tilde{\gamma}=\tilde{X}vvuv$ and 
$$\abs{\tilde{X}}=\abs{\tilde{\gamma}}-\abs{vvuv}=3n-\abs{\tilde{\omega_1}} 
-\abs{vvuv}>3n-\abs{uv} 
-\abs{vvuv}\geq 3n-2\abs{uuv}>n>\abs{u}.$$
Then as $Xbb(ab)^iab|a(ab)^{i+1}aY=Xb(ba)^{i+1}b|a(ab)^{i+1}aab$ is bad (where $X$ is like $\tilde{X}$ but replacing the $u$'s by $a$ and the $v$'s by $b$), Corollary \ref{1} implies that we have the bad cut $b u^+ (vu)^{i+1} v^- b | a u^+ (uv)^{i+1} v^- a$ in this case. If $\abs{u}<\abs{v}$ then $v=u\tilde{v}=\tilde{v}^-abu^+$ where $(u,\tilde{v})$ is an ordered alphabet, and then in this case we have again the bad cut $b u^+ (vu)^{i+1} v^- b | a u^+ (uv)^{i+1} v^- a$. Then, for $j\geq l(k)+1$ and $j\neq l(k)+3 $ we conclude that 
\begin{eqnarray*}
\lambda(\sigma^j(\dots,\tilde{\alpha}_{-2},\tilde{\alpha}_{-1};\tilde{\alpha}_{0},\dots ,\tilde{\alpha}_{l(k)-\abs{v^-}},v\overline{u(uv)^{i+1}}))&=&\\
\lambda(\sigma^j(\dots,\tilde{\alpha}_{-1};\tilde{\alpha}_{0},\dots,vv(uv)^iu,v\overline{u(uv)^{i+1}}))&<&3< \max f|_{\tilde{\Lambda}}+\epsilon/2.
\end{eqnarray*}

For $j=l(k)+3$, let us consider $(uv)^r$ the biggest power of $uv$ that appears in the $u$-continuation. We need to consider some cases: 

Suppose first that $r\geq i$ and after $(uv)^r$ there is the word $uu$. If $\abs{u}\geq \abs{v}$, note that the word of size $3n$ beginning with the last $uvuu$ of $(uv)^ruu$ is $(u,v)$-semi renormalizable with renormalization kernel $\tilde{\gamma}\neq uvu^s$ (Lemma \ref{rere}), then $\tilde{\gamma}=uvuu\tilde{Y}$ where $\tilde{Y}$ is a word in the alphabet $\{u,v\}$ with at least one $v$. As the cut $Xb(ba)^ib|a(ab)^raaY=Xb(ba)^ib|a(ab)^i(ab)^{r-i}aaY$ is bad, where $X$ is defined as in the previous case and $Y$ is like $\tilde{Y}$ but replacing the $u$'s by $a$ and the $v$'s by $b$, then by Corollary \ref{1} we have in the $u$-continuation the bad cut $b u^+ (vu)^i v^- b | a u^+ (uv)^i v^- a$. If $\abs{u}<\abs{v}$ then $v=u\tilde{v}=\tilde{v}^-abu^+$ where $(u,\tilde{v})$ is an ordered alphabet, and then in this case we have again the bad cut $b u^+ (vu)^i v^- b | a u^+ (uv)^i v^- a$. So, by Lemma \ref{compare} one has 
$$\lambda(\sigma^{l(k)+3}(\dots,\tilde{\alpha}_{-2},\tilde{\alpha}_{-1};\tilde{\alpha}_{0},\dots ,\tilde{\alpha}_{l(k)-\abs{v^-}},v\overline{u(uv)^{i+1}}))<\max f|_{\tilde{\Lambda}}+\epsilon/2.$$

If $r< i$ and after $(uv)^r$ there is the word $uu$ then, as before, we can find the word $Y$ in any case, if $\abs{u}\geq \abs{v}$, the word $X$ and if $\abs{u}<\abs{v}$, there is nothing to do. As the cut $b(ba)^{i-r-1}ba(ba)^rb|a(ab)^raa$ is bad, in any case, we have in the $u$-continuation the bad cut $b u^+ u(vu)^r v^- b | a u^+ (uv)^ru v^- a$ and then, by Lemma \ref{compare} one has 
$$\lambda(\sigma^{l(k)+3}(\dots,\tilde{\alpha}_{-2},\tilde{\alpha}_{-1};\tilde{\alpha}_{0},\dots ,\tilde{\alpha}_{l(k)-\abs{v^-}},v\overline{u(uv)^{i+1}}))<\max f|_{\tilde{\Lambda}}+\epsilon/2.$$

If $i+1\leq r$ and after $(uv)^r$ there is one $v$ then, as the cut 
$b(ba)^ib|a(ab)^iaY=b(ba)^ib|a(ab)^i(ab)(ab)^{r-i-1}b=b(ba)^ib|a(ab)^rb$
is bad, we can consider the same cases as before: if $\abs{u}\geq\abs{v}$ we can find the word $X$ and use Lemma \ref{1} or if $\abs{u}<\abs{v}$ there is nothing to do. Then, we have in the $u$-continuation the bad cut $b u^+ (vu)^i v^- b | a u^+ (uv)^i v^- a$ and then, by Lemma \ref{compare} one has 
$$\lambda(\sigma^{l(k)+3}(\dots,\tilde{\alpha}_{-2},\tilde{\alpha}_{-1};\tilde{\alpha}_{0},\dots ,\tilde{\alpha}_{l(k)-\abs{v^-}},v\overline{u(uv)^{i+1}}))<\max f|_{\tilde{\Lambda}}+\epsilon/2.$$

Finally, if $i\geq r$ and after $(uv)^r$ there is one $v$ then, as the cut $ a(ba)^{r-1}a|b(ab)^{r-1}b \subseteq b(ba)^{i+1}a|b(ab)^{r-1}b=bb(ab)^iaa|b(ab)^{r-1}b$ is bad where by Corollary \ref{2} implies that we have in the $u$-continuation the bad cut $au^+(vu)^{r-1}v^-a|bu^+(uv)^{r-1}v^-b$ and then, by Lemma \ref{transpuesto} one has the cut 
\begin{eqnarray*}
    (au^+(vu)^{r-1}v^-a|bu^+(uv)^{r-1}v^-b)^*&=&b (u^+ (uv)^{r-1} v^-)^*b|a (u^+ (vu)^{r-1} v^-)^* a\\ &=& b (u^+ (vu)^{r-1} v^-)b|a (u^+ (uv)^{r-1} v^-) a \\&=& b u^+ (vu)^{r-1} v^-b|a u^+ (uv)^{r-1} v^- a
    \end{eqnarray*}
in the transposed of $\alpha=(\alpha_n)_{n\in \mathbb{Z}}\in \Pi(\tilde{\Lambda})\subset\Sigma(\max f|_{\tilde{\Lambda}})$. As $\Sigma(\max f|_{\tilde{\Lambda}})$ is closed by transpositions, by Lemma \ref{compare} one has
$$\lambda(\sigma^{l(k)+3}(\dots,\tilde{\alpha}_{-2},\tilde{\alpha}_{-1};\tilde{\alpha}_{0},\dots ,\tilde{\alpha}_{l(k)-\abs{v^-}},v\overline{u(uv)^{i+1}}))<\max f|_{\tilde{\Lambda}}+\epsilon/2,$$
again. Then taking $x=\Pi^{-1}((\dots,\tilde{\alpha}_{-2},\tilde{\alpha}_{-1};\tilde{\alpha}_{0},\dots ,\tilde{\alpha}_{l(k)-\abs{v^-}},v\overline{u(uv)^{i+1}}))$ one would have
$$x\in W^u(\tilde{\Lambda})\cap W^s(\psi_{u(uv)^{i+1}})\subset W^u(\tilde{\Lambda})\cap W^s(\bar{\Lambda})\ \mbox{and}\ m_{\varphi,f}(x)\leq \max f|_{\tilde{\Lambda}}+\epsilon/2$$
that is a contradiction. We conclude that after the $vu$ follows other $u$ and then the $v$-continuation begins with $vuu$. Repeating the same argument of Scheme $1$, one has the same conclusions for the $v$-continuation here, as we claimed. 

\end{proof}

\begin{lemma}[Scheme 6]\label{esquema4}

Let $(u,v)\in \mathscr{A}$ be an ordered alphabet. Suppose one has an $(u,v)$-bifurcation in $u(uv)^i$ where $i\geq 2$ and $\abs{u^2v}\leq n$. Then, the $u$-continuation begins with $uvv$.
\end{lemma}
\begin{proof}
First observe that by Lemma \ref{rere} the $u$-continuation begins with $u^{\tilde{j}}v$ for some $\tilde{j}\geq 1$. If $\tilde{j}>1$ then, as $u^{\tilde{j}} v$ and $uv$ begin with $v^-a$ we would have 
$$[v\dots]<[\overline{uv}]<[uu^{\tilde{j}-1}v\dots ]$$
and therefore, as before, by Lemma \ref{lemao} for all $j\leq l(k)$
\begin{eqnarray*}
\lambda(\sigma^j(\dots,\tilde{\alpha}_{-2},\tilde{\alpha}_{-1};\tilde{\alpha}_{0},\dots ,\tilde{\alpha}_{l(k)-\abs{v^-}},\overline{uv}))< \max f|_{\tilde{\Lambda}}+\epsilon/2.
\end{eqnarray*}
As $(u,v)$ is an alphabet, as before, $\lambda$ has value less than $3$ in any position of the $uv$ in the middle of $uvuvuv$.
Then, for $j\geq l(k)+1$ we conclude that 
\begin{eqnarray*}
\lambda(\sigma^j(\dots,\tilde{\alpha}_{-1};\tilde{\alpha}_{0},\dots ,\tilde{\alpha}_{l(k)-\abs{v^-}},\overline{uv}))=
\lambda(\sigma^j(\dots,\tilde{\alpha}_{-1};\tilde{\alpha}_{0},\dots,u(uv)^i,\overline{uv}))< \max f|_{\tilde{\Lambda}}+\frac{\epsilon}{2}.
\end{eqnarray*}
Then taking $x=\Pi^{-1}((\dots,\tilde{\alpha}_{-2},\tilde{\alpha}_{-1};\tilde{\alpha}_{0},\dots ,\tilde{\alpha}_{l(k)-\abs{v^-}},\overline{uv}))$ one would have
$$x\in W^u(\tilde{\Lambda})\cap W^s(\psi_{uv})\subset W^u(\tilde{\Lambda})\cap W^s(\bar{\Lambda})\ \mbox{and}\ m_{\varphi,f}(x)\leq \max f|_{\tilde{\Lambda}}+\epsilon/2$$
that is a contradiction. We conclude that the $u$-continuation begins with $uv$.

Again, by Lemma \ref{rere} we conclude that the $u$-continuation begins with $uvu^{\tilde{j}}v$ for some $\tilde{j}\geq 0$. If $\tilde{j}>0$ then we would have 
$$[v\dots]<[uvv\overline{(uv)^{i+1}v}]<[uvu^{\tilde{j}}v\dots ]$$
and therefore, as before, by Lemma \ref{lemao} for all $j\leq l(k)$
\begin{eqnarray*}
\lambda(\sigma^j(\dots,\tilde{\alpha}_{-2},\tilde{\alpha}_{-1};\tilde{\alpha}_{0},\dots ,\tilde{\alpha}_{l(k)-\abs{v^-}},uvv\overline{(uv)^{i+1}v}))< \max f|_{\tilde{\Lambda}}+\epsilon/2.
\end{eqnarray*}

As $(uv,(uv)^iv)$ is an alphabet, $\lambda$ has value less than $3$ in any position of the $uv(uv)^iv$ in the middle of $uv(uv)^ivuv(uv)^ivuv(uv)^iv$. As the cut $(ab)^iabb|ab(ab)^ib=a(ba)^ibb|ab(ab)^ib$ is good, corollary \ref{4} imply that we have the good cut $u(vu)^ivv|uv(uv)^iv$. 

Additionally, as $a(ab)^ia|bb(ab)^i=a(ab)^ia|b(ba)^ib$ is bad, Corollary \ref{2} implies that we have the bad cut $au^+(uv)^iv^-a|bu^+(vu)^iv^-b$. Then, for $j\geq l(k)+1$ and $j\neq l(k)+2 $ we conclude that 
\begin{eqnarray*}
\lambda(\sigma^j(\dots,\tilde{\alpha}_{-2},\tilde{\alpha}_{-1};\tilde{\alpha}_{0},\dots ,\tilde{\alpha}_{l(k)-\abs{v^-}},uvv\overline{(uv)^{i+1}v}))&=&\\
\lambda(\sigma^j(\dots,\tilde{\alpha}_{-1};\tilde{\alpha}_{0},\dots,u(uv)^i,uvv\overline{(uv)^{i+1}v}))&<&3< \max f|_{\tilde{\Lambda}}+\epsilon/2.
\end{eqnarray*}

For $j=l(k)+2$, let us consider $(uv)^r$ the biggest power of $uv$ that appears after $v$ in the $v$-continuation. We need to consider some cases: 

If $r\geq i-1$, as the cut $a(ab)^{i-1}a|bb(ab)^{i-1}=a(ab)^{i-1}a|b(ba)^{i-1}b$ is bad, Corollary \ref{2} implies that we have in the $v$-continuation the bad cut $au^+(uv)^{i-1}v^-a|bu^+(vu)^{i-1}v^-b$ and then, by Lemma \ref{compare} one has 

$$\lambda(\sigma^{l(k)+2}(\dots,\tilde{\alpha}_{-2},\tilde{\alpha}_{-1};\tilde{\alpha}_{0},\dots ,\tilde{\alpha}_{l(k)-\abs{v^-}},uvv\overline{(uv)^{i+1}v}))<\max f|_{\tilde{\Lambda}}+\epsilon/2.$$

If $r+2\leq i$ and after $(uv)^r$ there is one $v$ then, as the cut $ab(ab)^ra|bb(ab)^rb=ab(ab)^ra|b(ba)^rbb$ is bad, Corollary \ref{2} implies that we have in the $v$-continuation the bad cut $au^+v(uv)^rv^-a|bu^+(vu)^rvv^-b$ and then, by Lemma \ref{compare} one has 
$$\lambda(\sigma^{l(k)+2}(\dots,\tilde{\alpha}_{-2},\tilde{\alpha}_{-1};\tilde{\alpha}_{0},\dots ,\tilde{\alpha}_{l(k)-\abs{v^-}},uvv\overline{(uv)^{i+1}v}))<\max f|_{\tilde{\Lambda}}+\epsilon/2.$$

Finally, if $r+2\leq i$ and after $(uv)^r$ there is the word $uu$, note that the word of size $3n$ beginning with the last $uvuu$ of $v(uv)^ruu$ is $(u,v)$-semi renormalizable with renormalization kernel $\tilde{\gamma}=uvuuY$ and 
$$\abs{Y}=\abs{\tilde{\gamma}}-\abs{uvuu}=3n-\abs{\tilde{\omega_2}} 
-\abs{uvuu}>3n-\abs{uv} 
-\abs{uvuu}=3n-2\abs{uuv}>n>\abs{v}.$$
As the cut $Xb(ab)^rb|ab(ab)^{r-1}aaY=ab(ab)^rb|a(ba)^raY$ is bad, Corollary \ref{1} implies that we have in the $v$-continuation the bad cut $b u^+ (uv)^r v^-b |a u^+ (vu)^r v^- a$ and then, by Lemma \ref{transpuesto} one has the cut 
\begin{eqnarray*}
    (b u^+ (uv)^r v^-b |a u^+ (vu)^r v^- a)^*&=&a (u^+ (vu)^r v^-)^*a|b (u^+ (uv)^r v^-)^* b\\ &=& a (u^+ (uv)^r v^-)a|b (u^+ (vu)^r v^-) b \\&=& au^+(uv)^r v^-a|b u^+ (vu)^r v^-b
    \end{eqnarray*}
in the transposed of $\tilde{\alpha}=(\tilde{\alpha}_n)_{n\in \mathbb{Z}}\in \Pi(\tilde{\Lambda})\subset\Sigma(\max f|_{\tilde{\Lambda}})$. As $\Sigma(\max f|_{\tilde{\Lambda}})$ is closed by transpositions, by Lemma \ref{compare} one has
$$\lambda(\sigma^{l(k)+2}(\dots,\tilde{\alpha}_{-2},\tilde{\alpha}_{-1};\tilde{\alpha}_{0},\dots ,\tilde{\alpha}_{l(k)-\abs{v^-}},uvv\overline{(uv)^{i+1}v}))<\max f|_{\tilde{\Lambda}}+\epsilon/2,$$
again. Then taking $x=\Pi^{-1}((\dots,\tilde{\alpha}_{-2},\tilde{\alpha}_{-1};\tilde{\alpha}_{0},\dots ,\tilde{\alpha}_{l(k)-\abs{v^-}},uvv\overline{(uv)^{i+1}v}))$ one would have
$$x\in W^u(\tilde{\Lambda})\cap W^s(\psi_{(uv)^{i+1}v})\subset W^u(\tilde{\Lambda})\cap W^s(\bar{\Lambda})\ \mbox{and}\ m_{\varphi,f}(x)\leq \max f|_{\tilde{\Lambda}}+\epsilon/2$$
that is a contradiction. We conclude that $\tilde{j}=0$ and then that the $u$-continuation begins with $uvv$ as we claimed.

\end{proof}
 
Now we are ready to consider the cases presented in the introduction of this section. Remember that we are considering the word $\tilde{\kappa}$ of size $3n$ just before the bifurcation and the ordered alphabet $(\alpha,\beta)\in \mathscr{A}$ is such that $\abs{\alpha},\abs{\beta}<n$, $\abs{\alpha \beta}\geq n$ and $\tilde{\kappa}$ is $(\alpha,\beta)$-semi renormalizable. Here $\kappa$ is the word of biggest length in $\{\alpha,\beta\}$ and $\gamma$ the renormalization kernel, we have the following cases:

\subsection{$\kappa=\alpha$ and $\kappa$ appears in $\gamma$}
In this case we have $\alpha=\tilde{\alpha}\beta^s$ where $\abs{\beta}\geq \abs{\tilde{\alpha}}$ and $s\geq 1$. Suppose first that $\abs{\beta}> \abs{\tilde{\alpha}}$, then for some $k\geq 1$, $\beta=\tilde{\alpha}^k\tilde{\beta}$ with $(\tilde{\alpha},\tilde{\beta})$ an ordered alphabet. Therefore $\alpha=u(uv)^s$ where $u=\tilde{\alpha}$, $v=\tilde{\alpha}^{k-1}\tilde{\beta}$ and then $\beta=uv$. Consider the last appearance of $\alpha$ in $\gamma$ and let $\omega_2$ as in the definition of $(\alpha,\beta)$-weakly renormalizable, then the end of $\tilde{\kappa}$ (see the comments at the beginning of this subsection) is $\alpha\beta^{i_0}\omega_2=u(uv)^{s+i_0}\omega_2$ where $\abs{\omega_2}<\abs{\alpha}=\abs{u(uv)^s}$, $i_0\geq 0$ and $\omega_2$ is a prefix of $\alpha$. Observe that by Lemma \ref{14}, the word of size $3n$ beginning in the last $uv$ of $(uv)^{s+i_0}\omega_2$ is $(u,v)$-semi renormalizable. Consider the largest subword $\eta$ of $u(uv)^s$ in the alphabet $\{ u,v\}$ contained in $\omega_2$. We have the following subcases: 

\textbf{a)} $\eta=\emptyset :$ This case is not possible because we would have an $(u,v)$-bifurcation after $u(uv)^s$ and as the word of size $3n$ beginning in the last $uv$ of $(uv)^{s+i_0}$ is $(u,v)$-semi renormalizable, by Lemma \ref{comienzo} and Corollary \ref{9} the $u$-continuation begins with $v^-a$ and as the $v$-continuation begins with $v^-b$ we conclude that $\omega_2=v^-$ (as we claimed before). But this is a contradiction by Lemma \ref{alfacuadrado} because 
$\abs{u(uv)^{s+i_0}v}=\abs{u(uv)^{s+i_0}\omega_2}+2\leq 3n+2$. In general, $\omega_2=\eta v^-$.

\textbf{b)} $\eta=u :$ Let us suppose first that $\abs{uv}\leq n/5$. As
$$\frac{n}{5}\cdot (s+2)\geq (s+2)\abs{uv}=\abs{(uv)^{s+2}}>\abs{u(uv)^{s+1}}=\abs{\alpha\beta}\geq n,$$
we conclude that $s\geq 4$.

We have an $(u,v)$-bifurcation in $u(uv)^{s+i_0}u$. By Scheme $1$ (Lemma \ref{esquema1}) the $u$-continuation begins with $(uv)^{s+i_0}$ and the $v$-continuation with $vuuvuv$. Then we determine in the continuation starting with $2$, at least  
$$\abs{(uv)^{s+i_0}}-\abs{v}\geq \abs{(uv)^s}-\abs{v}=\abs{u(uv)^{s+1}}-\abs{(uv)^2}=\abs{\alpha\beta}-2\abs{uv}\geq n-2n/5=3n/5$$
letters. 

Now, as the $v$-continuation begins with $vuuvuv$ one can consider the word os size $3n$ starting in the last $uuvuv$, this word is $(uuv,uv)$-semi renormalizable and then after $uuvuv$ it follows or $uv$ or $uuv$. Here we have some options, either the renormalization kernel $\tilde{\gamma}$ is of the form $u(uv)^{\tilde{s}}$ and then 
$$\abs{u(uv)^{\tilde{s}}}=3n-\abs{\tilde{\omega}_2}>3n-\abs{uuv}>2n$$
where $\tilde{\gamma}$ and $\tilde{\omega}_2$ are as in Definition \ref{def:weaklyrenormalizable}. Or if $(uv)^{\tilde{s}}$ is the biggest power of $uv$ that follows the first $u$ in the $v$-continuation and $\tilde{s}\geq s+i_0$ then
$$\abs{u(uv)^{\tilde{s}}}\geq\abs{u(uv)^{s+i_0}}\geq\abs{u(uv)^s}=\abs{u(uv)^{s+1}}-\abs{uv}=\abs{\alpha\beta}-\abs{uv}>4n/5,$$
or if $\tilde{s}\leq s+i_0-1= s+i_0+1-2$ then one has the subword $(uv)^{s+i_0+1}u(uv)^{\tilde{s}}uuv$ and then by Lemma \ref{s}, $\abs{u(uv)^{\tilde{s}+1}}\geq 3(n-1)$ which let us conclude that 
$$\abs{u(uv)^{\tilde{s}}}=\abs{u(uv)^{\tilde{s}+1}}-\abs{uv}\geq 3(n-1)-\abs{uv}>2n$$
Then we determine in the continuation starting with $1$, at least  
$$\abs{vu(uv)^{\tilde{s}}}-\abs{v}=\abs{u(uv)^{\tilde{s}}}>4n/5$$
letters in any case.

Now, if $n/5\leq \abs{uv}<\abs{uuv}< n$ and $s+i_0\geq 2$, Lemma \ref{rere} let us conclude as in Scheme $1$ (Lemma \ref{esquema1}) that the $u$-continuation begins with $uv$. If $s+i_0=1$ then, as $\abs{uuv}<n$, $\gamma$ contains other letter $\alpha$ or $\beta$ before $uuv$ and then we have an $(u,v)$-bifurcation in $uvuuvu$ and the $u$-continuation begins with $u^{\tilde{j}}v$ by lemma \ref{rere2}. If $\tilde{j}>1$, using that 
$$[v\dots]<[uv\overline{uuv}]<[u^{\tilde{j}}v\dots ]$$
and that $(u,uv)$ is alphabet, we conclude again that 
$$\lambda(\sigma^j(\dots,\tilde{\alpha}_{-2},\tilde{\alpha}_{-1};\tilde{\alpha}_{0},\dots ,\tilde{\alpha}_{l(k)-\abs{v^-}},uv\overline{uuv}))< \max f|_{\tilde{\Lambda}}+\epsilon/2$$
and then $\tilde{j}=1$. Then, the $u$-continuation begins with $uv$ in any case. Now, the letter that follows $uv$ in the $u$-continuation cannot be $v$ because in other case we would have $u^2v^2$ that contradicts Lemma \ref{alfacuadrado}. Using again the alphabet $(u,uv)$ we see that after $uv$ it follows either $u^{\tilde{j}}$ or $u^{\tilde{j}}v$ for some $\tilde{j}$ and in any case $v^-a$ by Lemma \ref{comienzo}. Then we determine in the continuation starting with $2$, at least  
$$\abs{uvv^-a}-\abs{v}=\abs{uv}\geq n/5$$
letters. 
On the other hand, arguing as before, it follows that the $v$-continuation begins with $vuuvu$ and then we determine in the continuation starting with $1$, at least  
$$\abs{vuuvu}-\abs{v}>\abs{uv}\geq n/5.$$

\textbf{c)} $\eta=u(uv)^{i_1}$ where $i_1>0:$ This case is not possible because on would have in the $v$-continuation $u(uv)^{i_1}v=u^2(vu)^{i_1-1}v^2$ and $\abs{u(uv)^{i_1}v}<2n$.

\textbf{d)} $\eta=u^2 :$ Let us suppose first that $s+i_0\geq 2$, then we have an $(u,v)$-bifurcation in $u(uv)^{s+i_0}uu$. For the $u$-continuation one has the sequence $u(uv)^{s+i_0}uuu$ and then the word of size $3n$ starting in the last $uv$ of $u(uv)^{s+i_0}uuu$ is $(u,v)$-semi renormalizable. Let $\tilde{\gamma}=uvuuuX$ and $\tilde{\omega_2}$ as in the definition or renormalizartion. As $\abs{\tilde{\omega_2}}<\abs{uv}$, we conclude that $\abs{u\tilde{\omega}_2}<\abs{uuv}\leq \abs{\alpha}< n$ and then
$$3n=\abs{\tilde{\gamma}\tilde{\omega}_2}=\abs{\tilde{\gamma}}+\abs{\tilde{\omega}_2}<\abs{uvu}+\abs{u\tilde{\omega}_2}+\abs{uX}<\abs{uvu}+n+\abs{uX}<2n+\abs{uX}.$$
But this implies that $\abs{uv}<n<\abs{uX}$ and by Lemmas \ref{comienzo} and \ref{9}, $uX$ begins with $v^-a$, so $bu^+uvuuv^-a\subseteq uvuvuuuX \subseteq 
(uv)^{s+i_0}uuuX$ which contradicts Lemma \ref{bwords}. We conclude that  $s+i_0=s=1$. 



Let us consider first the $u$-continuation. Note that as $\abs{uuv}<n$ in $\gamma$ there is other letter before $uuv$ and then there is other $uv$, observe that the word of size $3n$ beginning in the last $uv$ of $uvuuvuuu$ is $(u,v)$-semi renormalizable. Let $\tilde{\gamma}=uvuuuX$ and $\tilde{\omega_2}$ as in the definition or renormalization. As before, $\abs{u\tilde{\omega}_2}<\abs{uuv}< n$ and then
$$3n=\abs{\tilde{\gamma}\tilde{\omega}_2}=\abs{\tilde{\gamma}}+\abs{\tilde{\omega}_2}<\abs{uvu}+\abs{u}+\abs{u\tilde{\omega}_2}+\abs{X}<2\abs{uvu}+\abs{u}+\abs{X}<2n+\abs{u}+\abs{X},$$
this implies that $\abs{uuv}<n<\abs{u}+\abs{X}$ and then $\abs{uv}<\abs{X}$. By Lemmas \ref{comienzo} and \ref{9}, if $X$ begins with $u$ then actually begins with $v^-a$, but then $bu^+u^2vu^3v^-a\subseteq uvuuvu^3X$ contradicts Lemma \ref{bwords}. Then, in that case the $u$-continuation begins with $uv$ and the word of size $3n$ starting in this $uv$ is $(u,v)$-semi renormalizable. Let $\tilde{\gamma}=uvX$ and $\tilde{\omega_2}$ as in the definition or renormalizartion. As $\abs{\tilde{\omega_2}}<\abs{uv}<n$, we conclude that 
$$3n=\abs{\tilde{\gamma}\tilde{\omega}_2}=\abs{\tilde{\gamma}}+\abs{\tilde{\omega}_2}<\abs{uv}+n+\abs{X}<2n+\abs{X},$$
this implies that $\abs{uv}<n<\abs{X}$ and by Lemmas \ref{comienzo} and \ref{9}, $X$ begins with $v^-$, so the $u$-continuation begins with $uvv^-$ in this case.

Observe that $n\leq \abs{\alpha\beta}=\abs{u(uv)^2}<3\abs{uv}$ and then $\abs{uv}>n/3$. So, we determine in the continuation beginning with $2$, at least $\abs{uvv^-}-\abs{v^-}=\abs{uv}>n/3$ letters.

For the $v$-continuation, one has that the word beginning in the last $uv$ of $uuvuuv$ is $(u,v)$-semi renormalizable, if after $v$ we have other $v$, then we would have 
$$[vv\dots]<[v\overline{uuv}]<[(uv)^{s+i_0}\dots ]$$
and again, by Lemma \ref{lemao} for all $j\leq l(k)$
\begin{eqnarray*}
\lambda(\sigma^j(\dots,\tilde{\alpha}_{-2},\tilde{\alpha}_{-1};\tilde{\alpha}_{0},\dots ,\tilde{\alpha}_{l(k)-\abs{v^-}},v\overline{uuv}))< \max f|_{\tilde{\Lambda}}+\epsilon/2.
\end{eqnarray*}
As $(u,uv)$ is an alphabet, $\lambda$ has value less than $3$ in any position of the $uuv$ in the middle of $uuvuuvuuv$.
Then, for $j\geq l(k)+1$ we conclude that 
\begin{eqnarray*}
\lambda(\sigma^j(\dots,\tilde{\alpha}_{-2},\tilde{\alpha}_{-1};\tilde{\alpha}_{0},\dots ,\tilde{\alpha}_{l(k)-\abs{v^-}},v\overline{uuv}))&=&
\lambda(\sigma^j(\dots,\tilde{\alpha}_{-1};\tilde{\alpha}_{0},\dots,uuvuu,v\overline{uuv}
 ))\\ &<& \max f|_{\tilde{\Lambda}}+\epsilon/2.
\end{eqnarray*}
Then taking $x=\Pi^{-1}((\dots,\tilde{\alpha}_{-2},\tilde{\alpha}_{-1};\tilde{\alpha}_{0},\dots ,\tilde{\alpha}_{l(k)-\abs{v^-}},v\overline{uuv}))$ one would have
$$x\in W^u(\tilde{\Lambda})\cap W^s(\psi_{uv})\subset W^u(\tilde{\Lambda})\cap W^s(\bar{\Lambda})\ \mbox{and}\ m_{\varphi,f}(x)\leq \max f|_{\tilde{\Lambda}}+\epsilon/2$$
that is a contradiction. We conclude that the $v$-continuation begins with $vu$.

Let $\tilde{\gamma}=uvuX$ and $\tilde{\omega_2}$ as in the definition or renormalizartion. Again, as $\abs{\tilde{\omega_2}}<n$, we conclude that 
$$3n=\abs{\tilde{\gamma}\tilde{\omega}_2}=\abs{\tilde{\gamma}}+\abs{\tilde{\omega}_2}<\abs{uvu}+n+\abs{X}<2n+\abs{X},$$
which implies that $\abs{uv}<n<\abs{X}$ and by Lemmas \ref{comienzo} and \ref{9}, $X$ begins with $v^-$, so the $v$-continuation begins with $vuv^-$ and as before, we determine in the continuation beginning with $1$, at least $\abs{vuv^-}-\abs{v^-}=\abs{uv}>n/3$ letters.

\textbf{e)} $\eta=u(uv)^{i_1}u :$ This case is similar with case \textbf{b)}. Let us suppose first that $i_1\neq 0$, then we have an $(u,v)$-bifurcation in $u(uv)^{s+i_0}u(uv)^{i_1}u$ and using the same arguments of Scheme $1$ (Lemma \ref{esquema1}) it follows that the $u$-continuation begins with $uv$. Now, using the second item of Lemma \ref{force} we get that $i_1\geq s+i_0-1\geq s-1$.

As before, if $\abs{uv}\leq n/5$, we conclude that $s\geq 4$ and then $s-1\geq 3$. Then, by Scheme $1$ (Lemma \ref{esquema1}), the $u$-continuation begins with $(uv)^{i_1}$ and the $v$-continuation with $vuuvuv$.  Then, we determine in the continuation beginning with $2$, at least  
$$\abs{(uv)^{i_1}}-\abs{v}\geq \abs{(uv)^{s-1}}-\abs{v}=\abs{u(uv)^{s+1}}-\abs{(uv)^3}=\abs{\alpha\beta}-3\abs{uv}\geq n-3n/5=2n/5$$
letters. 

Again, as the $v$-continuation begins with $vuuvuv$, one can consider the word of size $3n$ starting in the last $uuvuv$, this word is $(uuv,uv)$-semi renormalizable and we have some options: either the renormalization kernel $\tilde{\gamma}$ is of the form $u(uv)^{\tilde{s}}$ and then 
$$\abs{u(uv)^{\tilde{s}}}=3n-\abs{\tilde{\omega}_2}>3n-\abs{uuv}>2n$$
or if $(uv)^{\tilde{s}}$ is the biggest power of $uv$ that follows the first $u$ in the $v$-continuation and $\tilde{s}\geq i_1$ then
$$\abs{u(uv)^{\tilde{s}}}\geq\abs{u(uv)^{i_1}}\geq\abs{u(uv)^{s-1}}=\abs{u(uv)^{s+1}}-\abs{(uv)^2}=\abs{\alpha\beta}-2\abs{uv}>3n/5.$$
The case $\tilde{s}\leq i_1-1= i_1+1-2$ is not possible by Lemma \ref{s} because one would have the subword $(uv)^{i_1+1}u(uv)^{\tilde{s}}uuv$ and $\abs{u(uv)^{i_1+1}}<\abs{\omega_2}+\abs{v}<2n$.
Then we determine in the continuation starting with $1$, at least  
$$\abs{vu(uv)^{\tilde{s}}}-\abs{v}=\abs{u(uv)^{\tilde{s}}}>3n/5$$
letters in any case.

Finally, if $\abs{uv}>n/5$, the same arguments of case $\textbf{b)}$ let us determine at least $n/5$ letters in both continuations.

The case where $i_1=0$ was already considered in the case $\textbf{d)}$.

Remember that we assumed that $\abs{\beta}> \abs{\tilde{\alpha}}$, where $\alpha=\tilde{\alpha}\beta^s$. Now, if $\abs{\beta}= \abs{\tilde{\alpha}}$, then $(\tilde{\alpha},\beta)=(a,b)$ and therefore $\alpha=ab^s$. As before, consider the last appearance of $\alpha$ in $\gamma$ and let $\omega_2$ as in the definition of $(\alpha,\beta)$ weakly renormalizable, then the end of $\tilde{\kappa}$ is, unless one letter at the end of $\omega_2$, $ab^{s+i_0}\omega_2$ where $\abs{\omega_2}<\abs{\alpha}=\abs{ab^s}$ and $\omega_2$ is a prefix of $\alpha$. Consider the largest subword $\eta$ of $ab^s$ in the alphabet $\{ a,b\}$ contained in $\omega_2$. We have the following subcases: 

\textbf{f)} $\eta=\emptyset :$ In this case, we have an $(a,1)$-bifurcation in $a(1)^i$ where 
$$\abs{(1)^i}\geq \abs{b^{s+i_0}}-1= \abs{ab^{s+i_0+1}}-5\geq \abs{\alpha\beta}-5\geq n-5.$$
By Scheme $4$ (Lemma \ref{muchos1}) the
$1$-continuation begins with $ba$. Observe that in any sequence of $\Sigma(3)$, after the sequence $b^ra$ always follows $b^{r-1}$ because in other case we would have the bad cut $bb^{j}b|ab^{j}a$ where $j<r-1$. Additionally, as by Theorem \ref{palabras}, $\Sigma(3+6^{-3n},n)=\Sigma(3,n)$, if we consider the word of size $n$ starting in the first $b$ of the smallest power $b^r$ of $b$ at the end of $(1)^i$ such that $\abs{b^r}\geq n/3$, then we determine in both continuations at least $\abs{b^{r-1}}+2=\abs{b^r}\geq n/3$ letters.

\textbf{g)} $\eta=a :$ This case is not possible because one would have the bad cut $bb|aa$ in the $a$-continuation.

\textbf{h)} $\eta=ab^{i_1}$ where $i_1>0 :$ In this case, we have an $(a,1)$-bifurcation in $a(1)^i$, where as before, one has $\abs{(1)^i}>n/3$ and we can consider $r$ as in the case \textbf{f)} to conclude that we determine in both continuations at least $\abs{b^r}\geq n/3$ letters.

\subsection{$\kappa=\beta$ and $\kappa$ appears in $\gamma$}

In this case we have $\beta=\alpha^s\tilde{\beta}$ where $\abs{\alpha}\geq \abs{\tilde{\beta}}$ and $s\geq 1$. Suppose first that $\abs{\alpha}> \abs{\tilde{\beta}}$, then for some $k\geq 1$, $\alpha=\tilde{\alpha}\tilde{\beta}^k$ with $(\tilde{\alpha},\tilde{\beta})$ alphabet. Then $\beta=(uv)^sv$ where $u=\tilde{\alpha}\tilde{\beta}^{k-1}$, $v=\tilde{\beta}$ and then $\alpha=uv$. Consider the last appearance of $\beta$ in $\gamma$ and let $\omega_2$ as in the definition of $(\alpha,\beta)$ weakly renormalizable, then the end of $\tilde{\kappa}$ is $\beta\alpha^{i_0}\omega_2=(uv)^sv(uv)^{i_0}\omega_2$ where $\abs{\omega_2}<\abs{\beta}=\abs{(uv)^sv}$ and $\omega_2$ is a prefix of $(uv)^{s+1}$. Consider the biggest subword $\eta$ of $(uv)^{s+1}$ in the alphabet $\{ u,v\}$ contained in $\omega_2$. We have the following subcases: 

\textbf{a)} $\eta=(uv)^{i_1} :$ Observe that by Lemma \ref{14}, the word of size $3n$ beginning in the last $uv$ of $(uv)^{i_0+i_1}$ is $(u,v)$-semi renormalizable. Note that, as before, $\omega_2=\eta v^-$. By hypothesis, we have an $(u,v)$-bifurcation in $(uv)^sv(uv)^{i_0+i_1}$. Lets suppose first that $\abs{uv}\leq n/8$. As
$$\frac{n}{2}\cdot (s-1)>(s-1)\abs{uv}=\abs{(uv)^{s+1}v}-\abs{(uv)^2v}>\abs{\alpha\beta}-3\abs{uv} \geq n-3n/8>n/2$$
we conclude that $s > 3$. Additionally, as $$\abs{(uv)^sv(uv)^{i_0+i_1}v}=\abs{(uv)^sv(uv)^{i_0}\omega_2}+2\leq 3n+2$$ 
the first item of Lemma \ref{force} let us conclude that $i_0+i_1\geq s-1>2$.

By Scheme $3$ (Lemma \ref{esquema2}), the $u$-continuation begins with $uvv$. Also, as $\abs{(uv)^{s-1}v}\\<\abs{\beta}<n$, if we consider the alphabet $(uv,(uv)^{s-2}v)$, one has in the $u$-continuation and $v$-continuation the word $uv(uv)^{s-2}v=(uv)^{s-1}v$ and then by Lemma \ref{14} the word of size $3n$ beginning in that $(uv)^{s-1}v$ is $(uv,(uv)^{s-2}v)$-semirenormalizable and then in both continuations it appears the word $(uv)^{s-2}$ after $(uv)^{s-1}v$. That is, the $u$-continuation begins with $uvv(uv)^{s-2}$ and the $v$-continuation begins with $v(uv)^{s-2}$. 
Then we determine in the continuations starting with $1$ and $2$, at least 
$$\abs{v(uv)^{s-2}}-\abs{v}=\abs{(uv)^{s-2}}=\abs{(uv)^{s+1}v}-\abs{(uv)^3v}>\abs{\alpha\beta}-4\abs{uv}>n-n/2=n/2$$
letters.

If $\abs{uv}>n/8$ and $i_0+i_1\geq 2$, then as before we can force in the continuations starting with $1$ and $2$ at least  
$$\abs{v(uv)^{i_0+i_1-1}}-\abs{v}=\abs{(uv)^{i_0+i_1-1}}\geq \abs{uv}>n/8$$
letters. 

If $i_0+i_1=1$, we have an $(u,v)$-bifurcation in $uvvuv$ and arguing as in Scheme $3$ (Lemma \ref{esquema2}), the $u$-continuation begins with $uvv$. For the $v$-continuation if after $v$ there is other $v$ then we would have 
$$[vv\dots]<[v\overline{uvv}]<[uvv\dots ]$$
and therefore, as before, by Lemma \ref{lemao} for all $j\leq l(k)$
\begin{eqnarray*}
\lambda(\sigma^j(\dots,\tilde{\alpha}_{-2},\tilde{\alpha}_{-1};\tilde{\alpha}_{0},\dots ,\tilde{\alpha}_{l(k)-\abs{v^-}},v\overline{uvv}))< \max f|_{\tilde{\Lambda}}+\epsilon/2.
\end{eqnarray*}
As $(uv,v)$ is an alphabet, $\lambda$ has value less than $3$ in any position of the $uvv$ in the middle of $uvvuvvuvv$. Also, as the cut $abb|abb$ is good, Corollary \ref{4} implies that we have the good cut $uvv|uvv$.

Then, for $j\geq l(k)+1$ we conclude that 
\begin{eqnarray*}
\lambda(\sigma^j(\dots,\tilde{\alpha}_{-2},\tilde{\alpha}_{-1};\tilde{\alpha}_{0},\dots ,\tilde{\alpha}_{l(k)-\abs{v^-}},v\overline{uvv}))&=&
\lambda(\sigma^j(\dots,\tilde{\alpha}_{-1};\tilde{\alpha}_{0},\dots,uvvuv,v\overline{uvv}))\\ &<& \max f|_{\tilde{\Lambda}}+\epsilon/2.
\end{eqnarray*}
Then, taking $x=\Pi^{-1}((\dots,\tilde{\alpha}_{-2},\tilde{\alpha}_{-1};\tilde{\alpha}_{0},\dots ,\tilde{\alpha}_{l(k)-\abs{v^-}},v\overline{uvv}))$ one would have
$$x\in W^u(\tilde{\Lambda})\cap W^s(\psi_{uvv})\subset W^u(\tilde{\Lambda})\cap W^s(\bar{\Lambda})\ \mbox{and}\ m_{\varphi,f}(x)\leq \max f|_{\tilde{\Lambda}}+\epsilon/2$$
that is a contradiction. We conclude that the $v$-continuation begins with $vu$. 

Now, the word of size $3n$ beginning with the last $uvv$ of $uvvuvvu$ is $(uv,v)$-semirenormalizable. Let $\tilde{\gamma}=uvvuX$ and $\tilde{\omega_2}$ as in the definition or renormalizartion. Again, as $\abs{\tilde{\omega_2}v}<\abs{uvv}<n$, we conclude that 
$$3n=\abs{\tilde{\gamma}\tilde{\omega}_2}=\abs{\tilde{\gamma}}+\abs{\tilde{\omega}_2}<\abs{uvv}+\abs{\tilde{\omega_2}v}+\abs{X}<2n+\abs{X},$$
which implies that $\abs{uv}<n<\abs{X}$ and then we have enough space for other letter. As $uv$ and $v$ both begin with $v^-$, so the $v$-continuation begins with $vuv^-$ and as before, we determine in the continuation beginning with $1$ and $2$, at least 
$$\abs{uvv^-}-\abs{v^-}=\abs{vuv^-}-\abs{v^-}=\abs{uv}>n/8$$ 
letters.

Suppose $i_0+i_1=0$. Note that this case is only possible if $s=1$ because in other case we would have $(uv)^2v^2$ and $\abs{(uv)^2v^2}=\abs{\beta^2}<2n$. Observe that, in both continuations, the word of size $3n$ beginning in the last $uvv$ is $(uv,v)$-semirenormalizable and in both continuations, if $\tilde{\gamma}=uvvX$ and $\tilde{\omega_2}$ are as in the definition or renormalization, again, we can conclude that $\abs{uvv}<n<\abs{X}$. Then, the $u$-continuation begins with $uvv^-$ because the first letter of $X$ in this continuation must be $uv$ and as $X\neq uv$ the following letter begins with $v^-$ in any case. For the $v$-continuation, the argument is like the previous one: if after $v$ there is other $v$ then we would have 
$$[vv\dots]<[v\overline{uvvv}]<[uvv^-\dots ]$$
and therefore, by Lemma \ref{lemao} for all $j\leq l(k)$
\begin{eqnarray*}
\lambda(\sigma^j(\dots,\tilde{\alpha}_{-2},\tilde{\alpha}_{-1};\tilde{\alpha}_{0},\dots ,\tilde{\alpha}_{l(k)-\abs{v^-}},v\overline{uvvv}))< \max f|_{\tilde{\Lambda}}+\epsilon/2.
\end{eqnarray*}
As $(uvv,v)$ is an alphabet, $\lambda$ has value less than $3$ in any position of the $uvvv$ in the middle of $uvvvuvvvuvvv$. Also, as the cut $abbb|abbb$ is good, Corollary \ref{4} implies that we have the good cut $uvvv|uvvv$.

Then, for $j\geq l(k)+1$ we conclude that 
\begin{eqnarray*}
\lambda(\sigma^j(\dots,\tilde{\alpha}_{-2},\tilde{\alpha}_{-1};\tilde{\alpha}_{0},\dots ,\tilde{\alpha}_{l(k)-\abs{v^-}},v\overline{uvvv}))&=&
\lambda(\sigma^j(\dots,\tilde{\alpha}_{-1};\tilde{\alpha}_{0},\dots,uvv,v\overline{uvvv}))\\&<& \max f|_{\tilde{\Lambda}}+\epsilon/2.
\end{eqnarray*}
Then taking $x=\Pi^{-1}((\dots,\tilde{\alpha}_{-2},\tilde{\alpha}_{-1};\tilde{\alpha}_{0},\dots ,\tilde{\alpha}_{l(k)-\abs{v^-}},v\overline{uvvv}))$ one would have
$$x\in W^u(\tilde{\Lambda})\cap W^s(\psi_{uvvv})\subset W^u(\tilde{\Lambda})\cap W^s(\bar{\Lambda})\ \mbox{and}\ m_{\varphi,f}(x)\leq \max f|_{\tilde{\Lambda}}+\epsilon/2$$
that is a contradiction. We conclude that the $v$-continuation begins with $vu$ and as $X\neq uv$ the following letter begins with $v^-$ in any case. Then, we determine in the continuations beginning with $1$ and $2$, at least 
$$\abs{uvv^-}-\abs{v^-}=\abs{vuv^-}-\abs{v^-}=\abs{uv}>n/8$$ 
letters. 
Actually, in this case, as $3\abs{uv}>\abs{(uv)^2v}=\abs{\alpha\beta}\geq n$ one determine at least $\abs{uvv^-}-\abs{v^-}=\abs{uv}>n/3$ letters in both continuations.

\textbf{b)} $\eta=(uv)^{i_1}u :$ By hypothesis, we have an $(u,v)$-bifurcation in $(uv)^sv(uv)^{i_0+i_1}u$. Observe that the word, in the sequence determined by the $u$-continuation, of size $3n$ beginning in the last $uv$ of $(uv)^sv(uv)^{i_0+i_1}u$ is $(u,v)$-semi renormalizable. Let $\tilde{\gamma}=uvuuX$ and $\tilde{\omega_2}$ as in the definition or renormalization. If $\abs{u}\leq \abs{v}$, one has $\abs{\tilde{\omega}_2}<\abs{v}$ and then 
\begin{eqnarray*}
3n=\abs{\tilde{\gamma}\tilde{\omega}_2}=\abs{\tilde{\gamma}}+\abs{\tilde{\omega}_2}&<&\abs{uvu}+\abs{\tilde{\omega}_2}+\abs{uX}<\abs{uvu}+\abs{v}+ \abs{uX}\\&=&2\abs{uv}+\abs{uX}<2n+\abs{uX},
\end{eqnarray*}
but this implies that $\abs{uv}<n<\abs{uX}$ and by Lemma \ref{comienzo} and Corollary \ref{9}, $uX$ begins with $v^-a$. If $\abs{v}<\abs{u}$ it is also true that $uX$ begins with $v^-a$ because in that situation $u=\tilde{u}v^k=\tilde{u}vv^{k-1}$ which begins with $v^-a$ where $(\tilde{u},v)$ is an ordered alphabet and $k\geq 1$. Then one has in the $u$-continuation the sequence $bu^+v(uv)^{i_0+i_1}uv^-a\subseteq (uv)^sv(uv)^{i_0+i_1}uuX$, but this contradicts Lemma \ref{alfacuadrado} because by Lemma \ref{transpuesto} 
$$bu^+v(uv)^{i_0+i_1}uv^-a=(u(v(uv)^{i_0+i_1}u)^*v)^*=(uu(vu)^{i_0+i_1}vv)^*,$$
where the internal transposed is taken respect to the alphabet $\{u,v\}$ and 
\begin{eqnarray*}
\abs{bu^+v(uv)^{i_0+i_1}uv^-a}&=&\abs{uv(uv)^{i_0+i_1}uv}=\abs{uvv(uv)^{i_0+i_1}u}\\ &\leq& \abs{(uv)^sv(uv)^{i_0}}+\abs{(uv)^{i_1}u}\leq 3n.
\end{eqnarray*}

Remember that we assumed that $\abs{\alpha}> \abs{\tilde{\beta}}$, where $\beta=\alpha^s\tilde{\beta}$. Now, if $\abs{\alpha}= \abs{\tilde{\beta}}$, then $(\alpha,\tilde{\beta})=(a,b)$ and therefore $\beta=a^sb$. As before, consider the last appearance of $\beta$ in $\gamma$ and let $\omega_2$ as in the definition of $(\alpha,\beta)$ weakly renormalizable, then the end of $\tilde{\kappa}$ is, unless one letter at the end of $\omega_2$, $\beta\alpha^{i_0}\omega_2=a^sba^{i_0}\omega_2$ where $\abs{\omega_2}<\abs{\beta}=\abs{a^sb}$ and $\omega_2$ is a prefix of $a^{s+1}$. Then one has a $(2, b)$-bifurcation in $b(2)^i$ where $i\in \mathbb{N}$ and by Scheme $2$ (Lemma \ref{muchos2}) the $2$-continuation begins with $ab$. Observe that in any sequence of $\Sigma(3)$, after the sequence $a^rb$ always follows $a^{r-1}$ because in other case we would have the bad cut $aa^{j}a|ba^{j}b$ where $j<r-1$. Additionally, note that 
$$\abs{a^s}\geq \abs{a^{s+1}b}-4\geq \abs{\alpha\beta}-4\geq n-4.$$
As by Theorem \ref{palabras}, $\Sigma(3+6^{-3n},n)=\Sigma(3,n)$, if we consider the word of size $n$ starting in the first $a$ of the smallest power $a^r$ of $a$ at the end of $a^s$ such that $\abs{a^r}\geq n/3+2$, we conclude that $(2)^i$ contains $a^{r-1}$ and $\abs{a^{r-1}} \geq  \abs{a^{r}}-2 \geq n/3$. By repeating the same argument with the first $a$ of the last $a^{r-1}$ of $(2)^i$, one concludes that we determine in both continuations at least $\abs{a^{r-2}}+2=\abs{a^{r-1}}\geq n/3$ letters.


\subsection{$\kappa=\alpha$ and $\kappa$  does not appear in $\gamma$}

As before, in this case we have $\alpha=\tilde{\alpha}\beta^s$ where $\abs{\beta}\geq \abs{\tilde{\alpha}}$ and $s\geq 1$. Suppose first that $\abs{\beta}> \abs{\tilde{\alpha}}$, then for some $k\geq 1$, $\beta=\tilde{\alpha}^k\tilde{\beta}$ with $(\tilde{\alpha},\tilde{\beta})$ alphabet. Then $\alpha=u(uv)^s$ where $u=\tilde{\alpha}$, $v=\tilde{\alpha}^{k-1}\tilde{\beta}$ and then $\beta=uv$. In this case $\gamma=\beta^{\tilde{i}}=(uv)^{\tilde{i}}$ for some $\tilde{i}\geq 2$ and let $\omega_2$ as in the definition of $(\alpha,\beta)$ weakly renormalizable, then the end of $\tilde{\kappa}$ is $\beta^{\tilde{i}}\omega_2=(uv)^{\tilde{i}}\omega_2$ where $\abs{\omega_2}<\abs{\alpha}=\abs{u(uv)^s}$ and $\omega_2$ is a prefix of $\alpha$. Observe that by Lemma \ref{14}, the word of size $3n$ beginning in the last $uv$ of $(uv)^{\tilde{i}}\omega_2$ is $(u,v)$-semi renormalizable. Consider the biggest subword $\eta$ of $u(uv)^s$ in the alphabet $\{ u,v\}$ contained in $\omega_2$. Let $(uv)^i$ where $i\geq \tilde{i}$, the biggest power of $uv$ that comes before $\omega_2$. We have the following subcases: 

\textbf{a)} $\eta=\emptyset :$ One has an $(u,v)$-bifurcation in $u(uv)^i$ or in $vv(uv)^i$ and, in any case, by schemes $3$ and $6$ (Lemmas \ref{esquema2} and \ref{esquema4}), the $u$-continuation begins with $uvv$. Suppose first that $\abs{uv}\leq n/4$ and let $i^*\geq 1$ be such that $\abs{(uv)^{i^*}v}<n$ but $\abs{(uv)^{i*+1}v}\geq n$. As
$$\frac{n}{4}\cdot i^*\geq i^*\abs{uv}=\abs{(uv)^{i^*}}=\abs{(uv)^{i^*+1}v}-\abs{uvv}> n-2\abs{uv} \geq n-n/2=n/2$$
we conclude that $i^* > 3$. 

If we consider the alphabet $(uv,(uv)^{i^*-1}v)$, one has in the $u$-continuation and $v$-continuation the word $uv(uv)^{i^*-1}v=(uv)^{i^*}v$ and then by Lemma \ref{14} the word of size $3n$ beginning in that $(uv)^{i^*}v$ is $(uv,(uv)^{i^*-1}v)$-semi renormalizable and then in both continuations it appears the word $(uv)^{i^*-1}$ after $(uv)^{i^*}v$. That is, the $u$-continuation begins with $uvv(uv)^{i^*-1}$ and the $v$-continuation begins with $v(uv)^{i^*-1}$. 
Then, we determine in the continuations starting with $1$ and $2$, at least 
$$\abs{(uv)^{i^*-1}v}-\abs{v}=\abs{(uv)^{i^*-1}}=\abs{(uv)^{i^*+1}v}-\abs{(uv)^2v}>n-3\abs{uv}>n-3n/4=n/4$$
letters.

If $\abs{uv}>n/4$, the same argument of the previous paragraph let us show that the $v$-continuation begins with $vuv$ and then, we can force in the continuations starting with $1$ and $2$ at least  
$$\abs{v(uv)}-\abs{v}=\abs{uv}>n/4$$
letters.

\textbf{b)} $\eta=u :$ Let us suppose first that $\abs{uv}\leq n/5$. In this case we have an $(u,v)$-bifurcation in $u(uv)^iu$ or $vv(uv)^iu$. By schemes $1$ and $5$ (lemmas \ref{esquema1} and \ref{esquema3}) the $u$-continuation begins with $(uv)^2$ and the $v$-continuation with $vuuvuv$. 
Let us consider the word of size $3n$ beginning in the $uuvuv$ of the $u$-continuation, this word is $(uuv,uv)$-semi renormalizable and hence, either the renormalization kernel $\tilde{\gamma}$ is of the form $u(uv)^{\tilde{s}}$ with
$$\abs{(uv)^{\tilde{s}}}=3n-\abs{\tilde{\omega}_2}-\abs{u}>3n-\abs{uuv}-\abs{u}>2n$$
or if $(uv)^{\tilde{s}}$ is the biggest power of $uv$ in the $u$-continuation and $\tilde{s}\geq i-1$ then
$$\abs{(uv)^{\tilde{s}}}\geq\abs{(uv)^{i-1}}=\abs{(uv)^i}-\abs{uv}> n-n/5=4n/5,$$
or if $\tilde{s}\leq i-2$ then one has the subword $(uv)^iu(uv)^{\tilde{s}}uuv$ and then by Lemma \ref{s}, $\abs{u(uv)^{\tilde{s}+1}}\geq 3(n-1)$ which let us conclude that 
$$\abs{(uv)^{\tilde{s}}}=\abs{u(uv)^{\tilde{s}+1}}-\abs{uuv}\geq 3(n-1)-\abs{uv}>2n$$
Then we determine in the continuation starting with $2$, at least  
$$\abs{(uv)^{\tilde{s}}}-\abs{v}>4n/5-n/5=3n/5$$
letters in any case.

In the same way, as the $v$-continuation begins with $vuuvuv$ one can consider the word os size $3n$ starting in the last $uuvuv$, this word is $(uuv,uv)$-semi renormalizable and we can consider similar cases: if the renormalization kernel $\tilde{\gamma}$ is of the form $u(uv)^{\tilde{s}}$ then 
$$\abs{u(uv)^{\tilde{s}}}=3n-\abs{\tilde{\omega}_2}>3n-\abs{uuv}>2n$$
if $(uv)^{\tilde{s}}$ is the biggest power of $uv$ that follows the first $u$ in the $v$-continuation and $\tilde{s}\geq i$ then
$$\abs{u(uv)^{\tilde{s}}}\geq\abs{u(uv)^i}> n,$$
if $\tilde{s}\leq i-1= i+1-2$ then one has the subword $(uv)^{i+1}u(uv)^{\tilde{s}}uuv$ and then by Lemma \ref{s}, $\abs{u(uv)^{\tilde{s}+1}}\geq 3(n-1)$ which let us conclude that 
$$\abs{u(uv)^{\tilde{s}}}=\abs{u(uv)^{\tilde{s}+1}}-\abs{uv}\geq 3(n-1)-\abs{uv}>2n$$
Then we determine in the continuation starting with $1$, at least  
$$\abs{vu(uv)^{\tilde{s}}}-\abs{v}=\abs{u(uv)^{\tilde{s}}}>n$$
letters in any case.

Finally, if $n/5<\abs{uv}<n$, schemes $1$ and $5$ (lemmas \ref{esquema1} and \ref{esquema3}) let us conclude again that the $u$-continuation begins with $uvu$ and the $v$-continuation with $vuuvu$. Using the alphabet $(u,v)$ we see that the word of size $3n$ beginning in the $uv$ of $uvu$ is $(u,v)$-semi renormalizable and if the renormalization kernel is $\tilde{\gamma}=uvuX$ then $\abs{uv}<n<\abs{uX}$ and by lemmas \ref{comienzo} and \ref{9}, $uX$ begins with $v^-a$. Then we determine in the continuation starting with $2$, at least  
$$\abs{uvv^-a}-\abs{v}=\abs{uv}\geq n/5$$
letters. And in the continuation starting with $1$, at least  
$$\abs{vuuvu}-\abs{v}>\abs{uv}> n/5$$
letters.

\textbf{c)} $\eta=u(uv)^{i_1}$ where $i_1>0:$ This case is not possible as in the case \textbf{c)} of subsection 5.2.

\textbf{d)} $\eta=u^2 :$ This case is not possible as in the case \textbf{d)} of subsection 5.2 because $i\geq 2$.

\textbf{e)} $\eta=u(uv)^{i_1}u :$ In this case, we have an $(u,v)$-bifurcation in $(uv)^iu(uv)^{i_1}u$. Let us suppose first that $\abs{uv}\leq n/5$, then, as $\frac{n}{5}\cdot i\geq i\abs{uv}> n$ we conclude that $i > 5$. This implies that $i_1\geq 2$, otherwise we would have in the $u$-continuation the sequence $(uv)^3uuvuu$ and then the word of size $3n$ starting in the last $uv$ of $(uv)^3uuvuu$ is $(u,v)$-semi renormalizable. Let $\tilde{\gamma}=uvuuX$ and $\tilde{\omega_2}$ as in the definition or renormalizartion. As $\abs{\tilde{\omega_2}}<\abs{uv}$, we conclude that
$$3n=\abs{\tilde{\gamma}\tilde{\omega}_2}=\abs{\tilde{\gamma}}+\abs{\tilde{\omega}_2}<\abs{uvu}+\abs{\tilde{\omega}_2}+\abs{uX}<2n+\abs{uX}.$$
But this implies that $\abs{uv}<n<\abs{uX}$ and by lemmas \ref{comienzo} and \ref{9}, $uX$ begins with $v^-a$, so one has the sequence $u^+buvuvuuvuv^-a\subseteq (uv)^3uuvuuX$ that contradicts Lemma \ref{bwords}. Then one can apply Scheme $1$ (Lemma \ref{esquema1}) to conclude that the $u$-continuation begins with $uv$ and then one has in the $u$-continuation the sequence $(uv)^iu(uv)^{i_1}uuv$. If $i_1\leq i-2$ we would have by Lemma \ref{s}
$$ 3(n-1)\leq \abs{u(uv)^{i_1+1}}=\abs{u(uv)^{i_1}v^-abu^+}=\abs{u(uv)^{i_1}uv^-}+2=\abs{\omega_2}+2<\abs{\alpha}+2<n+2$$
which is absurd. Then $i_1\geq i-1$ (observe that if $i$ is big, this case could not happen).  

Then, by Lemma \ref{esquema1}, the $u$-continuation begins with $(uv)^{i_1}$ and the $v$-continuation with $vuuvuv$.  Then, we determine in the continuation beginning with $2$, at least  
$$\abs{(uv)^{i_1}}-\abs{v}\geq \abs{(uv)^{i-1}}-\abs{v}=\abs{(uv)^i}-\abs{uuv}\geq n-2\abs{uv}\geq n-2n/5=3n/5$$
letters. 

Again, as the $v$-continuation begins with $vuuvuv$, one can consider the word of size $3n$ starting in the last $uuvuv$, this word is $(uuv,uv)$-semi renormalizable and we have some options: either the renormalization kernel $\tilde{\gamma}$ is of the form $u(uv)^{\tilde{s}}$ and then 
$$\abs{u(uv)^{\tilde{s}}}=3n-\abs{\tilde{\omega}_2}>3n-\abs{uuv}>2n$$
or if $(uv)^{\tilde{s}}$ is the biggest power of $uv$ that follows the first $u$ in the $v$-continuation and $\tilde{s}\geq i_1$ then
$$\abs{u(uv)^{\tilde{s}}}\geq\abs{u(uv)^{i_1}}\geq\abs{u(uv)^{i-1}}=\abs{u(uv)^i}-\abs{uv}>n-n/5=4n/5,$$
The case $\tilde{s}\leq i_1-1= i_1+1-2$ is not possible by Lemma \ref{s} because one would have the subword $(uv)^{i_1+1}u(uv)^{\tilde{s}}uuv$ and $\abs{u(uv)^{i_1+1}}<\abs{\omega_2}+\abs{v}<2n$. Then we determine in the continuation starting with $1$, at least  
$$\abs{vu(uv)^{\tilde{s}}}-\abs{v}=\abs{u(uv)^{\tilde{s}}}>4n/5$$
letters in any case.

Finally, if $\abs{uv}>n/5$, the same arguments of case $\textbf{b)}$ of subsection 5.2 let us determine at least $n/5$ letters in both continuations.

The case where $i_1=0$ was already considered in the case $\textbf{d)}$.

Remember that we assumed that $\abs{\beta}> \abs{\tilde{\alpha}}$, where $\alpha=\tilde{\alpha}\beta^s$. Now, if $\abs{\beta}= \abs{\tilde{\alpha}}$, then $(\tilde{\alpha},\beta)=(a,b)$ and therefore $\alpha=ab^s$. In this case $\gamma=b^{\tilde{i}}$ for some $\tilde{i}\geq 2$ and if $\omega_2$ is as in the definition of $(\alpha,\beta)$ weakly renormalizable, then the end of $\tilde{\kappa}$ is, unless one letter at the end of $\omega_2$, $b^{\tilde{i}}\omega_2$ where $\abs{\omega_2}<\abs{\alpha}=\abs{ab^s}$ and $\omega_2$ is a prefix of $\alpha$. Consider the biggest subword $\eta$ of $ab^s$ in the alphabet $\{ a,b\}$ contained in $\omega_2$. Let $(1)^i$ where $i\geq 2\tilde{i}$, the biggest power of $1$ that comes before $\omega_2$. As $\abs{(1)^i}\geq \abs{(1)^{\tilde{i}}}\geq n,$ the same arguments of cases \textbf{f)}, \textbf{g)} and \textbf{h)} of section 5.2 let us determine in both continuations at least $n/3$ letters.






\subsection{$\kappa=\beta$ and $\kappa$ does not appear in $\gamma$}

As before, in this case we have $\beta=\alpha^s\tilde{\beta}$ where $\abs{\alpha}\geq \abs{\tilde{\beta}}$ and $s\geq 1$. Suppose first that $\abs{\alpha}>\abs{\tilde{\beta}}$, then for some $k\geq1$, $\alpha=\tilde{\alpha}\tilde{\beta}^k$ with $(\tilde{\alpha},\tilde{\beta})$ an ordered alphabet. Then $\beta=(uv)^sv$ where $u=\tilde{\alpha}\tilde{\beta}^{k-1}$, $v=\tilde{\beta}$ and then $\alpha=uv$. In this case $\gamma=\alpha^{\tilde{i}}=(uv)^{\tilde{i}}$ for some $\tilde{i}\geq 2$ and let $\omega_2$ as in the definition of $(\alpha,\beta)$-weakly renormalizable, then the end of $\tilde{\kappa}$ is $\alpha^{\tilde{i}}\omega_2=(uv)^{\tilde{i}}\omega_2$ where $\abs{\omega_2}<\abs{\beta}=\abs{(uv)^sv}$ and $\omega_2$ is a prefix of $(uv)^{s+1}$. Consider the biggest subword $\eta$ of $(uv)^{s+1}$ in the alphabet $\{ u,v\}$ contained in $\omega_2$. Let $(uv)^i$ where $i\geq \tilde{i}$, the biggest power of $uv$ that comes before $\omega_2$. We have the following subcases: 

\textbf{a)} $\eta=(uv)^{i_1}$ where $i_1\geq 0 :$ In this case we have an $(u,v)$-bifurcation in $(uv)^{i+i_1}$. Then, we determine in the continuation starting with $1$ and $2$, at least $n/4$ letters. The argument is the same of case \textbf{a)} of subsection 5.3.

\textbf{b)} $\eta=(uv)^{i_1}u$ where $i_1\geq 0 :$ In this case we have an $(u,v)$-bifurcation in $(uv)^{i+i_1}u$. Then, we determine in the continuation starting with $1$ and $2$, at least $n/5$ letters. The argument is the same of case \textbf{b)} of subsection 5.3.

Now, if $\abs{\alpha}= \abs{\tilde{\beta}}$, then $(\alpha,\tilde{\beta})=(a,b)$ and therefore $\beta=a^sb$. In this case $\gamma=a^{\tilde{i}}$ for some $\tilde{i}\geq 2$ and if $\omega_2$ is as in the definition of $(\alpha,\beta)$-weakly renormalizable, then the end of $\tilde{\kappa}$ is, unless one letter at the end of $\omega_2$, $a^{\tilde{i}}\omega_2$ where $\abs{\omega_2}<\abs{\beta}=\abs{a^sb}$ and $\omega_2$ is a prefix of $a^{s+1}$. Consider the biggest subword $\eta$ of $a^{s+1}$ in the alphabet $\{ a,b\}$ contained in $\omega_2$. Let $(2)^i$ where $i\geq 2\tilde{i}$, the biggest power of $2$ that comes at the end of $\tilde{\kappa}$. As $\abs{(2)^i}\geq \abs{(2)^{\tilde{i}}}\geq n,$ if we consider the word of size $n$ starting in the first $a$ of the smallest power $a^r$ of $a$ at the end of $a^i$ such that $\abs{a^r}\geq n/3$, we conclude that we determine in both continuations at least $\abs{a^{r-1}}+2=\abs{a^r}\geq n/3$ letters.

\subsection{End of the proof of Proposition \ref{fundamental}}

Summarizing what we did until now, one has that in any case, if $(a_0, a_1, \dots,a_{l(k)})$ has two continuations, $ \gamma_{l(k)+1}=\\ (2, a_{l(k)+2},\dots)$ and $\beta_{l(k)+1}=(1, a^*_{l(k)+2},\dots)$, then $p_1=(2,a_{l(k)+2},\dots,a_{l(k)+n/5})$ and $p_2=(2,a^*_{l(k)+2},\dots,a^*_{l(k)+n/5})$ are uniquely determined, as we claimed before. In particular, we can refine the cover $\mathcal{C}_k$ by replacing the interval $I^u(a_0; a_1, \dots,a_k)$ with the two intervals $I^u(a_0;a_1, \dots, a_k,p_1)$ and  $ I^u(a_0;a_1, \dots, a_k,p_2)$. Indeed, we affirm that for some constant $c>0$ this procedure does not increase the $\frac{c}{n}$-sum, $H_{\frac{c}{n}}(\mathcal{C}_k)= \sum \limits_{I\in \mathcal{C}_k} \abs{I}^{\frac{c}{n}}$ of the cover $\mathcal{C}_k$ of $K^u(\tilde{\Lambda})$. That is, we need to prove that
$$ \abs{I^u(a_1, \dots, a_k,p_1)}^{\frac{c}{n}} + \abs{I^u(a_1, \dots, a_k,p_2)}^{\frac{c}{n}} <  \abs{I^u(a_1, \dots, a_k)}^{\frac{c}{n}} $$
or
\begin{equation}\label{sum}
 \left(\frac{\abs{I^u(a_1, \dots, a_k,p_1)}}{\abs{I^u(a_1, \dots, a_k)}}\right)^{\frac{c}{n}}+ \left(\frac{\abs{I^u(a_1, \dots, a_k,p_2)}}{\abs{I^u(a_1, \dots, a_k)}}\right)^{\frac{c}{n}}<1.
 \end{equation}

Remember that in our context of dynamically defined Cantor sets, we can relate the length of the unstable intervals determined by an admissible word to its length as a word in the alphabet $\mathcal{A}$ via the {\it bounded distortion property} that let us conclude that for some constant $c_1>0$, and admissible words $\alpha$ and $\beta$
\begin{equation}\label{bdp1}
e^{-c_1}|I^u(\alpha)|\cdot|I^u(\beta)|\le |I^u(\alpha\beta)|\le e^{c_1}|I^u(\alpha)|\cdot|I^u(\beta)|,
\end{equation}
and also that for some positive constants $\lambda_1,\lambda_2<1$, one has
\begin{equation}\label{bdp3}
e^{-c_1} \lambda_1^{\abs{\alpha}}\leq\abs{I^u(\alpha)}\leq e^{c_1} \lambda_2^{\abs{\alpha}}.
\end{equation} 
We conclude for $i=1,2$, that if $c=\frac{-5\log4}{\log \lambda_2}$ and $n$ is large
\begin{eqnarray*}
\left(\frac{\abs{I^u(a_1, \dots, a_k,p_i)}}{\abs{I^u(a_1, \dots, a_k)}}\right)^{\frac{c}{n}}&\leq& \left(\frac{ e^{c_1}|I^u(a_1,\dots,a_k )|\cdot|I^u(p_i)|}{|I^u(a_1,\dots,a_k )|}\right)^{\frac{c}{n}}\\ &\leq& ( e^{2c_1}\lambda_2^{\abs{p_i}}) ^{\frac{c}{n}}=( e^{2c_1}\lambda_2^{\frac{n}{5}}) ^{\frac{c}{n}}<1/2
\end{eqnarray*}
that proves (\ref{sum}) and so let us conclude that $HD(K^u(\tilde{\Lambda}))\leq \frac{c}{n}$ for $n$ large. Finally, as we are in the conservative setting
$$HD(\tilde{\Lambda})=2HD(K^u(\tilde{\Lambda}))\leq \frac{2c}{n}.$$
This finishes the proof of the proposition with $C_0=2c$.

\section{Proof of the main theorems}

\subsection{Proof of Theorem \ref{teo1}}

Remember that we are considering the horseshoe
$\Lambda(2)=C(2)\times\tilde{C}(2)$ equipped with the diffeomorphism $\varphi$ and the map $f$ as in Section \ref{Lambda}. Given $\epsilon>0$, take $r(\epsilon)\in \mathbb{N}$ sufficiently large such that if $\alpha=(c_{-r(\epsilon)}, \dots,c_0,\dots, c_{r(\epsilon)}) \in \{1,2\}^{2r(\epsilon)+1}$
and $x, y\in R(\alpha;0)=\Pi^{-1}\{(x_{n})\in \Sigma:(x_{-r(\epsilon)},\dots, x_0,\dots,x_{r(\epsilon)})=\alpha \}$  then $\abs{f(x)-f(y)}<\epsilon/4.$ 

Now, we will define the sequence $\{a_r \}_{r\in \mathbb{N}}$ as in the statement of the theorem: Let $n_1\in \mathbb{N}$ such that for any $n\geq n_1$ proposition \ref{fundamental} holds, $3+6^{-3n}<t_1$ and $\frac{6C_0\cdot C_2\cdot \log 6}{C_1^2}<\log (3n\cdot \log 6)$, where $C_0$ comes from Proposition \ref{fundamental} and $C_1, C_2$ are given by equations \ref{cunocdos} and  \ref{cunocdos2}. Define $a_1=3+6^{-3n_1}$ and once we have defined $a_r\in \mathbb{R}$, set $a_{r+1}=3+6^{-3n_{r+1}}$ where $n_{r+1}=\min \{n\in \mathbb{N}:d(3+6^{-3n})<d(a_r) \}$. Note that, by definition, the sequence $\{d(a_r) \}_{r\in \mathbb{N}}$ is strictly decreasing. Fix $r\in \mathbb{N}$ and consider $t\in(a_{r+1},a_r)\cap \mathcal{L}$ and $\epsilon>0$ such that $t+\epsilon<a_r$ let
$$C(t,\epsilon)=\{\alpha=(c_{-r(\epsilon)}, \dots,c_0,\dots, c_{r(\epsilon)})\in \{1,2\}^{2r(\epsilon)+1}:R(\alpha;0)\cap (\Lambda(2))_{t+\epsilon/4}\neq \emptyset \}.$$
Define
$$M(t,\epsilon):=\bigcap \limits_{n \in \mathbb{Z}} \varphi ^{-n}(\bigcup \limits_{\alpha \in C(t,\epsilon)}  R(\alpha;0)).$$
Note that by construction, $(\Lambda(2))_{t+\epsilon/4}\subset M(t,\epsilon)\subset (\Lambda(2))_{t+\epsilon/2}$ and being $M(t,\epsilon)$ a hyperbolic set of finite type (see Subsection \ref{tipofinito} for the corresponding definitions and results), it admits a decomposition 
$$M(t,\epsilon)=\bigcup \limits_{x\in \mathcal{X}(t,\epsilon)} \tilde{\Lambda}_x $$
 where $\mathcal{X}(t,\epsilon)$ is a finite index set and for $x\in \mathcal{X}(t,\epsilon)$,\ $\tilde{\Lambda}_x$ is a subhorseshoe or a transient set of the form $\mathcal{T} _{\tilde{\Lambda}_{x_1},\tilde{\Lambda}_{x_2}}$, where $\tilde{\Lambda}_{x_1}$ and $\tilde{\Lambda}_{x_2}$ with $x_1, x_2 \in \mathcal{X}(t,\epsilon)$ are subhorseshoes. As for every transient set $\mathcal{T} _{\tilde{\Lambda}_{x_1},\tilde{\Lambda}_{x_2}}$ as before, we have by Equation \ref{transient}
$$HD(\mathcal{T}_{\tilde{\Lambda}_{x_1},\tilde{\Lambda}_{x_2}})=\frac{HD(\tilde{\Lambda}_{x_1})+HD(\tilde{\Lambda}_{x_2})}{2} \leq \max\{HD(\tilde{\Lambda}_{x_1}),HD(\tilde{\Lambda}_{x_2})\},$$
we conclude that 
\begin{equation}
  HD(M(t,\epsilon))=\max\limits_{x\in \mathcal{X}(t,\epsilon)} HD(\tilde{\Lambda}_{x})=\max \limits_{\substack{x\in \mathcal{X}(t,\epsilon): \ \tilde{\Lambda}_x \ is\\ subhorseshoe }}HD(\tilde{\Lambda}_{x}).  
\end{equation}
With this in mind, let us consider
$$\tilde{M}(t,\epsilon)=\bigcup\limits_{\substack{x\in \mathcal{X}(t,\epsilon): \ \tilde{\Lambda}_x \ is\\ subhorseshoe }}\tilde{\Lambda}_x= \bigcup \limits_{i\in \mathcal{I}(t,\epsilon)} \tilde{\Lambda}_i \cup \bigcup \limits_{i\in \mathcal{J}(t,\epsilon)} \tilde{\Lambda}_j$$ where
 $$\mathcal{I}(t,\epsilon)=\{i\in \mathcal{X}(t,\epsilon): \tilde{\Lambda}_i \ \mbox{is a subhorseshoe and it connects with}\  \psi_b \ \mbox{before}\ \max f|_{\tilde{\Lambda}_i}+\frac{\epsilon}{2}\}$$
and $\mathcal{J}(t,\epsilon)=\{x\in \mathcal{X}(t,\epsilon): \ \tilde{\Lambda}_x \ \mbox{is a subhorseshoe}\}\setminus \mathcal{I}(t,\epsilon)$. Note that for any $j\in \mathcal{J}(t,\epsilon)$
$$\max f|_{\tilde{\Lambda}_j}+\epsilon/2< t+\epsilon/2+\epsilon/2<a_r= 3+6^{-3n_r}$$  
and then, by Proposition \ref{fundamental}, one has $HD(\tilde{\Lambda}_j)\leq \frac{C_0}{n_r}$.

On the other hand, by definition, for $i\in \mathcal{I}(t,\epsilon)$, $\tilde{\Lambda}_i$ connects with  $\psi_b$ before $\max f|_{\tilde{\Lambda}_i}+\epsilon/2<t+\epsilon$, then we can apply Corollary 3.4 at most $\abs{\mathcal{I}(t,\epsilon)}-1$ times to see that there exists a subhorseshoe $\tilde{\Lambda}(t,\epsilon)\subset \Lambda(2)$ and some $q(t,\epsilon)<t+\epsilon$ such that 
$$\bigcup \limits_{i\in \mathcal{I}(t,\epsilon)} \tilde{\Lambda}_i\subset \tilde{\Lambda}(t,\epsilon)\subset (\Lambda(2))_{q(t,\epsilon)}.$$

As any subhorseshoe $\tilde{\Lambda}_x$ is locally maximal we have (see page $18$ of \cite{GC})

$$\ell_{\varphi,f}(M(t,\epsilon))=\ell_{\varphi,f}(\tilde{M}(t,\epsilon))=\bigcup \limits_{i\in \mathcal{I}(t,\epsilon)} \ell_{\varphi,f}(\tilde{\Lambda}_i)\cup\bigcup \limits_{j\in \mathcal{J}(t,\epsilon)} \ell_{\varphi,f}(\tilde{\Lambda}_j).$$

We are ready to define the set of Theorem \ref{teo1}:
$$\tilde{\mathcal{B}}_r:=\{t\in (a_{r+1},a_r)\cap \mathcal{L}:\forall \epsilon>0,\ \  (t-\epsilon/4,t+\epsilon/4)\cap \bigcup \limits_{x\in \mathcal{I}(t,\epsilon)} \ell_{\varphi,f}(\tilde{\Lambda}_x)\neq \emptyset \}.$$
Observe that, given $t\in ((a_{r+1},a_r)\cap \mathcal{L})\setminus \tilde{\mathcal{B}}_r$, one can find some $\epsilon(t)$ such that
$$(t-\epsilon(t)/4,t+\epsilon(t)/4)\cap \mathcal{L}\subset  \bigcup \limits_{j\in \mathcal{J}(t,\epsilon(t))} \ell_{\varphi,f}(\tilde{\Lambda}_j).$$
Consider some countable sequence $\{t_n\}_{n\in\mathbb{N}}\subset ((a_{r+1},a_r)\cap \mathcal{L})\setminus \tilde{\mathcal{B}}_r$ such that 
$$((a_{r+1},a_r)\cap \mathcal{L})\setminus \tilde{\mathcal{B}}_r\subset \bigcup \limits_{n\in \mathbb{N}} (t_n-\epsilon(t_n)/4,t_n+\epsilon(t_n)/4)\cap \mathcal{L}$$
then, one concludes
\begin{eqnarray*}
    HD(((a_{r+1},a_r)\cap \mathcal{L})\setminus \tilde{\mathcal{B}}_r)&\leq& \sup \limits_{n\in \mathbb{N}}HD((t_n-\epsilon(t_n)/4,t_n+\epsilon(t_n)/4)\cap \mathcal{L})\\ &\leq& \sup \limits_{n\in \mathbb{N}} HD(\bigcup \limits_{j\in \mathcal{J}(t_n,\epsilon(t_n))} \ell_{\varphi,f}(\tilde{\Lambda}_j))\\ &\leq&  
    \sup \limits_{n\in \mathbb{N}} HD(\bigcup \limits_{j\in \mathcal{J}(t,\epsilon(t))} f(\tilde{\Lambda}_j))\\ &=&  \sup \limits_{n\in \mathbb{N}} \max \limits_{j\in \mathcal{J}(t_n,\epsilon(t_n))} HD(f(\tilde{\Lambda}_j))\\ &\leq& \sup \limits_{n\in \mathbb{N}} \max \limits_{j\in \mathcal{J}(t_n,\epsilon(t_n))} HD(\tilde{\Lambda}_j) \leq \frac{C_0}{n_r}.  
    \end{eqnarray*}

Given $r\in \mathbb{N}$ define 
$$I_r=\{t\in \tilde{\mathcal{B}}_r:\exists s>0 \ \text{such that}\ (t-s,t)\cap \tilde{\mathcal{B}}_r= \emptyset \}.$$
Then, $I_r$ is the enumerable set of points of $\tilde{\mathcal{B}}_r$ isolated on the left. Note that $\mathcal{B}_r:=\tilde{\mathcal{B}}_r\setminus  I_r\subset \mathcal{L}^{'}$ and $HD(((a_{r+1},a_r)\cap \mathcal{L})\setminus \tilde{\mathcal{B}}_r)=HD(((a_{r+1},a_r)\cap \mathcal{L})\setminus \mathcal{B}_r)$. 

We are ready to prove the theorem. Let us show first that for any $t\in \mathcal{B}_r$ one has $D(t)=HD(\ell^{-1}(t))$: Given $s\in \tilde{\mathcal{B}_r}$ and $\epsilon>0$ small, we can find $i_0\in \mathcal{I}(s,\epsilon)$ and $r_0\in \tilde{\Lambda}_{i_0}$ such that 
$\ell_{\varphi,f}(r_0)\in (s-\epsilon/4,s+\epsilon/4)$. Also, as 
$$\ell_{\varphi,f}(\tilde{\Lambda}_{i_0})\subset \ell_{\varphi,f}(\tilde{\Lambda}(s,\epsilon))\subset f(\tilde{\Lambda}(s,\epsilon))\subset f((\Lambda(2))_{q(s,\epsilon)})\subset (-\infty,s+\epsilon],$$
we conclude that $s-\epsilon/4<\max f|_{\tilde{\Lambda}(s,\epsilon)}\leq s+\epsilon$. On the other hand, given $r\in \mathbb{N}$, by definition of $a_r$ and Equation \ref{cunocdos}, one has
\begin{eqnarray*}
\frac{C_0}{n_r}<\frac{C_1^2}{2C_2}\cdot\frac{\log (3n_r\cdot \log 6)}{3n_r\cdot\log 6}&\leq&\frac{C_1}{2C_2}\cdot d(3+6^{-3n_r})=\frac{C_1}{2C_2}\cdot d(3+6^{-3(n_{r+1}-1)})\\ &\leq& \frac{C_1}{2}\cdot\frac{\log (3(n_{r+1}-1)\cdot \log 6)}{3(n_{r+1}-1)\cdot\log 6}\\&<& \frac{C_1\log (3n_{r+1})\cdot \log 6)}{3n_{r+1}\cdot\log 6}\\ &\leq& d(3+6^{-3n_{r+1}})=d(a_{r+1}).
\end{eqnarray*}
From this, we conclude that 
\begin{eqnarray*}
HD(\bigcup \limits_{j\in \mathcal{J}(s,\epsilon)} \tilde{\Lambda}_j)=\max\limits_{j\in \mathcal{J}(s,\epsilon)}HD(\tilde{\Lambda}_j)&\leq&\frac{C_0}{n_r}<d(a_{r+1})\leq HD((\Lambda(2))_s)\leq HD(M(s,\epsilon))\\ &=&HD(\tilde{M}(s,\epsilon))=HD(\bigcup \limits_{x\in \mathcal{X}(s,\epsilon)} \tilde{\Lambda}_x)
\end{eqnarray*}
and then 
$$HD((\Lambda(2))_s)\leq HD(\bigcup \limits_{i\in \mathcal{I}(s,\epsilon)} \tilde{\Lambda}_i)\leq HD(\tilde{\Lambda}(s,\epsilon)).$$
Finally, if some subhorseshoe $\tilde{\Lambda}\subset (\Lambda(2))_s$ satisfies that $HD(\tilde{\Lambda})>\frac{C_0}{n_r}$ then 
$$\tilde{\Lambda}\subset \bigcup \limits_{i\in \mathcal{I}(s,\epsilon)} \tilde{\Lambda}_i\subset \tilde{\Lambda}(s,\epsilon).$$
These conclusions are similar to the hypothesis of proposition 3.3 of \cite{GC}. Indeed, given some $t\in \mathcal{B}_r$ we consider a strictly increasing sequence $\{s_n\}_{n\in \mathbb{N}}$ of elements of $\tilde{\mathcal{B}_r}$ with $s_0$ arbitrarily close to $t$ and some sequence $\{\epsilon_n \}_{n\in \mathbb{N}}$ of small positive numbers. Using that the sequence of subhorseshoes $\{\tilde{\Lambda}(s_n,\epsilon_n) \}_{n\in \mathbb{N}}$ is increasing and the limits: $\lim\limits_{{n \to \infty}} \max f|_{{\Lambda}(s_n,\epsilon_n)}=t$ and $\lim\limits_{n \to \infty} HD({\Lambda}(s_n,\epsilon_n))=HD((\Lambda(2))_t)$, one is able to construct a homeomorphism $\theta:K^u({\Lambda}(s_0,\epsilon_0)) \rightarrow \ell^{-1}(t)$ whose inverse is H\"older with exponent arbitrarily close to one. Lettin first the H\"older exponent tend to $1$ and then $s_0$ tend to t, one gets easily that $D(t)=HD(\ell^{-1}(t))$. The details are in Section 3.3 of \cite{GC}. 

Now, the spectral decomposition theorem and Corollary 3.9 of \cite{GC} let us conclude the following proposition
\begin{proposition}
  Given two subhorseshoes $\tilde{\Lambda}_1$ and $\tilde{\Lambda}_2$ of $\Lambda$ such that $\tilde{\Lambda}_1\nsubseteq \tilde{\Lambda}_2$, we have
    $$HD(\tilde{\Lambda}_1)<HD(\tilde{\Lambda}_2).$$  
\end{proposition}
That we use to show that the function $D|_{\mathcal{B}_r}$ is strictly increasing: Observe first that for any $t\in \tilde{\mathcal{B}}_r$ and $\epsilon>0$ small, by definition, $t\in \ell_{\varphi,f}(\tilde{\Lambda}(t,\epsilon))\subset f(\tilde{\Lambda}(t,\epsilon))$ (this can be proved by letting $0<\tilde{\epsilon}<\epsilon$ tends to $0$ in the definition of $\tilde{\mathcal{B}}_r$ and use that $\ell_{\varphi,f}(\tilde{\Lambda}(t,\epsilon))$ is a closed set). If $t_1$ and $t_2$ are two elements of $\mathcal{B}_r
 =\tilde{\mathcal{B}}_r\setminus I_r$ such that $t_1<t_2$ then, consider any $t_3\in \tilde{\mathcal{B}}_r$ such that $t_1<t_3<t_2$ and $\epsilon>0$ small enough. Then, as $\max f|_{\tilde{\Lambda}(t_1,\epsilon)}<t_1+\epsilon$, $\tilde{\Lambda}(t_1,\epsilon) \subset \tilde{\Lambda}(t_3,\epsilon)$ and $t_3\in f(\tilde{\Lambda}(t_3,\epsilon))$ one has $\tilde{\Lambda}(t_1,\epsilon) \subsetneq \tilde{\Lambda}(t_3,\epsilon) $, and then
$$D(t_1)=\frac{1}{2}HD((\Lambda(2))_{t_1})\leq \frac{1}{2}HD(\tilde{\Lambda}(t_1,\epsilon))< \frac{1}{2}HD(\tilde{\Lambda}(t_3,\epsilon))\leq\frac{1}{2}HD((\Lambda(2))_{t_2})=D(t_2).$$

Finally, let us show that for $r\in \mathbb{N}$, one has $HD(\mathcal{B}_r)=HD((a_{r+1},a_r)\cap \mathcal{L})$: Note first that
$$d(a_r)=\max \{HD((a_{r+1},a_r)\cap \mathcal{L}),d(a_{r+1})\}=HD((a_{r+1},a_r)\cap \mathcal{L}).$$
On the other hand, one has
$$HD(((a_{r+1},a_r)\cap \mathcal{L})\setminus \tilde{\mathcal{B}}_r)\leq \frac{C_0}{n_r}<d(a_{r+1})<d(a_r)$$
and then
$$HD((a_{r+1},a_r)\cap \mathcal{L})=\max \{HD(\tilde{\mathcal{B}}_r),HD(((a_{r+1},a_r)\cap \mathcal{L})\setminus \tilde{\mathcal{B}}_r
 )\}=HD(\tilde{\mathcal{B}}_r).$$
As $I_r$ is enumerable, we have the result.

\subsection{Proof of Theorem \ref{teo2}}

Let $n_1\in \mathbb{N}$ such that for any $n\geq n_1$, Proposition \ref{fundamental} holds and let $\rho>0$ small such that $6^{-{3(n+1)}}\leq\rho<6^{-3n}$ for some $n\geq n_1$ and inequalities \ref{cunocdos} and \ref{cunocdos2} hold. Given $t\in (\mathcal{M}\setminus \mathcal{L})\cap (-\infty,3+\rho)$ let $\tilde{x}\in \Lambda(2)$ be such that $f(\tilde{x})=m_{\varphi,f}(\tilde{x})=t$. For small $\epsilon>0$, as before, consider the hyperbolic set of finite type $M(t,\epsilon)$ such that 
$(\Lambda(2))_{t+\epsilon/4}\subset M(t,\epsilon)\subset (\Lambda(2))_{t+\epsilon/2}$ and write it as
$$M(t,\epsilon)=\bigcup \limits_{x\in \mathcal{X}(t,\epsilon)} \tilde{\Lambda}^{\epsilon}_x $$
 where for $x\in \mathcal{X}(t,\epsilon)$,\ $\tilde{\Lambda}^{\epsilon}_x$ is a subhorseshoe or a transient set of the form $\mathcal{T}_{\tilde{\Lambda}^{\epsilon}_{x_1}, \tilde{\Lambda}^{\epsilon}_{x_2}}$, where $\tilde{\Lambda}^{\epsilon}_{x_1}$ and $\tilde{\Lambda}^{\epsilon}_{x_2}$ with $x_1, x_2 \in \mathcal{X}(t,\epsilon)$ are subhorseshoes (here, for convenience, we specified the dependence on $\epsilon$). 

 Note that there exists $\epsilon_1>0$ such that for $0<\epsilon<\epsilon_1$, one has $\tilde{x}\notin \tilde{\Lambda}^{\epsilon}_x$ if $\tilde{\Lambda}^{\epsilon}_x$ is a subhorseshoe, otherwise we would have a sequence of positive numbers, $\{\epsilon_n\}_{n\in \mathbb{N}}$ such that $\lim \limits_{n\to \infty}\epsilon_n=0$, $\tilde{\Lambda}^{\epsilon_n}_{x_n}$ is subhorseshoe and $\tilde{x}\in \tilde{\Lambda}^{\epsilon_n}_{x_n}$ for some $x_n\in \mathcal{X}(t,\epsilon_n)$. As $t\leq \max f|_{ \tilde{\Lambda}^{\epsilon_n}_{x_n}}\leq t+\epsilon_n/2$ we would have $\lim \limits_{n\to \infty}\max f|_{ \tilde{\Lambda}^{\epsilon_n}_{x_n}}=t$. But maximums of subhorseshoes are always elements of $\mathcal{L}$ which is a closed set, then we would get the contradiction $t\in \mathcal{L}$.

 Then, one has for $0<\epsilon<\epsilon_1$ that $\tilde{x}\in  \mathcal{T}_{\tilde{\Lambda}^{\epsilon}_{x_1}, \tilde{\Lambda}^{\epsilon}_{x_2}} $ where $\tilde{\Lambda}^{\epsilon}_{x_1}$ and $\tilde{\Lambda}^{\epsilon}_{x_2}$ with $x_1, x_2 \in \mathcal{X}(t,\epsilon)$ are subhorseshoes. We affirm that we can find some $0<\epsilon_2<\epsilon_1$ such that for $0<\epsilon<\epsilon_2$ either  $x_1\in \mathcal{J}(t,\epsilon)$ or $x_2\in \mathcal{J}(t,\epsilon)$ where, as before, $\mathcal{J}(t,\epsilon)$ is the set of index $j\in \mathcal{X}(t,\epsilon)$ such that $\tilde{\Lambda}_j$ is a subhorseshoe that does not connect with $\psi_b$ before $\max f|_{\tilde{\Lambda}_j}+\epsilon/2$. Otherwise, we would have a sequence of positive numbers, $\{\epsilon_n\}_{n\in \mathbb{N}}$ such that $\lim \limits_{n\to \infty}\epsilon_n=0$ and the subhorsehoes $\tilde{\Lambda}^{\epsilon_n}_{x_1}$ and $\tilde{\Lambda}^{\epsilon_n}_{x_2}$ connect with $\psi_b$ before $\max \{\max f|_{\tilde{\Lambda}_{x_1}}, \max f|_{\tilde{\Lambda}_{x_2}}\}+\epsilon_n/2\leq t+\epsilon_n$ and in particular, we can find  some $\tilde{y}\in W^u(\tilde{\Lambda}^{\epsilon_n}_{x_2})\cap W^s(\tilde{\Lambda}^{\epsilon_n}_{x_1})$ with $m_{\varphi,f}(\tilde{y})<t+\epsilon_n$ and as $\tilde{x}\in W^u(\tilde{\Lambda}^{\epsilon_n}_{x_1})\cap W^s(\tilde{\Lambda}^{\epsilon_n}_{x_2})$, Proposition \ref{connection11} lets us find some subhorseshoe $\tilde{\Lambda}_n$ such that $\tilde{\Lambda}^{\epsilon_n}_{x_1}\cup \tilde{\Lambda}^{\epsilon_n}_{x_2}\cup \mathcal{O}(\tilde{x})\subset \tilde{\Lambda}_n\subset (\Lambda(2))_{t+\epsilon_n}$ which allows us to get, as before, the contradiction: $\lim \limits_{n\to \infty}\max f|_{ \tilde{\Lambda}_n}=t$. 

If $\mathcal{C}$ is the collection of pairs of subhorseshoes $(\Lambda^1,\Lambda^2)$ of $\Lambda(2)$ such that for some $t$ and $\epsilon>0$ that satisfy $t+\epsilon<3+\rho<3+6^{-3n}$ one has that $\mathcal{T}_{\Lambda^1, \Lambda^2}$ is a transient component of $M(t,\epsilon)$ and either $\Lambda^1$  does not connect with $\psi_b$ before $\max f|_{\Lambda^1}+\epsilon/2$ or $\Lambda^2$ does not connect with $\psi_b$ before $\max f|_{\Lambda^2}+\epsilon/2$.  Then, one conclude that
\begin{eqnarray}\label{mmenosl}
(\mathcal{M}\setminus \mathcal{L})\cap (-\infty,3+\rho)\subset \bigcup \limits_{(\Lambda^1,\Lambda^2)\in \mathcal{C}}f(\mathcal{T}_{\Lambda^1, \Lambda^2}).
\end{eqnarray}
As for every transient set $\mathcal{T}_{\Lambda^1, \Lambda^2}$ one has $$HD(\mathcal{T}_{\Lambda^1, \Lambda^2})=\frac{HD(\Lambda^1)+HD(\Lambda^2)}{2}$$
we conclude for $(\Lambda^1,\Lambda^2) \in \mathcal{C}$ that 
\begin{eqnarray*}
HD(\mathcal{T}_{\Lambda^1, \Lambda^2})&\leq& \frac{1}{2}HD((\Lambda(2))_{3+\rho})+\frac{C_0}{2n}=\frac{1}{2}d(3+\rho)+\frac{C_0}{2n}\\ &\leq& \frac{\log(\abs{\log \rho})-\log (\log(\abs{\log \rho}))+C_2}{\abs{\log \rho}}+\frac{C_0}{n+1}\\ &\leq& \frac{\log(\abs{\log \rho})-\log (\log(\abs{\log \rho}))+C_2}{\abs{\log \rho}}+\frac{3\log 6\cdot C_0}{\abs{\log \rho}}\\ &=&\frac{\log(\abs{\log \rho})-\log (\log(\abs{\log \rho}))+C_2+3\log 6\cdot C_0}{\abs{\log \rho}},
\end{eqnarray*}
because $\Lambda^1\cup \Lambda^2\subset (\Lambda(2))_{3+\rho}$ and for some $\epsilon>0$ either $\Lambda^1$  does not connect with $\psi_b$ before $\max f|_{\Lambda^1}+\epsilon/2$ or $\Lambda^2$ does not connect with $\psi_b$ before $\max f|_{\Lambda^2}+\epsilon/2$ (here we used Proposition \ref{fundamental}). By \ref{mmenosl} we get 
\begin{eqnarray*}
HD((\mathcal{M}\setminus \mathcal{L})\cap (-\infty,3+\rho))&\leq& \sup \limits_{(\Lambda^1,\Lambda^2)\in \mathcal{C}}HD(f(\mathcal{T}_{\Lambda^1, \Lambda^2}))\leq \sup \limits_{(\Lambda^1,\Lambda^2)\in \mathcal{C}}HD(\mathcal{T}_{\Lambda^1, \Lambda^2})\\ &\leq&  \frac{\log(\abs{\log \rho})-\log (\log(\abs{\log \rho}))+C}{\abs{\log \rho}},
\end{eqnarray*}
where $C=C_2+3\log 6\cdot C_0$. This finishes the proof of the theorem.

\end{document}